\documentclass[a4paper,11pt,reqno,noindent]{amsart}
\usepackage[centertags]{amsmath}
\usepackage{amsfonts,amssymb,amsthm,dsfont,cases,amscd,esint,enumerate}
\usepackage[T1]{fontenc}
\usepackage[english]{babel}
\usepackage[applemac]{inputenc}
\usepackage{newlfont}
\usepackage{color}
\usepackage[body={15cm,21.5cm},centering]{geometry} 
\usepackage{fancyhdr}
\usepackage{enumerate}
\usepackage[colorlinks,citecolor=black,urlcolor=black]{hyperref}

\theoremstyle{plain}
\newtheorem{theor0}{Theorem}
\newenvironment{theor}
  {\pushQED{\qed}\begin{theor0}}
  {\popQED\end{theor0}}
\newtheorem{lem0}{Lemma}[section]
\newenvironment{lem}
  {\pushQED{\qed}\begin{lem0}}
  {\popQED\end{lem0}}
\newtheorem{prop0}[lem0]{Proposition}
\newenvironment{prop}
  {\pushQED{\qed}\begin{prop0}}
  {\popQED\end{prop0}}
\newtheorem{cor0}[lem0]{Corollary}
\newenvironment{cor}
  {\pushQED{\qed}\begin{cor0}}
  {\popQED\end{cor0}}
\theoremstyle{definition}
\newtheorem{defin0}[lem0]{Definition}

\newtheorem{rems0}[lem0]{Remarks}
\newenvironment{rems}
  {\pushQED{\qed}\begin{rems0}}
  {\popQED\end{rems0}}
\newtheorem{rem0}[lem0]{Remark}
\newenvironment{rem}
  {\pushQED{\qed}\begin{rem0}}
  {\popQED\end{rem0}}

\numberwithin{equation}{section}
\mathchardef\emptyset="001F

\newcommand{\e}{\varepsilon}
\newcommand{\Log}{|\!\log\e|}

\newcommand{\calQ}{\mathcal{Q}}

\newcommand{\calP}{\mathcal{P}}
\newcommand{\calN}{\mathcal{N}}
\newcommand{\calS}{\mathcal S}
\newcommand{\calK}{\mathcal K}
\newcommand{\calT}{\mathcal T}

\newcommand{\R}{\mathbb R}
\newcommand{\Z}{\mathbb Z}

\newcommand{\B}{\mathcal B}
\newcommand{\Nc}{\mathcal N}
\newcommand{\G}{G}
\newcommand{\F}{F}
\newcommand{\Id}{\operatorname{Id}}
\newcommand{\E}{\mathbb{E}}

\newcommand{\ee}{e}
\newcommand{\Aa}{\boldsymbol a}
\newcommand{\dig}[1]{\mathrm{diag}\left[ #1\right]}
\newcommand{\Ld}{\operatorname{L}}
\newcommand{\step}[1]{\noindent \textit{Step} #1.}

\newcommand{\Pm}{\mathbb{P}}
\newcommand{\pr}[1]{\mathbb{P}\left[ #1 \right]}
\newcommand{\expec}[1]{\mathbb{E}\left[ #1 \right]}
\newcommand{\expecm}[1]{\mathbb{E}\big[ #1 \big]}
\newcommand{\var}[1]{\mathrm{Var}\left[#1\right]}
\newcommand{\varm}[1]{\mathrm{Var}\big[#1\big]}
\newcommand{\cov}[2]{\operatorname{Cov}\left[{#1};{#2}\right]}

\newcommand{\expeCm}[2]{\mathbb{E}\big[ #1 \big|\,#2\big]}
\newcommand{\dW}[2]{\operatorname{d}_{\operatorname{W}}\left({#1},{#2}\right)}
\newcommand{\dK}[2]{\operatorname{d}_{\operatorname{K}}\left({#1},{#2}\right)}
\newcommand{\dWK}[2]{(\operatorname{d}_{\operatorname{W}}+\operatorname{d}_{\operatorname{K}})\!\left({#1}\,,\,{#2}\right)}

\title[The structure of fluctuations in homogenization]{The structure of fluctuations \\in stochastic homogenization}
\author[M. Duerinckx]{Mitia Duerinckx}
\author[A. Gloria]{Antoine Gloria}
\author[F. Otto]{Felix Otto}
\address[Mitia Duerinckx]{Laboratoire de Mathématique d'Orsay, UMR 8628, Université Paris-Sud, F-91405 Orsay, France \& Universit\'e Libre de Bruxelles, Département de Mathématique, Brussels, Belgium}
\email{mduerinc@ulb.ac.be}
\address[Antoine Gloria]{Sorbonne Universit\'e, CNRS, Universit\'e de Paris, Laboratoire Jacques-Louis Lions (LJLL), F-75005 Paris, France \& Universit\'e Libre de Bruxelles, Département de Mathématique, Brussels, Belgium}
\email{gloria@ljll.math.upmc.fr}
\address[Felix Otto]{Max Planck Institute for Mathematics in the Sciences \\ Leipzig, Germany}
\email{otto@mis.mpg.de}
\begin{document}

\addtocontents{toc}{\protect\setcounter{tocdepth}{1}}

\begin{abstract}
Four quantities are fundamental in homogenization of elliptic systems in divergence form and in its applications: the field and the flux of the solution operator (applied to a general deterministic right-hand side), and the field and the flux of the corrector.
Homogenization is the study of the large-scale properties of these objects. In case of random coefficients, these quantities fluctuate and their fluctuations are a priori unrelated.
Depending on the law of the coefficient field, and in particular on the decay of its correlations on large scales, these fluctuations may display different scalings and different limiting laws (if any).
{\color{black}In this contribution, we identify another crucial intrinsic quantity,  motivated by H-convergence,
which we refer to as the \emph{homogenization commutator} and is related to variational quantities first considered by Armstrong and Smart.
In the simplified setting of the random conductance model,
we show what we believe to be a general principle, namely that the homogenization commutator drives at leading order the fluctuations of each of the four other quantities in a strong norm in probability, which is expressed in form of a suitable two-scale expansion and reveals the \emph{pathwise structure} of fluctuations in stochastic homogenization.
In addition, we show   that the 
(rescaled) homogenization commutator converges in law to a Gaussian white noise,
and we analyze to which precision the covariance tensor that characterizes the latter can be extracted from the representative volume element method. This collection of results constitutes a new theory of fluctuations in stochastic homogenization that holds in 
any dimension and yields optimal rates.
Extensions to the (non-symmetric) continuum setting are also discussed, the details of which are postponed to forthcoming works.}
\end{abstract}

\maketitle
\tableofcontents

\section{Introduction}

This article constitutes the first part of a series of works that develops a new theory of fluctuations in stochastic homogenization of elliptic (non-necessarily symmetric) systems.
{\color{black}In this first contribution, we provide a full picture of our theory with optimal convergence rates in the simplified setting of the random conductance model, that is, for discrete elliptic equations with independent and identically distributed~(iid) conductances.
The extension of our results to more general frameworks is shortly described in Section~\ref{chap:extensions} below and is postponed to forthcoming works, and a thorough discussion of the literature is provided in Section~\ref{chap:previous-lit}.}

\subsection{General overview}
Although in the sequel we shall focus on the case of discrete elliptic equations, we use non-symmetric continuum notation in this introduction.
Let $\Aa$ be a stationary and ergodic random coefficient field on $\R^d$ that satisfies
the boundedness and ellipticity properties 
\begin{equation*}
|\Aa(x)\xi|\le|\xi|,\qquad\xi\cdot\Aa(x)\xi\ge\lambda|\xi|^2,\qquad\mbox{for all}\;x,\xi\in\mathbb{R}^d,
\end{equation*}
for some $\lambda>0$.
For $\e>0$, we set $\Aa_\e:=\Aa(\frac \cdot \e)$,
and for all deterministic vector fields $f\in C^\infty_c(\R^d)^d$, we consider the random family $(u_\e)_{\e>0}$ of unique Lax-Milgram solutions in~$\R^d$ (which in the rest of this article means the unique weak solutions in $\dot H^1(\R^d)$) of the rescaled problems
\begin{align}\label{eq:first-def-ups}
-D\cdot \Aa_\e D u_\e\,=\, D\cdot f,
\end{align}
where $D$ denotes the continuum gradient (while the notation $\nabla$ is reserved in the sequel for the discrete gradient).
It is known since the pioneering work of Papanicolaou and Varadhan~\cite{PapaVara} and of Kozlov~\cite{Kozlov-79} that, almost surely, $u_\e$ converges weakly in $\dot H^1(\R^d)$
as $\e \downarrow 0$ to the unique Lax-Milgram solution $\bar u$ in $\R^d$ of 
\begin{align}\label{eq:first-def-ubar-intro}
-D\cdot \bar\Aa D \bar u\,=\, D\cdot f,
\end{align}
where $\bar\Aa$ is a deterministic and constant matrix that only depends on $\Aa$.
More precisely, for any direction $e\in \R^d$, the projection $\bar\Aa e$ is the expectation of the flux of the corrector in the direction $e$,
$$
\bar\Aa e\,=\,\expec{\Aa(D \phi_e+e)},
$$
where the corrector $\phi_e$ is the unique (up to a random additive constant) almost sure solution of the corrector equation in $\R^d$,
$$
-D \cdot \Aa(D  \phi_e+e)\,=\,0,
$$
in the class of functions the gradient of which is stationary, has vanishing expectation, and has finite second moment. We denote by $\phi=(\phi_i)_{i=1}^d$ the vector field the entries of which are the correctors $\phi_i$ in the canonical directions $e_i$ of $\R^d$.
Note that the convergence of $u_\e$ to $\bar u$ in $\dot H^1(\R^d)$ is only weak since $D u_\e$ typically displays spatial oscillations at scale $\e$, which are not captured by the limit $D \bar u$.
These oscillations are however well-described by those of the corrector field $D \phi(\frac \cdot \e)$ through the following two-scale expansion,\footnote{We systematically use Einstein's summation rule on repeated indices.}
\begin{align}\label{eq:2-scale-usual}
D u_\e\,\approx\, (D \phi_i(\tfrac\cdot\e)+\ee_i)D_i \bar u,
\end{align}
in the sense that $D u_\e- (D \phi_i(\frac\cdot\e)+\ee_i)D_i \bar u$
converges strongly to zero in $\Ld^2(\R^d)^d$. In the random setting, this theory of oscillations was recently optimally quantified in~\cite{GNO2,GNO-quant,GO4,AKM2} (see also~\cite[Chapter~6]{AKM-book}).

\medskip
As opposed to periodic homogenization, which boils down to the sole understanding of the (spatial) oscillations of $D u_\e$, the stochastic setting involves the (random) fluctuations of $D u_\e$ next to its oscillations.
More precisely, whereas oscillations reflect the (almost sure) lack of strong compactness for $Du_\e$ in $\Ld^2(\R^d)^d$, fluctuations are concerned with the leading-order probabilistic behavior of weak-type expressions of the form $\int_{\R^d}g\cdot Du_\e$ for $g\in C^\infty_c(\R^d)^d$.
Let us emphasize that in the case of a weakly correlated coefficient field~$\Aa$ the error in the two-scale expansion~\eqref{eq:2-scale-usual} is of order $\e$ in $\Ld^2(\R^d)^d$ (or $\e\Log^\frac12$ for $d=2$) while fluctuations of $Du_\e$ display the central limit theorem (CLT) scaling $\e^\frac d2$, so that~\eqref{eq:2-scale-usual} is not expected to be accurate in that scaling.
This was indeed first realized in dimension $d\ge3$ by Gu and Mourrat~\cite[Section~3.2]{GuM} (see also the last item in Remarks~\ref{rem:convcov-corr} below for $d=2$), who further argue that accuracy in~\eqref{eq:2-scale-usual} in the fluctuation scaling cannot even be reached by the use of higher-order correctors.
The corrector field $D\phi$ is therefore the driving quantity for oscillations but a priori not for fluctuations.

\medskip
In the present article, we develop a new theory of fluctuations in stochastic homogenization that builds on the known theory of oscillations, and our main achievement is the rigorous identification of the corresponding driving quantity for fluctuations.
The key consists in focusing on the \emph{homogenization commutator of the solution},
\begin{align}\label{intro:HC-sol}
\Aa_\e D u_\e-\bar\Aa D u_\e,
\end{align}
and in studying its relation to the \emph{(standard) homogenization commutator} $\Xi:=(\Xi_i)_{i=1}^d$ defined by
\begin{align}\label{intro:HC}
\Xi_i \,:=\,\Aa (D \phi_i+\ee_i)-\bar\Aa(D \phi_i+\ee_i),\qquad \Xi_{ij}:=(\Xi_i)_j.
\end{align}
{\color{black}We first briefly comment on the special form of these quantities.
As well-known in applications, homogenization is the rigorous version of averaging fields and fluxes in a consistent way, as made precise in the very definition of H-convergence by Murat and Tartar~\cite{MuratTartar}, which indeed requires both weak convergence of the fields $D u_\e \rightharpoonup D \bar u$ 
and of the fluxes $\Aa_\e D u_\e \rightharpoonup \bar\Aa D \bar u$ in $\Ld^2(\R^d)^d$, to the effect of
$$
\Aa_\e D u_\e - \bar\Aa D u_\e\, \rightharpoonup\, 0.
$$
This weak convergence of the homogenization commutator~\eqref{intro:HC-sol} is the mathematical formulation of the so-called Hill-Mandel relation in mechanics~\cite{H63,H72}.
Applied to the corrector, this justifies the definition of the (standard) homogenization commutator $\Xi$ in~\eqref{intro:HC}, which is thus seen as a natural and intrinsic measure of the accuracy of homogenization for large-scale averages. In addition, as pointed out in Section~\ref{chap:previous-lit}, this commutator turns out to be related to variational quantities first considered by Armstrong and Smart~\cite{AS}.
In these terms, our theory of fluctuations consists of the following three main principles.}
\begin{enumerate}[(I)]
\item First and most importantly,
the two-scale expansion of the {homogenization commutator} of the solution
\begin{align}\label{eq:2-scale-commut}
\Aa_\e D u_\e-\bar\Aa D u_\e-\expec{\Aa_\e D u_\e-\bar\Aa D u_\e}\,\approx\,\Xi_i(\tfrac \cdot \e) D_i \bar u
\end{align}
is (generically) accurate in the fluctuation scaling, in the sense of
\begin{multline}\label{eq:2-scale-commut-quant}
\hspace{1.1cm}\expec{\Big|\int_{\R^d}g\cdot \big(\Aa_\e D u_\e-\bar\Aa D u_\e-\expec{\Aa_\e D u_\e-\bar\Aa D u_\e}\big)-\int_{\R^d}g\cdot \Xi_i(\tfrac \cdot \e) D_i \bar u\Big|^2}^\frac12\\
 \le ~o(1)\,\expec{\Big|\int_{\R^d}g\cdot \Xi_i(\tfrac\cdot\e) D_i\bar u \Big|^2}^\frac12,
\end{multline}
where $o(1)\downarrow0$ as $\e\downarrow0$, for all $g\in C^\infty_c(\R^d)^d$. Let us emphasize again that this property is nontrivial and is due to the special form of the commutator, while it is not true for~\eqref{eq:2-scale-usual}.
\smallskip\item Second, both the fluctuations of the field $Du_\e$ and of the flux $\Aa_\e D u_\e$ can be recovered through \emph{deterministic} projections of the fluctuations of the homogenization commutator~\eqref{intro:HC-sol}, which shows that no information is lost by passing to the commutator. {\color{black}More precisely, the following elementary identities are easily checked,
\begin{eqnarray}
\int_{\R^d} g\cdot (Du_\e-D\bar u)&=&-\int_{\R^d} (\bar\calP_H^*g) \cdot (\Aa_\e Du_\e-\bar\Aa Du_\e),\nonumber\\
\int_{\R^d} g\cdot (\Aa_\e D u_\e-\bar\Aa D\bar u)&=&\int_{\R^d} (\bar\calP_L^*g)\cdot (\Aa_\e Du_\e-\bar\Aa Du_\e),\label{eq:rel-I1}
\end{eqnarray}
in terms of the Helmholtz and Leray projections in $\Ld^2(\R^d)^d$,}
\begin{gather}
\bar\calP_H:=D(D\cdot\bar\Aa D)^{-1}D\cdot,\qquad \bar\calP_L:=\Id-\bar\calP_H\bar\Aa,\nonumber\\
\bar\calP_H^*:=D(D\cdot\bar\Aa^* D)^{-1}D\cdot,\qquad \bar\calP_L^*:=\Id-\bar\calP_H\bar\Aa^*,\label{eq:proj-Helm-def}
\end{gather}
where $\bar\Aa^*$ denotes the transpose of $\bar\Aa$.
Similarly, the fluctuations of the field $D\phi$ and of the flux $\Aa D\phi$ of the corrector are determined by those of the standard commutator $\Xi$ itself: indeed, the definition of $\Xi$ yields $-D \cdot\bar \Aa D \phi_i=D \cdot \Xi_i$ and $\Aa (D \phi_i+\ee_i)-\bar\Aa \ee_i= \Xi_i+\bar\Aa D \phi_i$, to the effect of $D \phi_i=-\bar\calP_H\Xi_i$ and $\Aa(D \phi_i+\ee_i)-\bar\Aa\ee_i=(\Id-\bar\Aa\bar\calP_H)\Xi_i$ in the stationary sense, hence formally,
\begin{eqnarray}
\int_{\R^d}\F:D\phi(\tfrac\cdot\e)&=&-\int_{\R^d}(\bar\calP_H^*\F):\Xi(\tfrac\cdot\e),\nonumber\\
\int_{\R^d}\F:\big(\Aa_\e(D \phi(\tfrac\cdot\e)+\Id)-\bar\Aa\big)&=&\int_{\R^d}(\bar\calP_L^*\F):\Xi(\tfrac\cdot\e),\label{eq:rel-J1J2}
\end{eqnarray}
where $\bar\calP_H^*$ and $\bar\calP_L^*$ act on the second index of the tensor field $\F$; a suitable sense to these identities is given as part of Corollary~\ref{cor:pathwise} below.
\smallskip\item Third, the standard homogenization commutator $\Xi$ is an approximately local function of the coefficients $\Aa$, which allows to infer the large-scale behavior of $\Xi$ from the large-scale behavior of $\Aa$ itself.
{\color{black}While this locality property does not hold for the corrector $D\phi$ itself, it makes $\Xi$ a particularly relevant quantity for fluctuations.
This property is best seen when formally computing the so-called ``vertical'' derivatives of~$\Xi$ with respect to $\Aa$:} Letting $\phi^*$ denote the corrector associated with the pointwise transpose coefficient field $\Aa^*$, and letting $\sigma^*$ denote the corresponding flux corrector (cf.~\eqref{si.5} below), we obtain (cf.~\eqref{eq:der-Xi-claim})
\begin{multline}\label{eq:Xi-very-local}
\hspace{1cm}\frac{\partial}{\partial\Aa(x)}\Xi_{ij} = (D\phi_j^*+\ee_j)\cdot\frac{\partial\Aa}{\partial\Aa(x)}(D\phi_i+\ee_i)\\
-D\cdot\bigg(\phi_j^*\frac{\partial\Aa}{\partial\Aa(x)}(D\phi_i+\ee_i)\bigg)-D\cdot\bigg((\phi_j^*\Aa+\sigma^*_j)\frac{\partial D\phi_i}{\partial\Aa(x)}\bigg).
\end{multline}
In view of $\frac{\partial\Aa}{\partial\Aa(x)}=\delta(\cdot-x)$, the first right-hand side term reveals an exactly local dependence upon $\Aa$. The second term is exactly local as well, but since it is written in divergence form its contribution is negligible when integrating on large scales.
The only non-local effect comes from the last term due to $\frac{\partial D\phi}{\partial\Aa}$, which is given by the mixed derivative of the Green's function for $-D\cdot\Aa D$ and thus is expected to have only borderline integrable decay. However, it also appears inside a divergence, hence it is negligible when integrated on large scales.
In fact, in this work, the accuracy in~\eqref{eq:2-scale-commut} is established relying on a similar representation of the vertical derivative $\frac{\partial}{\partial\Aa}(\Aa_\e Du_\e-\bar\Aa D u_\e-\Xi_i(\tfrac\cdot\e) D_i\bar u)$, cf.~Lemma~\ref{lem:decompI3eps}.
\end{enumerate}

\medskip
{\color{black}
We first comment on the structure of fluctuations unravelled in~(I)--(II).
Combined with the two-scale expansion~\eqref{eq:2-scale-commut-quant} of commutators, identities~\eqref{eq:rel-I1} and~\eqref{eq:rel-J1J2} imply that the fluctuations of $Du_\e$, $\Aa_\e Du_\e$, $D\phi(\tfrac\cdot\e)$, and $\Aa_\e D\phi(\tfrac\cdot\e)$ are determined at leading order by those of the homogenization commutator $\Xi(\tfrac\cdot\e)$, with error estimated in a strong norm in probability.
In order to emphasize that it is not only a comparison of limiting laws,
this strong relation between fluctuations is henceforth referred to as the~\emph{pathwise structure} of fluctuations, in analogy to the language of SPDE.
While the classical two-scale expansion~\eqref{eq:2-scale-usual} provides a description of oscillations of a general solution by means of an off-line procedure using the corrector, a similar result is provided here for fluctuations, where the key driving quantity is now the standard commutator~$\Xi$;
this reduction of complexity for fluctuations in stochastic homogenization is bound to affect multi-scale computing and uncertainty quantification in an essential way.

\medskip
If $\Aa$ is a weakly correlated coefficient field, the locality property~(III) ensures that the homogenization commutator~$\Xi$ should also have weak correlations, so that one may expect the rescaling $\e^{-\frac d2}\Xi(\tfrac\cdot\e)$ to converge in law to a Gaussian white noise $\Gamma$.
In this case, 
the combination of~\eqref{eq:2-scale-commut-quant} with identities~\eqref{eq:rel-I1} and~\eqref{eq:rel-J1J2} leads to the joint convergence in law
\begingroup
\allowdisplaybreaks
\begin{align}
&\bigg(\e^{-\frac d 2}\int_{\R^d}\F:\Xi(\tfrac\cdot\e)~,~\e^{-\frac d 2}\int_{\R^d}g\cdot (D u_\e-\expec{D u_\e})~,~\e^{-\frac d 2}\int_{\R^d} g\cdot (\Aa_\e D u_\e-\expec{\Aa_\e D u_\e})\,,\nonumber\\
&\hspace{5cm}\e^{-\frac d 2}\int_{\R^d}\F:D\phi(\tfrac\cdot\e)~,~\e^{-\frac d 2}\int_{\R^d}\F:\big(\Aa_\e (D\phi(\tfrac\cdot\e)+\Id)-\bar\Aa\big)\bigg)\nonumber\\
&\hspace{0.5cm}\,\to\, \Big(\Gamma(\F)~,~\Gamma\big(\bar\calP_Hf\otimes\bar\calP_H^*g\big)~,~-\Gamma(\bar\calP_Hf\otimes\bar\calP_L^*g)~,~
-\Gamma\big(\bar\calP_H^*\F\big)~,~
\Gamma(\bar\calP_L^* F)
\Big).\label{eq:pathwise-res-qual}
\end{align}
\endgroup
The {pathwise structure} of fluctuations is manifested by the fact that the same random field~$\Gamma$ describes at the same time the different marginals of the limiting joint law.
For a weakly correlated coefficient field $\Aa$, we recover the previously known (or at least expected) scaling limit results for the different quantities of interest in stochastic homogenization.
As such, one could be tempted to reduce (I)--(III) to an (optimal) quantification of the joint convergence in law~\eqref{eq:pathwise-res-qual}.
However, as emphasized in Section~\ref{chap:extensions},
this is no longer true in a more general context, for instance in the case of coefficient fields with thick tails,
for which the pathwise structure (I)--(II) still holds true whereas the convergence in law~\eqref{eq:pathwise-res-qual} may fail.
From a broader perspective, the main novelty of the present contribution does not rely on such convergence results per se, but rather in uncovering the mechanism that leads to them, which is summarized in terms of properties~(I)--(III) above.
As discussed in Section~\ref{chap:previous-lit}, similar results had been conjectured in some form by Armstrong, Gu, and Mourrat in~\cite{GuM,MourratFun} based on the variational and renormalization perspective of~\cite{AS,Armstrong-Mourrat-16,AKM1,AKM2,AKM-book};
this heuristic has not been made rigorous yet.
}

\subsection{Main results}
In order to present our theory of fluctuations and address the above properties~(I)--(III), we place ourselves in the simplest setting possible and focus on the random conductance model, that is, the discrete setting with iid conductances, which we think of as the prototype for weakly correlated coefficient fields. Although conceptually simpler on the stochastic side, the discrete setting adds some technical inconveniences on the deterministic side, including a discretization error.
Our main result takes on the following guise; precise notation and assumptions on the model are postponed to Section~\ref{sec:main-results-fluc}, as well as many remarks and corollaries.
While items~(i) and~(ii) below (together with the non-degeneracy in~(iv)) imply property~(I) in the form~\eqref{eq:2-scale-commut-quant} with the optimal rate $o(1)\simeq_{f,g}\e\mu_d(\tfrac1\e)^\frac12$, items~(iii) and~(iv) are manifestations of the locality property~(III).
{\color{black}Throughout, the notation $\lesssim_\gamma$ (with possibly a subscript $\gamma$) stands for $\le$ up to a multiplicative constant $C_\gamma$ that only depends on $d$, $\lambda$, and~$\gamma$ (through a suitable norm of $\gamma$, should $\gamma$ be a function).}

\begin{theor}\label{th:main-1}
Consider the (iid) random conductance model, and assume that the law of conductances is non-degenerate. Then the following hold for all $\e>0$,
\begin{enumerate}[(i)]
\item \emph{CLT scaling:} For all $\F\in C^\infty_c(\R^d)^{d\times d}$,
\[\expec{\Big|\e^{-\frac d2}\int_{\R^d}\F: \Xi(\tfrac\cdot\e)\Big|^2}^\frac12\,\lesssim_\F\,1.\]
\item \emph{Pathwise structure (with optimal error estimate):} For all $f,g\in C^\infty_c(\R^d)^{d}$, letting $u_\e$ and $\bar u$ denote the solutions of (the discrete version of) \eqref{eq:first-def-ups} and of~\eqref{eq:first-def-ubar-intro},
\begin{multline}\label{eq:accuracy-2scale-comm}
\hspace{0.6cm}\expec{\Big|\e^{-\frac d2}\int_{\R^d}g\cdot \big(\Aa_\e \nabla_\e u_\e-\bar\Aa \nabla_\e u_\e-\expec{\Aa_\e \nabla_\e u_\e-\bar\Aa \nabla_\e u_\e}\big)-\e^{-\frac d2}\int_{\R^d}g\cdot\Xi_i(\tfrac \cdot \e) D_i \bar u\Big|^2}^\frac12 \\
\lesssim_{f,g}\,\e\mu_d(\tfrac1\e)^\frac12,
\end{multline}
where we set for all $r>0$,
\begin{align}\label{eq:def-mud}
\mu_d(r):=\begin{cases}
r&:~~d=1,\\
\log(2+r)&:~~d=2,\\1&:~~d>2.\end{cases}
\end{align}
\item \emph{Asymptotic normality (with optimal rate):} For all $\F\in C^\infty_c(\R^d)^{d\times d}$,
\begin{align*}
\delta_\calN\bigg(\e^{-\frac d2}\int_{\R^d}\F:\Xi(\tfrac\cdot\e)\bigg)\,\lesssim_\F \,\e^{\frac d2}\log(2+\tfrac1\e),
\end{align*}
where for a random variable $X\in\Ld^2(\Omega)$ its distance to normality is defined by
\begin{align}\label{eq:def-Delta-N}
\delta_\calN(X):=\dW{\frac{X}{\var{X}^\frac12}}\calN+\dK{\frac{X}{\var{X}^\frac12}}\calN,
\end{align}
with $\calN$ a standard Gaussian random variable and with $\dW\cdot\cdot$ and $\dK\cdot\cdot$ the Wasserstein and Kolmogorov metrics, {\color{black}that is, for random variables $X,Y\in\Ld^1(\Omega)$,
\begin{eqnarray*}
\dW{X}{Y}&:=&\sup\big\{\expec{f(X)-f(Y)}:f\in C^1(\R),\,\|\nabla f\|_{\Ld^\infty}\le1\big\},\\
\dK{X}{Y}&:=&\sup\big\{\big|\pr{X\le x}-\pr{Y\le x}\big|:x\in\R^d\big\}.
\end{eqnarray*}}
\item \emph{Convergence of the covariance structure (with optimal rate):} There exists a non-degenerate symmetric $4$-tensor $\calQ$ such that for all $\F\in C^\infty_c(\R^d)^{d\times d}$,
\[\bigg|\var{\e^{-\frac d2}\int_{\R^d}\F:\Xi(\tfrac\cdot\e)}-\int_{\R^d}\F:\calQ\,\F\bigg|\,\lesssim_{\F}\,\e\mu_d(\tfrac1\e)^\frac12.\]
In particular, combined with item~(iii), this yields the convergence in law of $\e^{-\frac d2}\Xi(\tfrac\cdot\e)$ to a (matrix-valued) Gaussian white noise $\Gamma$ with covariance structure $\calQ$, and the (discrete version of the) joint convergence result~\eqref{eq:pathwise-res-qual} follows.\qedhere
\end{enumerate}
\end{theor}

This fluctuation theory is complemented by the following characterization of the fluctuation tensor $\calQ$ by periodization in law. This characterization comes in form of a representative volume element (RVE) method, for which we give the optimal error estimate.
In particular, comparing with the results for the RVE approximation $\bar\Aa_{L,N}$ of the homogenized coefficients $\bar\Aa$ (cf.~\cite{GO1,GO2,GNO1}), and choosing $N\simeq L^d$ below, we may conclude that an RVE approximation for $\calQ$ with accuracy $O(L^{-\frac d2})$ (up to logarithmic corrections) is extracted at the same cost as an RVE approximation for $\bar\Aa$ with accuracy $O(L^{-d})$.
Precise assumptions and notation are again postponed to Section~\ref{sec:main-results-fluc}.

\begin{theor}\label{th:main-2}
Consider the (iid) random conductance model.
Define
\begin{align}\label{eq:def-barA-L}
\bar\Aa_{L}\ee_i:=\fint_{Q_L}\Aa_L(\nabla\phi_{L,i}+\ee_i),
\end{align}
in terms of the $L$-periodized coefficient field $\Aa_L$ and corrector $\phi_L$.
Then the fluctuation tensor $\calQ$ defined in Theorem~\ref{th:main-1}(iv) satisfies
\begin{align}\label{eq:calQ-charact}
\calQ=\lim_{L\uparrow\infty}\var{L^{\frac d2}\bar\Aa_L^*}.
\end{align}
In addition, considering iid realizations $(\Aa_L^{(n)})_{n=1}^N$ of $\Aa_L$ and setting $\bar\Aa_L^{(n)}:=\bar\Aa_L(\Aa_L^{(n)})$, we define the RVE approximation
as the square of the sample standard deviation 
\begin{align}\label{eq:def-QL-per}
\calQ_{L,N}:=\frac{L^d}{N-1}\sum_{n=1}^N\big(\bar\Aa_L^{(n)}-\bar\Aa_{L,N}\big)^*\otimes \big(\bar\Aa_L^{(n)}-\bar\Aa_{L,N}\big)^*,\qquad \bar\Aa_{L,N}:=\frac1N\sum_{n=1}^N\bar\Aa_{L}^{(n)},
\end{align}
and for all $L,N\ge2$ there holds
\[|\var{\calQ_{L,N}}|^\frac12\lesssim N^{-\frac12},\qquad|\expec{\calQ_{L,N}}-\calQ|\lesssim L^{-\frac d2}\log^{\frac d2}L.\qedhere\]
\end{theor}

{\color{black}
The ingredients to the proof of the above results are mainly twofold. The first one consists of specific concentration inequalities for iid conductances, which reduce various questions to sensitivity estimates, that is, to estimates  of ``vertical derivatives'' with respect to the coefficient field $\Aa$ (which quantify how random solutions are affected by changing the conductance at one edge). This line of argument in stochastic homogenization originates in an unpublished work by Naddaf and Spencer~\cite{NS} (see also~\cite{NS2}), and was considerably exploited in recent years starting with the early contributions of the last two authors~\cite{GO1,GO2}.
More precisely, items~(i) and~(ii) in Theorem~\ref{th:main-1} are established using a Poincaré inequality in the probability space (cf.\@ Lemma~\ref{lem:var}), item~(iii) using a second-order Poincaré inequality due to Chatterjee~\cite{C1,C2} (cf.\@ Lemma~\ref{lem:chat0}), and item~(iv) using (an iid version of) the so-called Helffer-Sjöstrand representation formula for variances~\cite{HS-94,Sjostrand-96,NS2} (cf.~Lemma~\ref{lem:cov-inequ-improved}).
Apart from these convenient tools specific to the random conductance model (which can be either extended~\cite{DG1,DG2,DG3} or avoided in some settings, cf.~Section~\ref{chap:extensions}), the proofs only rely on arguments that extend to the continuum setting and to the case of systems. 
The second main ingredient consists of large-scale regularity theory in form of a large-scale weighted Calder\'on-Zygmund theory for the random elliptic operator $-\nabla\cdot\Aa\nabla$ (cf.~\cite{Armstrong-Daniel-16,GNO-reg,AKM-book}), which is exploited here to properly estimate vertical derivatives of quantities of interest.
Large-scale regularity originates in the work of Armstrong and Smart~\cite{AS} (see also annealed regularity in~\cite{MaO}) and has been recently fully developed in~\cite{Armstrong-Mourrat-16,AKM1,AKM2,AKM-book} and in~\cite{GNO-reg,GNO-quant,GO4}.
}

\subsection{Extensions and robustness of results}\label{chap:extensions}
{\color{black}
The present contribution presents our new theory of fluctuations in the simplest possible setting of the iid conductance model.
As briefly described below, this theory will be extended to more general settings and in different directions in~the companion articles~\cite{DGO2,DFG1,DFG2,DO1,DG-19}. These developments highlight in particular the robustness of the pathwise structure of fluctuations unravelled here.}

\medskip
{\color{black}A first extension concerns the (non-symmetric) continuous setting.
In~\cite{DGO2,DFG1}, in the spirit of~\cite{DG1,DG2,DG3}, we consider the case when the coefficient field satisfies some proxy for the specific concentration inequalities of the iid setting; this actually covers all the models considered in the reference textbook~\cite{Torquato-02} on random heterogeneous materials. In addition, we analyze the case of strongly correlated coefficient fields.}
More precisely, we mainly focus on the model framework of a coefficient field given by (a local function of) a Gaussian field that has an algebraically decaying (not necessarily integrable) covariance function $c$, 
say at some fixed rate $c(x)\simeq(1+|x|)^{-\beta}$ parametrized by $\beta>0$.
For such coefficient fields, we establish in \cite{DGO2} the validity of the two-scale expansion~\eqref{eq:accuracy-2scale-comm} of the homogenization commutator in the suitable fluctuation scaling (that is, $\e^{-\frac d 2}$ is replaced by $\e^{-\frac12(\beta \wedge d)}$, with a logarithmic correction when $\beta=d$).
The proof relies on a weighted version of a Poincaré inequality in the probability space (cf.~\cite{DG1,DG2}), which can alternatively be reformulated in terms of Malliavin calculus, together with the available quantitative homogenization theory in that specific setting~\cite{GNO-reg,GNO-quant}.
This result illustrates the robustness of the pathwise structure with respect to the large-scale behavior of the homogenization commutator.
Indeed, in dimension $d=1$ (in which case the quantities under investigation are simpler and explicit), two typical behaviors have been identified in terms of the scaling limit of $\Xi$, depending on the parameter~$\beta$ (cf.~\cite{BGMP-08}):
\begin{itemize}
\item For $\beta>d=1$, the commutator $\Xi$ displays the CLT scaling and $\e^{-\frac d2}\Xi(\frac\cdot\e)$ converges to a Gaussian white noise (Gaussian fluctuations, local limiting covariance structure), {\color{black}but the convergence rate is arbitrarily slow as $\beta$ gets closer to~$d=1$.}
\item For $0<\beta<d=1$, the suitable rescaling $\e^{-\frac \beta2}\Xi(\frac\cdot\e)$ generically converges along a subsequence to a fractional Gaussian field (Gaussian fluctuations, nonlocal limiting covariance structure, potentially no uniqueness of the limit). (Note that a different, non-Gaussian behavior may also occur in this case, cf.~\cite{Gu-Bal-12,LNZH-17}.)
\end{itemize}
{\color{black}In particular, the pathwise result is shown to hold in both cases with rate $\e^{\frac\beta2\wedge1}$ (up to a logarithmic correction) whereas the rescaled homogenization commutator does not necessarily converge to white noise, may converge arbitrarily slowly, or may even  have no limit at all.
As already mentioned, this advocates that the pathwise structure of fluctuations in form of the two-scale expansion~\eqref{eq:accuracy-2scale-comm} and the scaling limit of the homogenization commutator are two separate properties that only \emph{partially} overlap.
The identification of the scaling limit in higher dimensions with optimal rates is addressed in~\cite{DFG1} for the whole range of values of $\beta>0$ (see also~\cite{DO1} for $\beta>d$),
where we further exploit the locality of $\Xi$ with respect to the coefficient field together with tools from Malliavin calculus;
this work extends~\cite{BGMP-08} to dimensions $d\ge2$.

\medskip
{\color{black}A second extension concerns the continuous setting without product space assumption (that is, without using any proxy for the specific concentration inequalities of the iid setting),
and more specifically we focus on random coefficient fields with finite range of dependence.
The convergence in law of the rescaled homogenization commutator $\e^{-\frac d2}\Xi(\tfrac\cdot\e)$ to a Gaussian white noise (albeit without optimal convergence rate) was obtained in that case independently in~\cite{AKM2} and~\cite{GO4}.
The proof of the validity of the two-scale expansion~\eqref{eq:accuracy-2scale-comm} of the homogenization commutator in the CLT scaling is more involved and will be presented in~\cite{DFG2} based on the semi-group approach of~\cite{GO4}, as well as an optimal convergence rate for the scaling limit of the commutator based on~\cite{Fischer-19}.}

\medskip
A third extension concerns higher-order corrections to fluctuation theory: in~\cite{DO1}, a suitable higher-order version of the homogenization commutator is identified and its higher-order two-scale expansion is shown to indeed have higher-order accuracy in the CLT scaling, which in particular leads to a higher-order description of fluctuations of the solution operator, while the scaling limit of the higher-order standard homogenization commutator is again easily computed in view of a corresponding higher-order locality property. In other words, the key properties~(I)--(III) are fully extended to higher order, and this theory is in line with the usual higher-order theory of oscillations based on higher-order correctors.

\medskip
A last extension concerns fluctuations for solutions of the wave equation with random coefficients. The homogenization commutator~\eqref{intro:HC-sol} is also relevant in that setting and similarly leads to a corresponding pathwise theory of fluctuations~\cite{DG-19}. This theory can further be extended to higher orders in terms of suitable higher-order hyperbolic correctors introduced in~\cite{DGR-19}.

\subsection{Relation to previous works}\label{chap:previous-lit}
The description of fluctuations in stochastic homogenization
has been the 
most central open question in the field since its very beginning in the late 1970s, and it has been a particularly active topic in the last few years.
We organize the discussion of previous works on fluctuations in separate parts addressing different aspects of the theory for the random conductance model.

\begin{enumerate}[$\bullet$]
\item \emph{CLT scaling.}\\
The CLT scaling for the solution operator in form of uniform moment bounds on $\e^{-\frac d2}\int_{\R^d} g\cdot (\nabla_\e u_\e-\expec{\nabla_\e u_\e})$ was first investigated by Conlon and Naddaf~\cite{CN} and by the second author~\cite{G},
and it was established in its optimal form in all dimensions by Marahrens and the third author~\cite{MaO} for the random conductance model.

\smallskip \item \emph{Asymptotic normality.} \\
The first asymptotic normality result in stochastic homogenization is due to Nolen~\cite{N,Nolen-16} (which contains the first use of second-order Poincaré inequalities in the field), see also Rossignol~\cite{R}, Biskup, Salvi, and Wolff~\cite{BSW}.

\smallskip\item \emph{Scaling limits.}\\
The covariance structure of the corrector was identified by Mourrat and the third author~\cite{MO} (for $d>2$), and the limiting variance of the energy of the periodized corrector (which a posteriori is nothing but the average of the periodized homogenization commutator) by the second author and Nolen~\cite{GN}.
These works were then extended to the solution operator by Gu and Mourrat~\cite{GuM} (for $d>2$), which essentially
yields as a corollary
the joint convergence in law~\eqref{eq:pathwise-res-qual}. These works were based on Chatterjee's second-order Poincar\'e inequalities~\cite{C1,C2} (cf.~Lemma~\ref{lem:chat0}) and on the Helffer-Sj\"ostrand representation formula~\cite{HS-94,Sjostrand-96,NS2} (cf.~Lemma~\ref{lem:cov-inequ-improved}) together with tools from the early quantitative homogenization theory as developed
by the last two authors~\cite{GO1,GO2} (inspired by the unpublished work by Naddaf and Spencer \cite{NS}), the last two authors and Neukamm~\cite{GNO1,GNO2}, and Marahrens and the last author~\cite{MaO}.
Note that the present contribution originated in the attempt to upgrade the results in~\cite{GN} into a functional CLT for a suitable energy density.

\smallskip\item\emph{Genesis of the homogenization commutator.}\\
Although strongly motivated by the very notion of H-convergence as defined by Murat and Tartar~\cite{MuratTartar}, the explicit apparition of homogenization commutators and the discovery of their remarkable properties are much more recent.
A variational quantity related to the standard commutator can first be traced back in the work of Armstrong and Smart~\cite{AS} and has been the driving quantity in the subsequent developments by Armstrong, Kuusi, and Mourrat~\cite{AKM1,AKM2,AKM-book};
the explicit link with the standard commutator in the form~\eqref{intro:HC} is provided e.g.~in~\cite[Lemma 4.25]{AKM-book}.
The idea that such a quantity might be important for fluctuations was originally formulated by Armstrong, Gu, and Mourrat in form of a heuristic~\cite{GuM,MourratFun}. The early version of the present work provided the first rigorous result in that direction, as well as a proof of the locality of the standard commutator and its convergence to white noise.
This scaling limit result for the standard commutator was also independently obtained in~\cite{AKM2} and in~\cite{GO4} for coefficient fields with finite range of dependence (although without optimal rates); in view of the obvious relation~\eqref{eq:rel-J1J2}, note that this characterizes fluctuations of the corrector.
Importantly, instead of a variational interpretation, the homogenization commutator is given here a more intrinsic and practical definition  motivated by H-convergence.

\smallskip\item\emph{Pathwise structure of fluctuations.}\\
The pathwise structure in the restricted form of the joint convergence in law~\eqref{eq:pathwise-res-qual} was essentially observed by Gu and Mourrat in~\cite{GuM}, and the first attempt at a general explanation of this structure was proposed by Armstrong, Gu, and Mourrat in form of a heuristic~\cite{GuM,MourratFun}.
No explanation was however suggested why fluctuations of the solution operator could be expressed in terms of the same intrinsic object as fluctuations of the corrector (that is, in terms of the standard commutator), which constitutes the main gap in that work.
As we show here, this question can be reformulated as the validity of a suitable two-scale expansion principle for fluctuations of the solution operator.
Since the two-scale expansion of the solution operator itself is known not to be accurate in the fluctuation scaling~\cite{GuM}, it was unclear whether a two-scale expansion principle could actually hold for fluctuations.
Our main contribution in this work precisely fills this gap in form of the accuracy of the two-scale expansion of commutators.
This key missing part in the heuristic proposed in~\cite{GuM,MourratFun} has still not been filled within a variational and renormalization perspective.
\end{enumerate}
\smallskip\noindent
To sum up, the  pathwise structure of fluctuations in form of the two-scale expansion of commutators is to fluctuations what the usual two-scale expansion of gradient fields is to oscillations.
The whole mechanism that drives fluctuations, as summarized in  properties~(I)--(III), is made precise and rigorous here for the first time in any setting.
This new theory of fluctuations turns out to be remarkably robust, as shown by the various extensions that it has triggered.
Even for the more specific question of scaling limits, this theory provides new optimal convergence rates.
}

\newpage
\section{Main results}\label{sec:main-results-fluc}

In this section, we introduce notation and assumptions on the random conductance model, we state precise versions of the main results (in particular including explicit norms of the test functions in the estimates), and we discuss various corollaries.

\subsection{Notation and assumptions}\label{chap:notations}
We start by introducing the random conductance model on the integer lattice $\Z^d$, which is the framework of our main results.
We denote by  $\{\ee_i\}_{i=1}^{d}$ the canonical basis of $\R^d$, and we regard $\Z^d$ as a graph with (unoriented) edge set $\B = \{ (x,z) \in \Z^d \times \Z^d :  |x-z|=1\}$. For edges $(x,z) \in \B$, we write $x\sim z$. We define the set of conductances $\{a(b)\}_{b\in \B}$  by $\Omega=[\lambda,1]^\B$ for some fixed $0<\lambda\le 1$. We endow $\Omega$ with the $\sigma$-algebra generated by cylinder sets and with a probability measure~$\mathbb P$. We denote by $\expec\cdot$, $\var\cdot$, and $\cov\cdot\cdot$ the associated expectation, variance, and covariance.
A realization $a \in \Omega$ is by definition a collection $\{a(b)\}_{b\in \B}$ of conductances.
A random field $u:\R^d\times\Omega\to\R$ is said to be stationary if it is shift-covariant, in the sense of $u(x,a(\cdot-z))=u(x-z,a)$ for all $x,z\in\R^d$ and $a\in\Omega$.
In this contribution, we focus on the case when the probability measure $\Pm$ is a product measure, that is, when the conductances $\{a(b)\}_{b\in\B}$ are iid random variables, and we shall make use of available functional inequalities in this product probability space.

Let $\nabla$ denote the forward discrete gradient $(u:\Z^d\to \R)\mapsto (\nabla u:\Z^d \to \R^d)$ defined componentwise by $\nabla_i u(x)=u(x+\ee_i)-u(x)$ for $1\le i\le d$, 
and let $\nabla^*$ denote the backward discrete gradient $(u:\Z^d\to \R)\mapsto (\nabla^* u:\Z^d \to \R^d)$ defined componentwise by $\nabla_i^* u(x)=u(x)-u(x-\ee_i)$ for $1\le i\le d$. The operator $-\nabla^*\cdot$ is thus the adjoint of $\nabla$ on $\ell^2(\Z^d)$,
and we consider the elliptic operator $-\nabla^*\cdot\Aa\nabla$ with coefficients 
\begin{equation*}
\Aa:x\mapsto \Aa(x):=\dig{a(x,x+\ee_1),\dots,a(x,x+\ee_d)},
\end{equation*}
acting on functions $u:\Z^d \to \R$ as
$$
-\nabla^*\cdot\Aa\nabla u(x)\,:=\,\sum_{z:z\sim x}a(x,z)(u(x)-u(z)).
$$
In order to state the standard qualitative homogenization result \cite{Ku,Kozlov-85} for the corresponding discrete elliptic equation, we consider for all $\e>0$ the rescaled operator
$-\nabla_\e^* \cdot \Aa_\e \nabla_\e$,
where $\Aa_\e(\cdot):=\Aa(\frac\cdot\e)$, and where $\nabla_\e$ and $\nabla^*_\e$ act on functions $u_\e: \Z_\e^d:=\e\Z^d\to \R$, and are defined componentwise by $\nabla_{\e,i} u_\e(x)=\e^{-1}(u_\e(x+\e\ee_i)-u_\e(x))$ and  $\nabla_{\e,i}^* u_\e(x)=\e^{-1}(u_\e(x)-u_\e(x-\e\ee_i))$ for all $i$.
We shall also let $\nabla_\e$ and $\nabla^*_\e$ act on continuous functions $u:\R^d \to \R$, so that $\nabla_\e u$ and $\nabla_\e^* u$ are continuous functions as well. If $u \in C^1(\R^d)$, then $\nabla_\e u(x)$ and $\nabla_\e^* u(x)$ converge to the continuum gradient $Du(x)$ for all $x\in \R^d$ as $\e \downarrow 0$.
In what follows, for all $m\ge1$, we systematically extend maps $v:\Z^d \to \R^m$ to piecewise constant maps $\R^d \to \R^m$ (still denoted by $v$) by setting $v|_{Q(x)}:= v(x)$ for all $x\in \Z^d$ (where $Q(x):=x+[-\frac12,\frac12)^d$ is the unit cube centered at $x$), and we use this notation e.g.\@ for $\Aa:\Z^d \to \R^{d \times d}$ (but also for $\phi:\Z^d \to \R^d$ and $\Xi:\Z^d \to \R^{d\times d}$ defined below). This systematic extension of functions from the lattice $\Z^d$ to $\R^d$ allows to state all discrete results in a form that would hold mutatis mutandis in the continuum setting.
In addition, although in the discrete setting it is more natural to consider a symmetric coefficient field~$\Aa$, we use non-symmetric notation in the statement of the results in view of the extension to the non-symmetric continuum setting, and we denote by $\Aa^*$ the pointwise transpose field associated with $\Aa$.

Qualitative stochastic homogenization~\cite{Ku,Kozlov-85} ensures that, for all $f,g\in C^\infty_c(\R^d)^d$, almost surely,
the unique Lax-Milgram solutions $u_\e$ and $v_\e$ in $\R^d$
of\footnote{
These equations are understood as follows: for all $x\in Q$ the function $u_\e(\e x+\cdot)$ on $\Z_\e^d$ is the solution of the discrete elliptic equation with coefficient $\Aa_\e$ and with right-hand side $\nabla_\e^*\cdot f(\e x+\cdot)$. This definition allows to state results in a form that holds in the continuum setting, and in terms of norms of the right-hand side that do not necessarily have to embed into the space of continuous functions.
}
\begin{equation}\label{e.def-ueps}
-\nabla_\e^* \cdot \Aa_\e \nabla_\e u_\e\,=\,\nabla_\e^*\cdot f,\qquad -\nabla_\e^* \cdot \Aa_\e^* \nabla_\e v_\e\,=\,\nabla_\e^*\cdot g,
\end{equation}
converge weakly as $\e\downarrow0$ to the unique Lax-Milgram solutions $\bar u$ and $\bar v$ in $\R^d$ of the (continuum) elliptic equations
\begin{align}\label{e.def-ubar}
-D\cdot \bar\Aa  D \bar u\,=\,D\cdot f,\qquad -D\cdot \bar\Aa^*  D \bar v\,=\,D\cdot g,
\end{align}
respectively, where $\bar\Aa$ is the homogenized matrix characterized by 
\begin{equation}\label{e:hom-coeff}
\bar\Aa \ee_i \,=\,\expec{\Aa (\nabla \phi_i+\ee_i)},
\end{equation}
for all $1\le i\le d$, and where $\phi_i$ is the so-called corrector in direction $\ee_i$. It is defined, for almost every realization $\Aa$, as the unique solution in $\Z^d$ of
\begin{equation}\label{e.corr}
-\nabla^* \cdot \Aa (\nabla\phi_i+\ee_i)\,=\,0,
\end{equation}
with $\nabla \phi_i$ stationary and having vanishing expectation and finite second moment, and with the anchoring $\phi_i(0)=0$ at the origin. We then set $\phi:=(\phi_i)_{i=1}^d$. Note that $\overline{(\Aa^*)}=(\bar\Aa)^*$. For symmetric coefficient fields, $\Aa^*=\Aa$ and $\bar\Aa^*=\bar\Aa$.

We consider the fluctuations of the field $\nabla u_\e$ and of the flux $\Aa_\e\nabla u_\e$, as encoded in the random linear functionals $I_1^\e:(f,g)\mapsto I_1^\e(f,g)$ and $I_2^\e:(f,g)\mapsto I_2^\e(f,g)$ defined
for all $f,g\in C^\infty_c(\R^d)^d$ by
\begin{eqnarray*}
I_1^\e(f,g)&:=&\e^{-\frac{d}{2}}\int_{\R^d} g\cdot\nabla_\e(u_\e-\expec{u_\e}),\\
I_2^\e(f,g)&:=&\e^{-\frac{d}{2}}\int_{\R^d} g\cdot\big(\Aa_\e\nabla_\e u_\e-\expec{\Aa_\e\nabla_\e u_\e}\big).
\end{eqnarray*}
We further encode the fluctuations of the corrector field $\nabla\phi$ and flux $\Aa(\nabla\phi+\Id)$ in the random linear functionals $J_1^\e:\F\mapsto J_1^\e(\F)$ and $J_2^\e:\F\mapsto J_2^\e(\F)$ defined for all $\F\in C^\infty_c(\R^d)^{d\times d}$ by
\begin{eqnarray*}
J_1^\e(\F)&:=&\e^{-\frac{d}{2}} \int_{\R^d} \F(x): \nabla \phi(\tfrac{x}{\e})\,dx,\\
J_2^\e(\F)&:=&\e^{-\frac{d}{2}} \int_{\R^d} \F(x):\big(\Aa_\e(x) (\nabla \phi(\tfrac{x}{\e})+\Id)-\bar\Aa\big)\,dx.
\end{eqnarray*}
As explained above, a crucial role is played by the (standard) homogenization commutator, which in the present discrete setting takes the form $\Xi:=(\Xi_i)_{i=1}^d$ with
\begin{equation}\label{eq:definition-Xi}
\Xi_i\,:=\,\Aa(\nabla \phi_i+\ee_i)-\bar\Aa (\nabla \phi_i+\ee_i),\qquad\Xi_{ij}:=(\Xi_i)_j,
\end{equation}
and by the error in the two-scale expansion of the homogenization commutator of the solution.
These quantities are encoded in the random linear functionals $J_0^\e:\F\mapsto J_0^\e(\F)$ and $E^\e:(f,g)\mapsto E^\e(f,g)$ defined for all $\F\in C^\infty_c(\R^d)^{d\times d}$ and all $f,g\in C^\infty_c(\R^d)^d$ by
\begin{eqnarray*}
J_{0}^\e(\F)\!\!&:=&\!\!\e^{-\frac{d}{2}}\int_{\R^d} \F(x) : \Xi(\tfrac{x}{\e})\,dx,\\
E^\e(f,g)\!\!&:=&\!\!\e^{-\frac d2}\int_{\R^d}g\cdot\big(\Aa_\e \nabla_\e u_\e-\bar\Aa\nabla_\e u_\e-\expec{\Aa_\e \nabla_\e u_\e-\bar\Aa \nabla_\e u_\e}\big)-\e^{-\frac d2}\int_{\R^d}g\cdot\Xi_i(\tfrac \cdot\e)D_i\bar u.
\end{eqnarray*}
Since the case $d=1$ is much simpler and well-understood~\cite{Gu-16}, we shall only focus in the sequel on dimensions $d\ge2$.

We first recall the following uniform boundedness result for $J_0^\e$, establishing the CLT scaling for the fluctuations of the homogenization commutator (cf.\@ Theorem~\ref{th:main-1}(i)). Although essentially contained in the main result of the first contribution~\cite{GO1} of the second and third authors to the field, a short proof with up-to-date tools is included for completeness in Section~\ref{chap:decompsol}.
Note that the norm of the test function is substantially weaker than $\Ld^1(\R^d)$ in terms of integrability and is thus compatible with the behavior of Helmholtz projections of smooth and compactly supported functions, which is necessary for the pathwise result of Corollary~\ref{cor:pathwise} below.

\begin{prop}\label{prop:bounded-I0}
Let $d\ge2$, let $\Pm$ be a product measure, and set $w_1(z):=1+|z|$.
For all $\e>0$ and all $\F\in C^\infty_c(\R^d)^{d\times d}$ we have
for all $0<p-1\ll1$ and all $\alpha>d\frac{p-1}{4p}$,
\[\expec{|J_0^\e(\F)|^2}^\frac12+\expec{|J_1^\e(\F)|^2}^\frac12+\expec{|J_2^\e(\F)|^2}^\frac12\,\lesssim_{\alpha,p}\,\|w_1^{2\alpha} \F\|_{\Ld^{2p}(\R^d)}.\]
Above, and in the rest of this article, $\ll$ stands for $\le$ up to a small enough multiplicative constant $C=C(d)>0$
that only depends on $d$.
\qedhere
\end{prop}

\subsection{Pathwise structure}

Our first main result establishes the smallness of the rescaled error $E^\e$ in the two-scale expansion of the homogenization commutator (cf.\@ Theorem~\ref{th:main-1}(ii)), which is the key to the pathwise structure~\eqref{eq:pathwise-res-qual}.
As for Proposition~\ref{prop:bounded-I0}, the proof relies on the Poincar\'e inequality in the probability space that is satisfied for iid coefficients.
From a technical point of view, we exploit the large-scale Calder\'on-Zygmund theory for the random elliptic operator $-\nabla^*\cdot\Aa\nabla$ as developed in~\cite{Armstrong-Daniel-16,GNO-reg} (see also~\cite[Section~7]{AKM-book}).

\begin{prop}\label{prop:pathwise}
Let $d\ge2$, let $\Pm$ be a product measure, let $\mu_d$ be defined in~\eqref{eq:def-mud}, and set $w_1(z):=1+|z|$. For all $\e>0$ and all $f,g\in C^\infty_c(\R^d)^d$ we have for all $0<p-1\ll1$ and all $\alpha>d\frac{p-1}{4p}$,
\[\expec{|E^\e(f,g)|^2}^\frac12\,\lesssim_{\alpha,p}\,\e\mu_d(\tfrac1\e)^\frac12
\Big(\|f\|_{\Ld^{4}(\R^d)}\|w_1^{\alpha}Dg\|_{\Ld^{4p}(\R^d)}+\|g\|_{\Ld^{4}(\R^d)}\|w_1^{\alpha}Df\|_{\Ld^{4p}(\R^d)}\Big).\qedhere\]
\end{prop}

\begin{rem}
For simplicity, the estimates in Propositions~\ref{prop:bounded-I0} and~\ref{prop:pathwise} above are stated and proved for second moments only, but the same arguments yield similar estimates for all algebraic (and even stretched exponential) moments (cf.\@ \cite{DGO2,DO1,DFG2}).
\end{rem}

In view of identity~\eqref{eq:rel-I1} (which indeed holds in the discrete setting up to a higher-order discretization error), the above result implies that the large-scale fluctuations of $I_1^\e$ and $I_2^\e$ are driven by those of $J_0^\e$ in a pathwise sense. Identity~\eqref{eq:rel-J1J2} (which again holds up to a discretization error) yields a similar pathwise result for $J_1^\e$ and $J^\e_2$.

\begin{cor}\label{cor:pathwise}
Let $d\ge2$, let $\Pm$ be a product measure, let $\bar\calP_H$, $\bar\calP_H^*$, and $\bar\calP_L^*$ be defined in~\eqref{eq:proj-Helm-def}, let $\mu_d$ be defined in~\eqref{eq:def-mud}, and set $w_1(z):=1+|z|$.
For all $\e>0$, all $f,g\in C^\infty_c(\R^d)^d$, and all $\F\in C^\infty_c(\R^d)^{d\times d}$,
we have for all $0<p-1\ll1$ and all $\alpha>d\frac{p-1}{4p}$,
\begin{multline}\label{eq:pathwise-I1-st}
\expec{|I_1^\e(f,g)-J_0^\e(\bar\calP_Hf\otimes \bar\calP_H^*g)|^2}^\frac12
+\expec{|I_2^\e(f,g)+J_0^\e(\bar\calP_Hf\otimes\bar\calP_L^*g)|^2}^\frac12\\
\lesssim_{\alpha,p}
\e\mu_d(\tfrac1\e)^\frac12\,
\Big(\|f\|_{\Ld^{4}(\R^d)}\|w_1^\alpha Dg\|_{\Ld^{4p}(\R^d)}
+\| g\|_{\Ld^{4}(\R^d)}\|w_1^\alpha Df\|_{\Ld^{4p}(\R^d)}\Big),
\end{multline}
and 
\begin{align}\label{eq:pathwise-I23-st}
\expec{|J_1^\e(\F)+J_{0}^\e(\bar\calP_H^*\F)|^2}^\frac12+\expec{|J_2^\e(\F)-J_{0}^\e(\bar\calP_L^*\F)|^2}^\frac12
\lesssim_{\alpha,p}\,\e\|w_1^{2\alpha} D\F\|_{\Ld^{2p}(\R^d)},
\end{align}
where by definition we have $\bar\calP_H f=-D \bar u$ and $\bar\calP_H^*g=-D\bar v$.
In particular, we give meaning to $J_0^\e(\bar\calP_H^*\F)$ and $J_{0}^\e(\bar\calP_L^*\F)$ in $\Ld^2(\Omega)$ for all $\F\in C^\infty_c(\R^d)^{d\times d}$, even when $\bar\calP_H^*\F$ and $\bar\calP_L^*\F$ do not have integrable decay.
\end{cor}

\begin{rem}
For all $f,g\in C^\infty_c(\R^d)^d$, we may also consider the unique Lax-Milgram solutions $u_\e^\circ$ and $v_\e^\circ$ in $\R^d$ of
\[-\nabla_\e^*\cdot\Aa_\e\nabla_\e u_\e^\circ=\nabla_\e^*\cdot\Aa_\e f,\qquad-\nabla_\e^*\cdot\Aa_\e^*\nabla_\e v_\e^\circ=\nabla_\e^*\cdot\Aa_\e^* g,\]
which, almost surely, converge weakly as $\e\downarrow0$ to the unique Lax-Milgram solutions $\bar u^\circ$ and $\bar v^\circ$ in $\R^d$ of
\[-D\cdot\bar\Aa D\bar u^\circ=D\cdot\bar\Aa f,\qquad-D\cdot\bar\Aa^* D\bar v^\circ=D\cdot\bar\Aa^* g,\]
respectively.
Similar considerations as in the proof of Proposition~\ref{prop:pathwise} and Corollary~\ref{cor:pathwise} then lead to a pathwise result for the fluctuations of the random linear functionals $I_3^\e:(f,g)\mapsto I_3^\e(f,g)$ and $I_4^\e:(f,g)\mapsto I_4^\e(f,g)$ defined for all $f,g\in C^\infty_c(\R^d)^d$ by
\begin{eqnarray*}
I_3^\e(f,g)&:=&\e^{-\frac d2}\int_{\R^d}g\cdot\nabla_\e(u_\e^\circ-\expec{u_\e^\circ}),\\
I_4^\e(f,g)&:=&\e^{-\frac d2}\int_{\R^d}g\cdot \big(\Aa_\e(\nabla_\e u_\e^\circ+f)-\expec{\Aa_\e(\nabla_\e u_\e^\circ+f)}\big),
\end{eqnarray*}
that takes the form
\begin{align*}
\expec{\big|I_3^\e(f,g)+J_0^\e\big(\bar\calP_Lf \otimes \bar\calP_H^*g\big)\big|^2}^\frac12
+\expec{\big|I_4^\e(f,g)-J_0^\e\big(\bar\calP_Lf\otimes\bar\calP_L^*g\big)\big|^2}^\frac12
\lesssim_{f,g}\e\mu_d(\tfrac1\e)^\frac12,
\end{align*}
where by definition $\bar\calP_H f=-D \bar u$, $\bar\calP_Lf=D\bar u^\circ+f$, $\bar\calP_H^*g=-D\bar v$, and $\bar\calP_L^*g=D\bar v^\circ+g$.
\end{rem}

Incidentally, as a consequence of our analysis, combining the two-scale expansion of the homogenization commutator~\eqref{eq:2-scale-commut} with identity~\eqref{eq:rel-I1}, we obtain a new (nonlocal) two-scale expansion for the solution $\nabla_\e u_\e$ that is not only accurate at order $1$ for the strong $\Ld^2(\R^d)$ topology but also at the order of the CLT scaling for the weak $\Ld^2(\R^d)$ topology, in contrast to the usual two-scale expansion~\eqref{eq:2-scale-usual} (cf.~\cite{GuM}).
(The second estimate below is a reformulation of Proposition~\ref{prop:pathwise}, whereas the first estimate is a corollary of~\cite[Theorem~3]{GNO-quant}.)
\begin{cor}\label{coro:2scale}
Let $d\ge2$, let $\Pm$ be a product measure, and let $\mu_d$ be defined in~\eqref{eq:def-mud}.
For all $\e>0$ and all $f\in C^\infty_c(\R^d)^d$, we set
$$
r_\e(f):=\nabla_\e u_\e-\Big(\underbrace{\expec{\nabla_\e u_\e}+\nabla_\e(-\nabla_\e^* \cdot \bar\Aa \nabla_\e)^{-1} \nabla_\e^* \cdot\big( \Xi_i(\tfrac\cdot\e) D_{i} \bar u\big)}_{\displaystyle \text{nonlocal two-scale expansion of $\nabla_\e u_\e$}} \Big).
$$
This (nonlocal) two-scale expansion correctly captures:
\begin{itemize} 
\item the \emph{spatial} oscillations of $\nabla_\e u_\e$ in a strong norm: for all $f\in C^\infty_c(\R^d)^d$,
\begin{equation*}
\expec{\|r_\e(f)\|_{\Ld^2(\R^d)}^2}^\frac12 \,\lesssim_{f} \, \e \mu_d(\tfrac1\e)^\frac12;
\end{equation*}
\item the \emph{random} fluctuations of $\nabla_\e u_\e$ in the CLT scaling:
for all $f,g \in C^\infty_c(\R^d)^d$,
\[\expec{\Big| \e^{-\frac{d}{2}} \int_{\R^d} g\cdot r_\e(f) \Big|^2}^\frac12 \,\lesssim_{f,g} \, \e \mu_d(\tfrac1\e)^\frac12.\qedhere\]
\end{itemize}
\end{cor}

\subsection{Approximate normality}
We turn to the normal approximation result for the homogenization commutator (cf.\@ Theorem~\ref{th:main-1}(iii)), which states that the fluctuations of $\e^{-\frac d2}\Xi(\frac\cdot\e)$ are asymptotically Gaussian (up to a non-degeneracy condition that is elucidated in Proposition~\ref{prop:conv-cov} below).
The approach is inspired by previous works by Nolen~\cite{N,Nolen-16}, based on a second-order Poincar\'e inequality \`a la Chatterjee~\cite{C1,LRP-15}, which is key to optimal convergence rates.
{\color{black}Such functional inequalities are not easily amenable to the use of large-scale Calder\'on-Zygmund theory for the random elliptic operator $-\nabla^*\cdot\Aa\nabla$, and we rather have to exploit optimal annealed estimates on mixed gradients of the Green's function~\cite{MaO} (see also \cite{Bella-Giunti-18} and \cite[Section~8.5]{AKM-book}). The proof exploits the approximate locality of the homogenization commutator $\Xi$, while the lack of exact locality precisely leads to the additional $\log(2+\frac1\e)$ factor in the convergence rate (we believe that this is optimal).}

\begin{prop}\label{prop:norm-dens}
Let $d\ge 2$, let $\Pm$ be a product measure, let $\mu_d$ and $\delta_\calN$ be defined in~\eqref{eq:def-mud} and~\eqref{eq:def-Delta-N}, and set $w_1(z):=1+|z|$.
For all $\e>0$ and all $\F\in C^\infty_c(\R^d)^{d\times d}$ with $\var{J_0^\e(\F)}>0$, we have for all $\alpha>0$,
\begin{align*}
&\delta_\calN(J_0^\e(\F))
\,\lesssim_\alpha \,
\e^{\frac{d}{2}}\,\frac{\|\F\|_{\Ld^3(\R^d)}^3+\|w_1^\alpha D\F\|_{\Ld^3(\R^d)}^3}{\var{J_0^\e(\F)}^\frac32}\\
&\hspace{5cm}+\e^{\frac d2}\log(2+\tfrac1\e)\,\frac{\|w_1^\alpha \F\|_{\Ld^4(\R^d)}^2+\|w_1^\alpha D\F\|_{\Ld^4(\R^d)}^2}{\var{J_0^\e(\F)}}.\qedhere
\end{align*}
\end{prop}

\begin{rem}\label{rem:Chat-Green}
{\color{black}In the case of iid conductances that are (smooth local transformations of) Gaussian random variables, a stronger version of a second-order Poincar\'e inequality is available (cf.~\cite[Theorem~2.2]{C2}), which in addition controls the total variation distance and allows to avoid Green's function estimates (see~\cite{DO1,DFG1} for such an argument in the continuum setting).}
\end{rem}

\subsection{Covariance structure}
Since $J_0^\e$ is asymptotically Gaussian, it remains to identify the limit of its covariance structure (cf.\@ Theorem~\ref{th:main-1}(iv)). The following shows that the limiting covariance is that of a (matrix-valued) white noise with some non-degenerate covariance tensor~$\calQ$. The convergence rate in~\eqref{eq:conv-cov-main-I0} below is new in any dimension and is expected to be optimal.
The proof crucially relies on the approximate locality of the homogenization commutator and on
(an iid version of) the Helffer-Sj\"ostrand representation formula for the variance~\cite{HS-94,Sjostrand-96,NS2}, which is a stronger tool than the Poincar\'e inequality.
As for the pathwise result, the proof exploits the large-scale Calder\'on-Zygmund theory for the random elliptic operator $-\nabla^*\cdot\Aa\nabla$.

\begin{prop}\label{prop:conv-cov}
Let $d\ge2$, let $\Pm$ be a product measure, let $\mu_d$ be defined in~\eqref{eq:def-mud}, and set $w_1(z):=1+|z|$.
Then the following hold,
\begin{enumerate}[(i)]
\item There exists a symmetric\footnote{Since $\calQ$ is a (limiting) covariance, it is of course symmetric in the sense of $\calQ_{ijkl}=\calQ_{klij}$. If the coefficients $\Aa$ are symmetric, then it has the additional symmetry $\calQ_{ijkl}=\calQ_{jikl}$ (hence also $\calQ_{ijkl}=\calQ_{ijlk})$.} $4$-tensor $\calQ$ such that for all $\e>0$ and all $\F,\G\in C^\infty_c(\R^d)^{d\times d}$
we have for all $0<p-1\ll1$ and all $\alpha>d\frac{p-1}{4p}$,
\begin{multline}\label{eq:conv-cov-main-I0}
\hspace{1cm}\bigg|\cov{J_0^\e(\F)}{J_0^\e(\G)}-\int_{\R^d}\F(x):\calQ\,\G(x)dx\bigg|
\,\lesssim_{\alpha,p}\,\e\mu_d(\tfrac1\e)^\frac12\\
\times\big(\|\F\|_{\Ld^2(\R^d)}+\|w_1^{2\alpha}D\F\|_{\Ld^{2p}(\R^d)}\big)\big(\|\G\|_{\Ld^2(\R^d)}+\|w_1^{2\alpha}D\G\|_{\Ld^{2p}(\R^d)}\big).
\end{multline}
Moreover, for all $1\le i,j,k,l\le d$ and all $\delta>0$, we have for all $L\ge1$,
\begin{align}\label{e.cov-struc.2}
\bigg|\calQ_{ijkl}-\int_{Q_{2L}}\frac{|Q_L\cap(x+Q_L)|}{|Q_L|}\cov{\Xi_{ij}(x)}{\Xi_{kl}(0)}dx\bigg|\,\lesssim_\delta\, L^{\delta-\frac12},
\end{align}
where $Q_L:=[\frac L2,\frac L2)^d$ denotes the cube of sidelength $L$ centered at the origin.
\smallskip\item If in addition $\Pm$ is nontrivial, then this \emph{effective fluctuation tensor} $\calQ$ is non-degenerate in the sense that $(\ee\otimes \ee):\calQ\,(\ee\otimes \ee)>0$ for all $\ee\in\R^d\setminus\{0\}$.\qedhere
\end{enumerate}
\end{prop}

\begin{samepage}
\begin{rems}\mbox{}
\begin{itemize}
\item When applying a covariance inequality (cf.\@ Lemma~\ref{lem:cov-inequ-improved} below) to the argument of the limit in the Green-Kubo formula~\eqref{e.cov-struc.2}, we end up with the bound 
$$
\int_{Q_{2L}}\frac{|Q_L\cap(x+Q_L)|}{|Q_L|}\,|\cov{\Xi_{ij}(x)}{\Xi_{kl}(0)}\!|\,dx\,\lesssim \,\log L,
$$
which is sharp. The main difficulty to characterize the limiting covariance structure is that, as usual for Green-Kubo formulas, the covariance of the homogenization commutator~$\Xi$ is not integrable and cancellations have to be unravelled.\smallskip
\item
The optimal rate~\eqref{eq:conv-cov-main-I0} for the convergence of the covariance structure of $J_0^\e$ owes to the very local structure of the homogenization commutator $\Xi$.
Combined with the pathwise result of Corollary~\ref{cor:pathwise}, it carries over to $I_1^\e$, $I_2^\e$, $I_3^\e$, $I_4^\e$, $J_1^\e$, and $J_2^\e$.
In~\cite{MO,GuM}, the Gaussian Helffer-Sj\"ostrand representation formula for the variance~\cite{HS-94,Sjostrand-96,NS2} was already used in order to prove the convergence of the covariance structure of $I_1^\e$ and $J_1^\e$ for $d>2$, but the obtained convergence rate was suboptimal in every dimension.
\smallskip
\item The non-degeneracy property (ii) already follows from \cite[Proposition~2.1]{GN} (modulo the identification~\eqref{eq:calQ-charact}); see also~\cite[Remark~2.3]{MO} in the Gaussian case.
\qedhere
\end{itemize}
\end{rems}
\end{samepage}

The combination of Propositions~\ref{prop:norm-dens} and~\ref{prop:conv-cov} leads to a complete scaling limit result for~$J_0^\e$, and proves the convergence in law to a Gaussian white noise.

\begin{cor}
Let $d\ge2$, let $\Pm$ be a product measure, and let $\mu_d$ be defined as in~\eqref{eq:def-mud}.
Let $\calQ$ be the $4$-tensor defined in Proposition~\ref{prop:conv-cov}(i), and let $\Gamma$ denote the $2$-tensor Gaussian white noise with covariance tensor $\calQ$, that is, $\Gamma$ is the Gaussian random linear functional with zero expectation $\expec{\Gamma(\F)}=0$ and with covariance structure $\cov{\Gamma(\F)}{\Gamma(\G)}=\int_{\R^d}\F:\calQ\,\G$ for all $\F,\G\in C^\infty_c(\R^d)^{d\times d}$.
Then
for all $\F\in C^\infty_c(\R^d)^{d\times d}$ the random variable $J_0^\e(\F)$ converges in law to $\Gamma(\F)$, and for $\int_{\R^d} \F:\calQ\,\F\neq 0$ there holds
\[\dWK{J_0^\e(\F)\!}{\!\Gamma(\F)} ~\lesssim_\F~\e\mu_d(\tfrac1\e).\qedhere\]
\end{cor}

\begin{samepage}
\begin{rems}\label{rem:convcov-corr}
\mbox{}\begin{itemize}
\item Combined with the pathwise result of Corollary~\ref{cor:pathwise}, this result leads to a proof of the joint convergence~\eqref{eq:pathwise-res-qual} and implies quantitative versions of the known scaling limit results for $I_1^\e$ and $J_1^\e$: For all $f,g \in C^\infty_c(\R^d)^d$ and all $F\in C^\infty_c(\R^d)^{d\times d}$ the random variables $I_1^\e(f,g)$ and $J_1^\e(\F)$ converge in law to $\Gamma(\bar\calP_Hf\otimes\bar\calP_H^*g)$ and $-\Gamma(\bar\calP_H^*\F)$, respectively. Moreover, whenever $\int_{\R^d} (\bar\calP_Hf\otimes\bar\calP_H^*g):\calQ\,(\bar\calP_Hf\otimes\bar\calP_H^*g) \neq 0$ and $\int_{\R^d}\bar\calP_H^*\F:\calQ\,\bar\calP_H^*\F\ne0$, we have
\begin{eqnarray*}
\dWK{I_1^\e(f,g)}{\Gamma(\bar\calP_Hf\otimes\bar\calP_H^*g)}&\lesssim_{f,g}&\e\mu_d(\tfrac1\e),\\
\dWK{J_1^\e(\F)}{-\Gamma(\bar\calP_H^*\F)}&\lesssim_{\F}&\e\mu_d(\tfrac1\e).
\end{eqnarray*}
This extends and unifies~\cite{GN,MO,MN,GuM}, and yields the first scaling limit results in the critical dimension $d=2$.
Convergence rates are new in any dimension and are expected to be optimal.

\smallskip\item \emph{SPDE representation for the scaling limit of the solution operator.}
The scaling limit result for $I_1^\e$ above indicates that $\e^{-\frac d2}\nabla_\e(u_\e-\expec{u_\e})$ (seen as a random linear functional) converges in law to the solution~$DU$ in $\R^d$ of
\[-D\cdot\bar\Aa DU=D\cdot(\Gamma_i D_i\bar u).\]
This justifies a posteriori the conclusion (although not the strategy) of the heuristics by Armstrong, Gu, and Mourrat~\cite{GuM,MourratFun} in dimensions $d\ge2$. (See also~\cite{Gu-16} for a rigorous treatment of the easier case of dimension $d=1$.)
\smallskip\item \emph{Scaling limit of the corrector.}
The scaling limit result for $J_1^\e$ above shows that the rescaled corrector field $\e^{-\frac d2}D\phi(\frac\cdot\e)$ (seen as a random linear functional) converges in law to $D(-D\cdot \bar\Aa D)^{-1} D \cdot \Gamma$,
that is, to the gradient of a variant of the so-called Gaussian free field. This variant involves both $\bar\Aa$ and $\calQ$. As pointed out in~\cite{MO}, it is easily checked in Fourier space that this variant does not coincide in general  with the standard Gaussian free field (unless the compatibility condition $\calQ_{ijkl}=\eta_{ik}\bar\Aa_{lj}$ is satisfied for some matrix $\eta$, which however does not hold in elementary examples, see e.g.~\cite[Section~3]{Gu-Mourrat-17} and~\cite[equation~(5.35)]{D-thesis}).
This variant of the Gaussian free field is studied in~\cite{Gu-Mourrat-17}, where it is shown to be Markovian only in the standard case.
In the critical dimension $d=2$, since the whole-space Gaussian free field is not well-defined (only its gradient is), this implies the non-existence of stationary correctors.
\smallskip\item \emph{Gu and Mourrat's observation.} With the above results at hand, we recover the observation by Gu and Mourrat~\cite{GuM} that the usual two-scale expansion~\eqref{eq:2-scale-usual} of $u_\e$ is not accurate in the CLT scaling.
The above indeed shows that the fluctuations of $\e^{-\frac d2}\int_{\R^d} g\cdot\nabla_\e(u_\e-\expec{u_\e})$ and of $\e^{-\frac d2} \int_{\R^d} g\cdot\nabla\phi_i(\frac{\cdot}{\e})\nabla_{\e,i} \bar u$ are asymptotically given by $\Gamma(\bar\calP_Hf\otimes\bar\calP_H^*g)$ and by $\Gamma(\bar\calP_H^*((\bar\calP_Hf)\otimes g))$, respectively, and therefore do not coincide in general.
\qedhere\end{itemize}
\end{rems}
\end{samepage}

\subsection{Approximation of the fluctuation tensor}
We finally turn to the representative volume element (RVE) method for approximation of $\calQ$ (cf.~Theorem~\ref{th:main-2}).
Indeed, the Green-Kubo formula~\eqref{e.cov-struc.2} for the fluctuation tensor $\calQ$ is of no practical use in applications since it requires to solve the corrector equation on the whole space and for every realization of the random coefficient field. It is therefore natural to seek a suitable RVE approximation.
It consists in introducing an artificial period $L>0$ and in considering an $L$-periodized coefficient field $\Aa_L$, typically given by a suitable periodization in law (cf.~\cite{EGMN-15}). In the present iid\@ setting, we simply define $\Aa_L(x+Ly):=\Aa(x)$ for all $y\in \Z^d$ and $x\in Q_L:=[-\frac L2,\frac L2)^d$. Note that the map $\Aa\mapsto\Aa_L$ on $\Omega$ pushes forward the measure $\Pm$ to a measure $\Pm_L$ concentrated on $L$-periodic coefficients, so that we may view $\Aa_L$ as an element of $\Omega_L=[\lambda,1]^{\B_L}$ equipped with the product measure $\Pm_L=\pi^{\otimes\B_L}$, where $\B_L:=\{(x,x+\ee_i):x\in Q_L\cap\Z^d,1\le i\le d\}$.
We then define the $L$-periodized corrector $\phi_{L,i}$ in the direction $\ee_i$ as the unique $L$-periodic solution in $Q_L\cap\Z^d$ of
\begin{align}\label{eq:corr-per-def}
-\nabla^*\cdot\Aa_L(\nabla\phi_{L,i}+\ee_i)=0,
\end{align}
satisfying $\sum_{z\in Q_L\cap\Z^d}\phi_{L,i}(z)=0$, and we set $\phi_L:=(\phi_{L,i})_{i=1}^d$ (which we implicitly extend as usual into a periodic piecewise constant map on $\R^d$). The spatial average of the flux,
\[\bar\Aa_{L}\ee_i:=\fint_{Q_L}\Aa_L(\nabla\phi_{L,i}+\ee_i),\]
is then an RVE approximation for the homogenized coefficient $\bar\Aa\ee_i=\expec{\Aa(\nabla\phi_i+\ee_i)}$.
The optimal numerical analysis of this approximation was originally performed in~\cite{GO1,GO2,GNO1}, where it was established that for all $L\ge2$ there holds
\begin{gather}\label{eq:res-RVE-GO-0}
|\var{\bar\Aa_{L}}\!|^\frac12\lesssim L^{-\frac d2},\qquad |\expec{\bar\Aa_{L}}-\bar\Aa|\lesssim L^{-d}\log^dL.
\end{gather}
In Theorem~\ref{th:main-2}, we claim that the fluctuation tensor $\calQ$ coincides with the limit of the rescaled variance of $\bar\Aa_L^*$. In addition, this characterization naturally leads to an RVE approximation $\calQ_{L,N}$ for $\calQ$, for which we obtain the optimal error estimate.

\begin{rems}
\mbox{}
\begin{itemize}
\item Definition~\eqref{eq:def-QL-per} for $\calQ_{L,N}$ is equivalent to
\begin{align}\label{eq:def-QL-per-bis}
\calQ_{L,N}=\frac{L^{d}}{N-1}\sum_{n=1}^N\Big(\fint_{Q_L}\Xi_{L,N}^{(n)}\Big)\otimes\Big(\fint_{Q_L}\Xi_{L,N}^{(n)}\Big),
\end{align}
where
\begin{align*}
\Xi_{L,N,i}^{(n)}:=\Aa_L^{(n)}(\nabla\phi_{L,i}^{(n)}+\ee_i)-\bar\Aa_{L,N}(\nabla\phi_{L,i}^{(n)}+\ee_i),
\end{align*}
with the obvious notation $\nabla\phi_L^{(n)}:=\nabla\phi_L(\Aa_L^{(n)})$.
Since by stationarity
\[\int_{Q_L}\cov{\Xi_{L,N}(x)}{\Xi_{L,N}(0)}dx=L^d\, \var{\fint_{Q_L}\Xi_{L,N}},\]
formula~\eqref{eq:def-QL-per-bis} is in the spirit of the Green-Kubo formula~\eqref{e.cov-struc.2}.

\smallskip\item
In~\eqref{eq:res-RVE-GO-0} the standard deviation $|\var{\bar\Aa_L}\!|^\frac12$ of the RVE approximation for $\bar\Aa$ is seen to be $O(L^{\frac d2})$ times larger than the systematic error $|\expec{\bar\Aa_L}-\bar\Aa|$ (up to a logarithmic correction).
In practice, we rather use $\bar\Aa_{L,N}$ as an approximation for $\bar\Aa$,
\begin{gather*}
|\var{\bar\Aa_{L,N}}\!|^\frac12\lesssim N^{-\frac12}L^{-\frac d2},\qquad |\expec{\bar\Aa_{L,N}}-\bar\Aa|\lesssim L^{-d}\log^dL,
\end{gather*}
since in the regime $N\simeq L^d$ the standard deviation becomes of the same order as the systematic error $O(L^{-d})$.
Combining this with the estimates in Theorem~\ref{th:main-2}, since $\calQ_{L,N}$ is extracted at no further cost than $\bar\Aa_{L,N}$ itself, we may infer that an RVE approximation for $\calQ$ with accuracy $O(L^{-\frac d2})$ is extracted at the same cost as an RVE approximation for $\bar\Aa$ with accuracy $O(L^{-d})$.

\smallskip\item
In~\cite{N,Nolen-16,GN} (see also~\cite{R,BSW}), the fluctuations of the RVE approximation $\bar\Aa_L$ for the homogenized coefficient $\bar\Aa$ was investigated. Combined with the characterization~\eqref{eq:calQ-charact} of the limit of the rescaled variance, the main result in~\cite{GN} takes on the following guise, for all $L\ge2$ and all $N\ge1$,
\[\quad\sup_{\ee\in\R^d\setminus\{0\}}\,\dWK{N^\frac12L^\frac d2\frac{\ee\cdot(\bar \Aa_{L,N}-\bar\Aa)\ee}{(\ee\otimes\ee:\calQ:\ee\otimes\ee)^\frac12}}\Nc\lesssim N^{-\frac12}L^{-\frac d2}\log^dL.\qedhere\]
\end{itemize}
\end{rems}

\section{Pathwise structure}\label{chap:decompsol}

Henceforth we place ourselves in the discrete setting of Section~\ref{sec:main-results-fluc}.
In the present section, we establish the pathwise result stated in Proposition~\ref{prop:pathwise}, that is, the main novelty of this contribution.
Similar estimates also lead to the CLT scaling result stated in Proposition~\ref{prop:bounded-I0}, and we further deduce
Corollary~\ref{cor:pathwise}.

\subsection{Structure of the proof and auxiliary results}
The main tool that we use to prove Propositions~\ref{prop:bounded-I0} and~\ref{prop:pathwise} is the following Poincar\'e inequality (or spectral gap estimate) in the probability space, which holds for any product measure $\Pm$ on $\Omega$ (see e.g.\@ \cite[Lemma~2.3]{GO1} for a proof).
Let us first fix some notation. Let $X=X(\Aa)$ be a random variable on $\Omega$, that is, a measurable function of $(a(b))_{b\in\B}$.
We choose an iid copy $\Aa'$ of $\Aa$,\footnote{Although we are then working on a product probability space $\Omega\times \Omega$, we use for simplicity the same notation~$\Pm$ (and $\E$) for the product probability measure (and expectation), that is, with respect to both $\Aa$ and $\Aa'$.} and for all $b\in\B$ we denote by $\Aa^b$ the random field that coincides with $\Aa$ on all edges $b'\neq b$ and with $\Aa'$ on edge $b$.
In particular, $\Aa$ and $\Aa^b$ have the same law.
{\color{black}We use the abbreviation $X^b=X(\Aa^b)$ and define the difference operator 
$\Delta_bX := X-X^b$, which we call the (Glauber or discrete) \emph{vertical derivative} at edge $b$.}

\begin{lem}[e.g.~\cite{GO1}]\label{lem:var}
Let $\Pm$ be a product measure. For all $X=X(\Aa)\in \Ld^2(\Omega)$ we have
\[\var{X}\,\leq \, \frac12 \,\expec{\sum_{b\in\B}|\Delta_b X|^2}.\qedhere\]
\end{lem}

Next to the corrector $\phi$, we need to recall the notion of flux corrector $\sigma$, which was recently introduced in~\cite{GNO-reg} in the continuum stochastic setting (see also~\cite[Lemma~4.4]{GuM} and~\cite[Proposition~III.2.2]{Nguyen-thesis} for its subsequent introduction in the discrete setting) and was crucially used in \cite{GNO-quant,GO4}. It allows to put the equation for the two-scale homogenization error in divergence form (cf.~\eqref{e.2s-wu}).
Let $\sigma=(\sigma_{ijk})_{i,j,k=1}^{d}$ be the $3$-tensor defined as the unique solution in $\Z^d$ of
\begin{equation} 
-\triangle\sigma_{ijk}\,:=\,-\nabla^*\cdot \nabla \sigma_{ijk}\,=\,\nabla_jq_{ik}-\nabla_kq_{ij},\label{si.5}
\end{equation}
with $\nabla \sigma$ stationary and having finite second moment, and with $\sigma(0)=0$, where $q_i$ denotes the flux of the corrector
\begin{equation}\label{si.6}
q_i=\Aa(\nabla\phi_i+\ee_i)-\bar\Aa\ee_i,\qquad q_{ij}:=(q_i)_j.
\end{equation}
Note that for all $i$ the $2$-tensor field $\sigma_i:=(\sigma_{ijk})_{j,k=1}^d$ is skew-symmetric, that is,
\begin{equation}\label{f.19}
\sigma_{ijk}=-\sigma_{ikj},
\end{equation}
and is shown to satisfy
\begin{equation}\label{f.20}
\nabla^* \cdot  \sigma_i \,:=\, \ee_j\nabla_k^* \sigma_{ijk} \,=\, q_i.
\end{equation}
Although considering a symmetric coefficient field, we use non-symmetric notation in view of the extension to the continuum setting, and we denote by $\phi^*$ and $\sigma^*$ the corrector and flux corrector associated with the pointwise transpose coefficient field $\Aa^*$. For symmetric coefficient fields, $\phi^*=\phi$ and $\sigma^*=\sigma$.

We now describe the string of arguments that leads to Proposition~\ref{prop:pathwise}.
We start with a suitable decomposition of the vertical derivative of $E^\e(f,g)$, which is key to the proof. Note that we rather consider a suitable version $E^\e_0(f,g)$ of $E^\e(f,g)$, which only coincides with $E^\e(f,g)$ up to some minor discretization error (in the continuum setting $\bar u_\e$ and $\bar v_\e$ would simply coincide with $\bar u(\e\cdot)$ and $\bar v(\e\cdot)$).
In the proofs, it is convenient to rescale all quantities down to scale $1$.

\begin{lem}\label{lem:decompI3eps}
For all $\e>0$ and all $f,g\in C^\infty_c(\R^d)^d$, setting $f_\e:=f(\e\cdot)$ and $g_\e:=g(\e\cdot)$, we denote by $\bar u_\e$ and $\bar v_\e$ the unique Lax-Milgram solutions in~$\R^d$ of
\begin{align}\label{e.def-utildeeps}
-\nabla^*\cdot\bar\Aa\nabla\bar u_\e=\nabla^*\cdot(\e f_\e),\qquad -\nabla^*\cdot\bar\Aa^*\nabla\bar v_\e=\nabla^*\cdot (\e g_\e),
\end{align}
and we define
\begin{multline}\label{eq:def-E0}
E^\e_0(f,g)\,:=\,\e^{\frac d2-1}\int_{\R^d}g_\e\cdot\big(\Aa \nabla (u_\e(\e\cdot))-\bar\Aa \nabla (u_\e(\e\cdot))-\expec{\Aa \nabla (u_\e(\e\cdot))-\bar\Aa \nabla(u_\e(\e\cdot))}\big)\\
-\e^{\frac d2-1}\int_{\R^d}g_\e\cdot\Xi_i \nabla_i\bar u_\e,
\end{multline}
as well as the two-scale expansion error $w_{f,\e}:=u_\e(\e\cdot)-(1+\phi_i\nabla_i)\bar u_\e$.
Then we have for all $b\in \B$,
\begin{multline}\label{eq:decomp-der-I0I1}
\Delta_bE_0^\e(f,g)
=\e^{\frac d2-1}\int_{\R^d}g_{\e,j}(\nabla\phi_j^*+\ee_j)\cdot\Delta_b\Aa(\nabla w^b_{f,\e}+\phi_i^b\nabla\nabla_i\bar u_\e)\\
+\e^{\frac d2-1}\int_{\R^d}\phi_j^*(\cdot+\ee_k)\nabla_kg_{\e,j}\ee_k\cdot\Delta_b\Aa\nabla (u_\e^b(\e\cdot))\\
-\e^{\frac d2-1}\int_{\R^d}\phi_j^*(\cdot+\ee_k)\nabla_k(g_{\e,j}\nabla_i\bar u_\e)\ee_k\cdot\Delta_b\Aa(\nabla\phi_i^b+\ee_i)\\
+\e^{\frac d2-1}\int_{\R^d}\nabla r_\e\cdot\Delta_b\Aa\nabla(u_\e^b(\e\cdot))-\e^{\frac d2-1}\int_{\R^d}\nabla R_{\e,i}\cdot\Delta_b\Aa(\nabla\phi_i^b+\ee_i),
\end{multline}
where the auxiliary fields $r_\e$ and $R_{\e}=(R_{\e,i})_{i=1}^d$ are the unique Lax-Milgram solutions in $\R^d$ of
\begin{eqnarray}
-\nabla^*\cdot\Aa^*\nabla r_\e&=&\nabla_l^*\big(\phi_j^*(\cdot+\ee_k)\Aa_{kl}\nabla_kg_{\e,j}+\sigma_{jkl}^*(\cdot-\ee_k)\nabla_k^*g_{\e,j}\big),\label{eq:aux-reps}\\
-\nabla^*\cdot\Aa^*\nabla R_{\e,i}&=&\nabla_l^*\big(\phi_j^*(\cdot+\ee_k)\Aa_{kl}\nabla_k(g_{\e,j}\nabla_i\bar u_\e)+\sigma_{jkl}^*(\cdot-\ee_k)\nabla_k^*(g_{\e,j}\nabla_i\bar u_\e)\big).\label{eq:aux-seps}
\end{eqnarray}
\end{lem}

By the spectral gap estimate of Lemma~\ref{lem:var}, the desired pathwise result~\eqref{eq:pathwise-I1-st} follows from a suitable estimate of the sum over $\B$ of the squares of the right-hand side terms in~\eqref{eq:decomp-der-I0I1}.
For that purpose, we make crucial use of the following moment bounds for the extended corrector $(\phi,\sigma)$ and its gradient. {\color{black}These bounds are a variation of~\cite{GO1,GO2,GNO1} and are the discrete versions of~\cite[Theorem~2]{GNO-reg}
and \cite[Theorem~2]{GNO-quant}, the proof of which extends to the discrete setting considered here. 
Similar bounds hold under the assumption of finite range of dependence, cf.~\cite{AKM2,AKM-book,GO4}.} 
\begin{lem}[\cite{GNO1,GNO-quant}]\label{lem:bd-sigma}
Let $d\ge2$, let $\Pm$ be a product measure, and let $\mu_d$ be defined in~\eqref{eq:def-mud}. For all $q<\infty$ and all $z\in \Z^d$ we have
\[\expec{|\phi(z)|^q}^\frac{1}{q}+\expec{|\sigma(z)|^q}^{\frac{1}{q}} \,\lesssim_q \, \mu_d(|z|)^\frac12,\]
and
\[\expec{|\nabla \phi(z)|^q}^\frac{1}{q}+\expec{|\nabla \sigma(z)|^q}^{\frac{1}{q}} \,\lesssim_q \, 1.\qedhere\]
\end{lem}

{\color{black}An additional key ingredient consists of large-scale regularity theory as originating in the work of Armstrong and Smart~\cite{AS} (see also the prior annealed regularity in~\cite{MaO}), which has recently been fully developed in~\cite{Armstrong-Mourrat-16,AKM1,AKM2,AKM-book} and in~\cite{GNO-reg,GO4}.
More precisely, we make use of the following large-scale weighted Calder\'on-Zygmund estimates for the random elliptic operator $-\nabla^*\cdot \Aa\nabla$ as established in~\cite{GNO-reg} (see also~\cite{Armstrong-Daniel-16} and \cite[Section~7]{AKM-book}).
A proof in the continuum setting was originally included in the first version of this article, see now~\cite[Corollary~5]{GNO-reg}; the adaptation to the discrete setting is straightforward since the argument is solely based on the energy and Caccioppoli estimates. Note that the following holds
under mere stationarity and ergodicity of the coefficients if the moment bounds on $r_*$ are replaced
by almost sure finiteness.}

\begin{lem}\cite[Corollary~5 \& Theorem~2]{GNO-reg}\label{lem:max-reg}
Let $d\ge1$ and let $\Pm$ be a product measure.
There exists a $\frac18$-Lipschitz stationary random field $r_*\ge1$ on $\R^d$ with $\expec{r_*^q}\lesssim_q1$ for all $q<\infty$, such that the following holds almost surely: For $\e>0$, $2\le p<\infty$, and $0\le \gamma<d(p-1)$,
for any (sufficiently fast) decaying scalar field $w$ and vector field $h$ related in $\R^d$ by
\[-\nabla^*\cdot \Aa \nabla w=\nabla^*\cdot h,\]
we have
\[\int_{\R^d}\big(1+\e(|x|+r_*(0))\big)^{\gamma}\Big(\fint_{B_*(x)}|\nabla w|^2\Big)^\frac p2dx\,\lesssim_{\gamma,p}\,\int_{\R^d}\big(1+\e(|x|+r_*(0))\big)^\gamma |h(x)|^{p}dx\]
with the short-hand notation $B_*(x):=B_{r_*(x)}(x)$.
\end{lem}

\subsection{Proof of Proposition~\ref{prop:bounded-I0}}
We focus on $J_0^\e$, while the proof is similar for~$J_1^\e$ and~$J_2^\e$.
Let $\F \in C^\infty_c(\R^d)^{d\times d}$, and set $\F_\e:=\F(\e\cdot)$.
We split the proof into two steps: we start by giving a suitable representation formula for the vertical derivative $\Delta_bJ_0^\e(\F)$, and then apply the spectral gap estimate.

\medskip
\step1 We prove the following representation formula for the vertical derivative $\Delta_bJ_0^\e(\F)$,
\begin{align}\label{eq:rep-form-DI0}
\Delta_bJ_0^\e(\F)=\e^{\frac d2}\int_{\R^d}\F_{\e,ij}\ee_j\cdot\Delta_b\Aa(\nabla\phi_i^b+\ee_i)+\e^{\frac d2}\int_{\R^d}\nabla s_{\e,i}\cdot\Delta_b\Aa(\nabla\phi_i^b+\ee_i),
\end{align}
where the auxiliary field $s_{\e}$ is the unique Lax-Milgram solution in $\R^d$ of
\begin{align}\label{eq:si-step}
-\nabla^*\cdot\Aa^*\nabla s_{\e,i}=\nabla^*\cdot\big(\F_{\e,ij}(\Aa-\bar\Aa)\ee_j\big).
\end{align}
By definition of the homogenization commutator,
\[\Delta_bJ_0^\e(\F)=\e^{\frac d2}\int_{\R^d}\F_{\e,ij}\ee_j\cdot\Delta_b\Aa(\nabla\phi_i^b+\ee_i)+\e^{\frac d2}\int_{\R^d}\F_{\e,ij}\ee_j\cdot(\Aa-\bar\Aa)\nabla\Delta_b\phi_i.\]
By definition~\eqref{eq:si-step} of $s_{\e,i}$, we find
\[\Delta_bJ_0^\e(\F)=\e^{\frac d2}\int_{\R^d}\F_{\e,ij}\ee_j\cdot\Delta_b\Aa(\nabla\phi_i^b+\ee_i)-\e^{\frac d2}\int_{\R^d}\nabla s_{\e,i}\cdot\Aa\nabla\Delta_b\phi_i.\]
Using then the vertical derivative of the corrector equation~\eqref{e.corr} in the form
\begin{align}\label{eq:corr-eq-der}
-\nabla^*\cdot\Aa\nabla\Delta_b\phi_i=\nabla^*\cdot\Delta_b\Aa(\nabla\phi_i^b+\ee_i),
\end{align}
the claim~\eqref{eq:rep-form-DI0} follows.

\medskip
\step2 Conclusion.\\
For $b\in\B$ we use the notation $b=(z_b,z_b+\ee_b)$. Inserting the representation formula~\eqref{eq:rep-form-DI0} in the spectral gap estimate of Lemma~\ref{lem:var}, and noting that $|\Delta_b \Aa(x)|\lesssim\mathds1_{Q(z_b)}(x)$, we obtain
\begin{multline*}
\var{J_0^\e(\F)}\,\lesssim\,\e^d\sum_{b\in\B}\expec{|\nabla\phi^b(z_b)+\Id|^2}\int_{Q(z_b)}|\F_{\e}|^2\\
+\e^d\,\expec{\sum_{b\in\B}|\nabla\phi^b(z_b)+\Id|^2\int_{Q(z_b)}|\nabla s_{\e}|^2},
\end{multline*}
and hence, by Lemma~\ref{lem:bd-sigma} in the form $\expec{|\nabla\phi^b|^2}=\expec{|\nabla\phi|^2}\lesssim1$,
\begin{align}\label{eq:first-bound-varI0}
\var{J_0^\e(\F)}\,\lesssim\,\e^d\|\F_\e\|_{\Ld^2(\R^d)}^2+\e^d\,\expec{\sum_{b\in\B}|\nabla\phi^b(z_b)+\Id|^2\int_{Q(z_b)}|\nabla s_{\e}|^2}.
\end{align}
It remains to estimate the last right-hand side term. Using equation~\eqref{eq:corr-eq-der} in the form $-\nabla^*\cdot\Aa^b\nabla(\phi^b-\phi)=\nabla^*\cdot(\Aa^b-\Aa)(\nabla\phi+\Id)$, an energy estimate yields
\begin{align}\label{eq:apriori-est-depert}
|\nabla(\phi^b-\phi)(z_b)|^2\le \int_{\R^d}|\nabla(\phi^b-\phi)|^2\lesssim \int_{\R^d}|\Aa^b-\Aa|^2|\nabla\phi+\Id|^2\lesssim |\nabla\phi(z_b)+\Id|^2,
\end{align}
so that $|\nabla\phi^b(z_b)+\Id|\lesssim |\nabla\phi(z_b)+\Id|$. Further estimating in~\eqref{eq:first-bound-varI0} integrals over unit cubes by integrals over balls at scale $r_*$ (cf.\@ Lemma~\ref{lem:max-reg}), smuggling in a power $\alpha\frac{p-1}p$ of the weight $w_\e(z):=1+\e|z|$, and applying H\"older's inequality in space with exponent $p$, we deduce for all $p>1$,
\begin{multline*}
\e^d\,\expec{\sum_{b\in\B}|\nabla\phi^b(z_b)+\Id|^2\int_{Q(z_b)}|\nabla s_{\e}|^2}\,\lesssim\,\e^d\,\expec{\int_{\R^d}|\nabla\phi(z)+\Id|^2\Big(\int_{Q_2(z)}|\nabla s_{\e}|^2\Big) dz}\\
\lesssim\,\e^d\,\E\bigg[\bigg(\int_{\R^d}|\nabla\phi(z)+\Id|^{\frac{2p}{p-1}}r_*(z)^{\frac{dp}{p-1}}w_\e(z)^{-\alpha}dz\bigg)^{\frac{p-1}p}\\
\times\bigg(\int_{\R^d}w_\e(z)^{\alpha(p-1)}\Big(\fint_{B_*(z)}|\nabla s_{\e}|^2\Big)^p dz\bigg)^\frac1p\bigg].
\end{multline*}
Applying H\"older's inequality in the probability space and Fubini's theorem, using Lemmas~\ref{lem:bd-sigma} and~\ref{lem:max-reg} in the form $\expecm{|\nabla\phi+\Id|^{q}+r_*^{q}}\lesssim_q1$ for all $q<\infty$, and noting that $\int_{\R^d}w_\e(z)^{-\alpha}dz\lesssim_\alpha\e^{-d}$ provided $\alpha>d$, we obtain for all $p>1$ and all $\alpha>d$,
\begin{multline}\label{eq:reduc-before-maxreg}
\e^d\,\expec{\sum_{b\in\B}|\nabla\phi^b(z_b)+\Id|^2\int_{Q(z_b)}|\nabla s_{\e}|^2}\\
\lesssim_{\alpha,p}\,\e^{\frac{d}p}\,\expec{\int_{\R^d}w_\e(z)^{\alpha(p-1)}\Big(\fint_{B_*(z)}|\nabla s_{\e}|^2\Big)^p dz}^\frac1p.
\end{multline}
By large-scale weighted Calder\'on-Zygmund theory (cf.\@ Lemma~\ref{lem:max-reg}) applied to equation~\eqref{eq:si-step} for $s_\e$ with $\alpha(p-1)<d(2p-1)$, and using again the moment bounds on $r_*$, we deduce for all $0<p-1\ll1$ and all $0<\alpha-d\ll1$,
\begin{eqnarray}
\e^d\,\expec{\sum_{b\in\B}|\nabla\phi^b(z_b)+\Id|^2\int_{Q(z_b)}|\nabla s_{\e}|^2}
&\lesssim_{\alpha,p}&\e^{\frac{d}p}\,\expec{r_*(0)^{\alpha(p-1)}\int_{\R^d}w_\e^{\alpha(p-1)}|\F_\e|^{2p}}^\frac1p\nonumber\\
&\lesssim_{\alpha,p}&\e^{\frac{d}p}\|w_\e^{\alpha\frac{p-1}{2p}}\F_\e\|_{\Ld^{2p}(\R^d)}^2.\label{eq:main-passage-first-lem}
\end{eqnarray}
Inserting this into~\eqref{eq:first-bound-varI0} and rescaling spatial integrals, we deduce for all $0<p-1\ll1$ and all $0<\alpha-d\ll1$,
\[\var{J_0^\e(\F)}\,\lesssim_{\alpha,p}\,\|\F\|_{\Ld^2(\R^d)}^2+\|w_1^{\alpha\frac{p-1}{2p}}\F\|_{\Ld^{2p}(\R^d)}^2.\]
Further using H\"older's inequality in the form
\[\|\F\|_{\Ld^2(\R^d)}\le\Big(\int_{\R^d} w_1^{-\alpha}\Big)^{\frac{p-1}{2p}}\Big(\int_{\R^d}w_1^{\alpha(p-1)}|\F|^{2p}\Big)^\frac1{2p}\lesssim_{\alpha,p}\|w_1^{\alpha\frac{p-1}{2p}}\F\|_{\Ld^{2p}(\R^d)},\]
the conclusion follows (after replacing the exponent $\alpha\frac{p-1}{2p}$ by $2\alpha$).

%%%%%%%%%%%%%%%%%%%%%%%%%%%%%%%%%%
%%%%%%%%%%%%%%%%%%%%%%%%%%%%%%%%%%

\subsection{Proof of Lemma~\ref{lem:decompI3eps}}
We split the proof into two steps. To simplify notation, in this proof (and only in this proof), we write $u:=u_\e(\e\cdot)$.
\nopagebreak

\smallskip
\step1 We prove the following representation formula for $\Delta_b((\Aa-\bar\Aa)\nabla u)$,
\begin{multline}\label{eq:der-rep-quasi-Eps}
\Delta_b\big(\ee_j\cdot(\Aa-\bar\Aa)\nabla u\big)
=(\nabla\phi_j^*+\ee_j)\cdot\Delta_b\Aa\nabla u^b
-\nabla_k^*\big(\phi_j^*(\cdot+\ee_k)\ee_k\cdot\Delta_b\Aa\nabla u^b\big)\\
-\nabla_k^*\big(\phi_j^*(\cdot+\ee_k)\ee_k\cdot\Aa\nabla\Delta_bu\big)-\nabla_k\big(\sigma_{jkl}^*(\cdot-\ee_k)\nabla_l\Delta_bu\big).
\end{multline}
In particular, replacing $x\mapsto u(x)$ by $x\mapsto \phi_i(x)+x_i$, this implies the following discrete version of~\eqref{eq:Xi-very-local},
\begin{multline}\label{eq:der-Xi-claim}
\Delta_b\Xi_{ij}=(\nabla\phi_j^*+\ee_j)\cdot\Delta_b\Aa(\nabla\phi_i^b+\ee_i)
-\nabla_k^*\big(\phi_j^*(\cdot+\ee_k)\ee_k\cdot\Delta_b\Aa(\nabla\phi_i^b+\ee_i)\big)\\
-\nabla_k^*\big(\phi_j^*(\cdot+\ee_k)\Aa_{kl}\nabla_l\Delta_b\phi_i\big)-\nabla_k\big(\sigma_{jkl}^*(\cdot-\ee_k)\nabla_l\Delta_b\phi_i\big).
\end{multline}
Using the definition~\eqref{f.20} of $\sigma_j^*$ in the form $(\Aa^*-\bar\Aa^*)\ee_j=-\Aa^*\nabla\phi_j^*+\nabla^*\cdot\sigma_j^*$, we find
\begin{eqnarray}
\Delta_b\big(\ee_j\cdot(\Aa-\bar\Aa)\nabla u\big)&=&\ee_j\cdot\Delta_b\Aa\nabla u^b+\ee_j\cdot(\Aa-\bar\Aa)\nabla\Delta_b u\nonumber\\
&=&\ee_j\cdot\Delta_b\Aa\nabla u^b+(\nabla^*\cdot\sigma_j^*)\cdot\nabla\Delta_b u-\nabla\phi_j^*\cdot\Aa\nabla\Delta_b u.\label{eq:first-der-Xi}
\end{eqnarray}
On the one hand, using the following discrete version of the Leibniz rule, for all $\chi_1,\chi_2:\Z^d \to \R$,
\begin{align}\label{pr:cov-str-2.4}
\nabla_l^*(\ee_l\chi_1(\cdot+\ee_l)\chi_2)=\chi_2\nabla\chi_1+\chi_1\nabla^*\chi_2,
\end{align}
we obtain
\[(\nabla^*\cdot\sigma_j^*)\cdot\nabla\Delta_bu\,=\,\nabla_l^*\big(\sigma_{jkl}^*\nabla_k\Delta_bu(\cdot+\ee_l)\big)-\sigma_{jkl}^*\nabla_k\nabla_l\Delta_bu,\]
so that the skew-symmetry~\eqref{f.19} of $\sigma_j$ leads to
\begin{eqnarray}\label{eq:der-Xi-01}
(\nabla^*\cdot\sigma_j^*)\cdot\nabla\Delta_bu\,=\,-\nabla_k^*\big(\sigma_{jkl}^*\nabla_l\Delta_bu(\cdot+\ee_k)\big)
\,=\,-\nabla_k\big(\sigma_{jkl}^*(\cdot-\ee_k)\nabla_l\Delta_bu\big).
\end{eqnarray}
On the other hand, using the vertical derivative of equation~\eqref{e.def-ueps} in the form $-\nabla^*\cdot\Aa\nabla\Delta_bu=\nabla^*\cdot\Delta_b\Aa\nabla u^b$, the discrete Leibniz rule~\eqref{pr:cov-str-2.4} yields
\begin{eqnarray}
\nabla\phi_j^*\cdot\Aa\nabla\Delta_bu&=&-\phi_j^*\nabla^*\cdot\Aa\nabla\Delta_bu+\nabla_k^*\big(\phi_j^*(\cdot+\ee_k)\ee_k\cdot\Aa\nabla\Delta_bu\big)\nonumber\\
&=&-\nabla\phi_j^*\cdot\Delta_b\Aa\nabla u^b+\nabla_k^*\big(\phi_j^*(\cdot+\ee_k)\ee_k\cdot\Delta_b\Aa\nabla u^b\big)\nonumber\\
&&\hspace{4.5cm}+\nabla_k^*\big(\phi_j^*(\cdot+\ee_k)\ee_k\cdot\Aa\nabla\Delta_bu\big).\label{eq:der-Xi-02}
\end{eqnarray}
Inserting~\eqref{eq:der-Xi-01} and~\eqref{eq:der-Xi-02} into~\eqref{eq:first-der-Xi},
the claim~\eqref{eq:der-rep-quasi-Eps} follows.

\medskip
\step2 Conclusion.\\
Integrating identities~\eqref{eq:der-rep-quasi-Eps} and~\eqref{eq:der-Xi-claim} with the test functions $g_\e$ and $\nabla\bar u_\e\otimes g_\e$, respectively, and integrating by parts, we obtain by definition of $E^\e_0$,
\begin{eqnarray*}
\Delta_bE_0^\e(f,g)
&=&\e^{\frac d2-1}\int_{\R^d}g_{\e,j}(\nabla\phi_j^*+\ee_j)\cdot\Delta_b\Aa\big(\nabla u^b-(\nabla\phi_i^b+\ee_i)\nabla_i\bar u_\e\big)\\
&&\hspace{-0.2cm}+\e^{\frac d2-1}\int_{\R^d}\phi_j^*(\cdot+\ee_k)\nabla_kg_{\e,j}\ee_k\cdot\Delta_b\Aa\nabla u^b\\
&&\hspace{-0.2cm}-\e^{\frac d2-1}\int_{\R^d}\phi_j^*(\cdot+\ee_k)\nabla_k(g_{\e,j}\nabla_i\bar u_\e)\ee_k\cdot\Delta_b\Aa(\nabla\phi_i^b+\ee_i)\\
&&\hspace{-0.2cm}+\e^{\frac d2-1}\int_{\R^d}\big(\phi_j^*(\cdot+\ee_k)\Aa_{kl}\nabla_kg_{\e,j}+\sigma_{jkl}^*(\cdot-\ee_k)\nabla_k^*g_{\e,j}\big)\nabla_l\Delta_bu\\
&&\hspace{-0.2cm}-\e^{\frac d2-1}\int_{\R^d}\big(\phi_j^*(\cdot+\ee_k)\Aa_{kl}\nabla_k(g_{\e,j}\nabla_i\bar u_\e)+\sigma_{jkl}^*(\cdot-\ee_k)\nabla_k^*(g_{\e,j}\nabla_j\bar u_\e)\big)\nabla_l\Delta_b\phi_i.
\end{eqnarray*}
The first right-hand side term is reformulated using the definition of $w_{f,\e}$ in the form $\nabla u^b-(\nabla\phi_i^b+\ee_i)\nabla_i\bar u_\e=\nabla w^b_{f,\e}+\phi_i^b\nabla\nabla_i\bar u_\e$.
It remains to post-process the last two right-hand side terms. Using equation~\eqref{eq:aux-reps} for $r_\e$ and using the vertical derivative of equation~\eqref{e.def-ueps} for $u_\e$ in the form $-\nabla^*\cdot\Aa\nabla\Delta_bu=\nabla^*\cdot\Delta_b\Aa\nabla u^b$, we find
\begin{multline*}
\int_{\R^d}\big(\phi_j^*(\cdot+\ee_k)\Aa_{kl}\nabla_kg_{\e,j}+\sigma_{jkl}^*(\cdot-\ee_k)\nabla_k^*g_{\e,j}\big)\nabla_l\Delta_bu\\
=-\int_{\R^d}\nabla r_\e\cdot\Aa\nabla\Delta_bu=\int_{\R^d}\nabla r_\e\cdot\Delta_b\Aa\nabla u^b.
\end{multline*}
Similarly, equations~\eqref{eq:aux-seps} and~\eqref{eq:corr-eq-der} lead to
\begin{multline*}
\int_{\R^d}\big(\phi_j^*(\cdot+\ee_k)\Aa_{kl}\nabla_k(g_{\e,j}\nabla_i\bar u_\e)+\sigma_{jkl}^*(\cdot-\ee_k)\nabla_k^*(g_{\e,j}\nabla_i\bar u_\e)\big)\nabla_l\Delta_b\phi_i\\
=-\int_{\R^d}\nabla R_{\e,i}\cdot\Aa\nabla\Delta_b\phi_i=\int_{\R^d}\nabla R_{\e,i}\cdot\Delta_b\Aa(\nabla\phi_i^b+\ee_i),
\end{multline*}
and the conclusion follows.

\subsection{Proof of Proposition~\ref{prop:pathwise}}

Using the representation formula~\eqref{eq:decomp-der-I0I1}, and recalling that for symmetric coefficients we have $(\phi^*,\sigma^*)=(\phi,\sigma)$,
the spectral gap estimate of Lemma~\ref{lem:var} leads to
\begin{align}\label{eq:bound-E0-pre-SG}
\var{E_0^\e(f,g)} \,\lesssim\, T_1^\e+T_2^\e+T_3^\e+T_4^\e+T_5^\e,
\end{align}
where we have set
\begingroup
\allowdisplaybreaks
\begin{eqnarray*}
T_1^\e&:=&\sum_{b\in\B}\expec{\Big(\e^{\frac d2-1}\int_{\R^d}g_{\e,j}(\nabla\phi_j+\ee_j)\cdot\Delta_b\Aa(\nabla w^b_{f,\e}+\phi_i^b\nabla\nabla_i\bar u_\e)\Big)^2},
\\
T_2^\e&:=&\sum_{b\in\B}\expec{\Big(\e^{\frac d2-1}\int_{\R^d}\phi_j(\cdot+\ee_k)\nabla_kg_{\e,j}\ee_k\cdot\Delta_b\Aa\nabla (u_\e^b(\e\cdot))\Big)^2},\\
T_3^\e&:=&\sum_{b\in\B}\expec{\Big(\e^{\frac d2-1}\int_{\R^d}\phi_j(\cdot+\ee_k)\nabla_k(g_{\e,j}\nabla_i\bar u_\e)\ee_k\cdot\Delta_b\Aa(\nabla\phi_i^b+\ee_i)\Big)^2},\\
T_4^\e&:=&\sum_{b\in\B}\expec{\Big(\e^{\frac d2-1}\int_{\R^d}\nabla r_\e\cdot\Delta_b\Aa\nabla(u_\e^b(\e\cdot))\Big)^2},\\
T_5^\e&:=&\sum_{b\in\B}\expec{\Big(\e^{\frac d2-1}\int_{\R^d}\nabla R_{\e,i}\cdot\Delta_b\Aa(\nabla\phi_i^b+\ee_i)\Big)^2},
\end{eqnarray*}
\endgroup
with the auxiliary fields $r_\e$ and $R_{\e}$ defined in~\eqref{eq:aux-reps} and in~\eqref{eq:aux-seps}.
The conclusion of Proposition~\ref{prop:pathwise} is a consequence of the following five estimates: for all $0<p-1\ll1$ and all $0<\alpha-d\ll1$,
\begin{eqnarray}
T_1^\e&\lesssim_{\alpha,p}&\e^{2}\mu_d(\tfrac1\e)\,\|g\|_{\Ld^4(\R^d)}^2\|w_1^{\alpha\frac{p-1}{4p}}\mu_d(|\cdot|)^\frac12Df\|_{\Ld^{4p}(\R^d)}^2,\label{eq:bound-T1}\\
T_2^\e&\lesssim&\e^{2}\mu_d(\tfrac1\e)\,
\|f\|_{\Ld^4(\R^d)}^2\|\mu_d(|\cdot|)^\frac12 Dg\|_{\Ld^4(\R^d)}^2,\label{eq:bound-T2}\\
T_3^\e&\lesssim&\e^{2}\mu_d(\tfrac1\e)\Big(\|f\|_{\Ld^4(\R^d)}^2\|\mu_d(|\cdot|)^\frac12Dg\|_{\Ld^4(\R^d)}^2\nonumber\\
&&\hspace{5cm}+\|g\|_{\Ld^4(\R^d)}^2\|\mu_d(|\cdot|)^\frac12Df\|_{\Ld^4(\R^d)}^2\Big),\label{eq:bound-T3}\\
T_4^\e&\lesssim_{\alpha,p}&\e^{2}\mu_d(\tfrac1\e)\,\|f\|_{\Ld^4(\R^d)}^2\|w_1^{\alpha\frac{p-1}{4p}}\mu_d(|\cdot|)^\frac12Dg\|_{\Ld^{4p}(\R^d)}^2,\label{eq:bound-T4}\\
T_5^\e&\lesssim_{\alpha,p}&\e^{2}\mu_d(\tfrac1\e)\Big(\|f\|_{\Ld^{4}(\R^d)}^2\|w_1^{\alpha\frac{p-1}{4p}}\mu_d(|\cdot|)^\frac12Dg\|_{\Ld^{4p}(\R^d)}^2\nonumber\\
&&\hspace{4cm}+\|g\|_{\Ld^{4}(\R^d)}^2\|w_1^{\alpha\frac{p-1}{4p}}\mu_d(|\cdot|)^\frac12 Df\|_{\Ld^{4p}(\R^d)}^2\Big).\label{eq:bound-T5}
\end{eqnarray}
We split the proof into three steps: we prove the above five estimates in the first two steps, and we conclude in the last step by controlling the discretization error.

\medskip
\step1
We establish the following equation for the two-scale expansion error $w_{f,\e}$ on $\R^d$,
\begin{align}\label{e.2s-wu}
-\nabla^*\cdot\Aa\nabla w_{f,\e}=\nabla_l^*\Big(\sigma_{jkl}(\cdot-\ee_k)\nabla_k^*\nabla_j\bar u_\e+\phi_j(\cdot+\ee_k)\Aa_{lk}\nabla_k\nabla_j\bar u_\e\Big),
\end{align}
which constitutes a discrete counterpart of similar identities in~\cite{GNO-reg,GNO-quant}.

\medskip\noindent
Using equations~\eqref{e.def-ueps} and~\eqref{e.def-utildeeps} in the form $-\nabla^*\cdot\Aa\nabla(u_\e(\e\cdot))=-\nabla^*\cdot\bar\Aa\nabla\bar u_\e$, and  using the following discrete version of the Leibniz rule, for all $\chi_1,\chi_2:\Z^d \to \R$,
\begin{align}\label{pr:cov-str-2.4-bis}
\nabla(\chi_1\chi_2)=\chi_1\nabla\chi_2+\ee_l\chi_2(\cdot+\ee_l)\nabla_l\chi_1,
\end{align}
we obtain
\begin{multline*}
-\nabla^*\cdot\Aa\nabla w_{f,\e}=-\nabla^*\cdot\Aa\nabla (u_\e(\e\cdot)-\bar u_\e-\phi_j\nabla_j\bar u_\e)\\
=-\nabla^*\cdot\bar\Aa\nabla\bar u_\e+\nabla^*\cdot\Aa\nabla\bar u_\e
+\nabla^*\cdot(\Aa\nabla\phi_j\nabla_j\bar u_\e)+\nabla^*\cdot(\Aa\ee_k\phi_j(\cdot+\ee_k)\nabla_k\nabla_j\bar u_\e).
\end{multline*}
Rearranging the terms and using the definition~\eqref{f.20} of $\sigma_j$, this turns into
\begin{eqnarray*}
-\nabla^*\cdot\Aa\nabla w_{f,\e}&=&\nabla^*\cdot\big((\Aa(\nabla\phi_j+\ee_j)-\bar\Aa\ee_j)\nabla_j\bar u_\e\big)+\nabla^*\cdot(\Aa\ee_k\phi_j(\cdot+\ee_k)\nabla_k\nabla_j\bar u_\e)\\
&=&\nabla^*\cdot\big((\nabla^*\cdot\sigma_j)\nabla_j\bar u_\e\big)+\nabla^*\cdot(\Aa\ee_k\phi_j(\cdot+\ee_k)\nabla_k\nabla_j\bar u_\e).
\end{eqnarray*}
Using again the discrete Leibniz rule~\eqref{pr:cov-str-2.4-bis} and the skew-symmetry~\eqref{f.19} of $\sigma_j$, we find
\begin{eqnarray*}
\nabla^*\cdot\big((\nabla^*\cdot\sigma_j)\nabla_j\bar u_\e\big)&=&\nabla_k^*(\nabla_l^*\sigma_{jkl}\nabla_j\bar u_\e)=\underbrace{\nabla_k^*\nabla_l^*\sigma_{jkl}\nabla_j\bar u_\e}_{=0}+\nabla_l^*\sigma_{jkl}(\cdot-\ee_k)\nabla_k^*\nabla_j\bar u_\e\\
&=&\nabla_l^*(\sigma_{jkl}(\cdot-\ee_k)\nabla_k^*\nabla_j\bar u_\e)-\underbrace{\sigma_{jkl}(\cdot-\ee_k-\ee_l)\nabla_k^*\nabla_l^*\nabla_j\bar u_\e}_{=0},
\end{eqnarray*}
and the conclusion~\eqref{e.2s-wu} follows.

\medskip
\step2 Proof of estimates~\eqref{eq:bound-T1}--\eqref{eq:bound-T5}.\\
We start with the first term $T_1^\e$. For $b\in\B$ we use the notation $b=(z_b,z_b+\ee_b)$. Since $|\Delta_b \Aa(x)|\lesssim\mathds1_{Q(z_b)}(x)$, the Cauchy-Schwarz inequality yields
\begin{eqnarray*}
T_1^\e&\lesssim&\e^{d-2}\sum_{b\in\B}\expec{|\nabla\phi(z_b)+\Id|^2\Big(\int_{Q(z_b)}|g_{\e}||\nabla w^b_{f,\e}+\phi_i^b\nabla\nabla_i\bar u_\e|\Big)^2}\\
&\lesssim&\e^{d-2}\,\expec{\sum_{b\in\B}|\nabla\phi(z_b)+\Id|^4\Big(\int_{Q(z_b)}|g_{\e}|^2\Big)^2}^\frac12\\
&&\hspace{5cm}\times\,\expec{\sum_{b\in\B}\Big(\int_{Q(z_b)}|\nabla w^b_{f,\e}+\phi_i^b\nabla\nabla_i\bar u_\e|^2\Big)^2}^\frac12,
\end{eqnarray*}
and hence, using the moment bounds of Lemma~\ref{lem:bd-sigma} and the exchangeability of $(\Aa,\Aa^b)$,
\begin{align*}
T_1^\e\lesssim\e^{d-2}\,\|g_\e\|_{\Ld^4(\R^d)}^2\,\expec{\sum_{b\in\B}\Big(\int_{Q(z_b)}|\nabla w_{f,\e}+\phi_i\nabla\nabla_i\bar u_\e|^2\Big)^2}^\frac12.
\end{align*}
We argue as in~\eqref{eq:reduc-before-maxreg}: We rewrite the second right-hand side factor as a norm of averages at the scale~$r_*$, smuggle in a suitable power of the weight $w_\e$, and apply H\"older's inequality, so that for all $p>1$ and all $\alpha>d$,
\begin{multline}\label{eq:treat-T1-Bstar}
T_1^\e\lesssim_{\alpha,p}\e^{\frac d2(1+\frac1p)-2}\,\|g_\e\|_{\Ld^4(\R^d)}^2\\
\times\,\expec{\int_{\R^d}w_\e(z)^{\alpha(p-1)}\Big(\fint_{B_*(z)}|\nabla w_{f,\e}|^2\Big)^{2p}dz+\int_{\R^d}w_\e^{\alpha(p-1)}|\phi|^{4p}|\nabla^2\bar u_\e|^{4p}}^\frac1{2p}.
\end{multline}
By large-scale weighted Calder\'on-Zygmund theory (cf.\@ Lemma~\ref{lem:max-reg})
applied to equation~\eqref{e.2s-wu} for $w_{f,\e}$, we deduce for all $0<p-1\ll1$ and all $0<\alpha-d\ll1$,
\begin{multline*}
T_1^\e\lesssim_{\alpha,p}\e^{\frac d2(1+\frac1p)-2}\,\|g_\e\|_{\Ld^4(\R^d)}^2\\
\times\,\expec{r_*(0)^{\alpha(p-1)}\int_{\R^d}w_\e^{\alpha(p-1)}\Big(|\sigma|^{4p}+|\phi|^{4p}+\sum_{k=1}^d|\phi(\cdot+\ee_k)|^{4p}\Big)|\nabla^2\bar u_\e|^{4p}}^\frac1{2p}.
\end{multline*}
By the moment bounds of Lemmas~\ref{lem:bd-sigma} and~\ref{lem:max-reg}, this yields
\begin{align*}
T_1^\e\lesssim_{\alpha,p}\e^{\frac d2(1+\frac1p)-2}\,\|g_\e\|_{\Ld^4(\R^d)}^2\|w_\e^{\alpha\frac{p-1}{4p}}\mu_d(|\cdot|)^{\frac12}\nabla^2\bar u_\e\|_{\Ld^{4p}(\R^d)}^2.
\end{align*}
We then apply  the standard weighted Calder\'on-Zygmund theory (e.g.~\cite[Section~V.4.2]{Stein-93}) to the discrete constant-coefficient equation~\eqref{e.def-utildeeps} for $\bar u_\e$ (cf.\@ Lemma~\ref{lem:max-reg} with $r_*=1$),
note that for all $\chi,\zeta\in C^\infty_c(\R^d)$ and all $q<\infty$
the inequality $|\nabla(\zeta(\e \cdot))|\le\e \int_0^1|D_k\zeta(\e(\cdot+ t\ee_k))|dt$ leads to
\begin{align}\label{eq:nabla-D}
\int_{\R^d}\chi|\nabla(\zeta(\e \cdot))|^q\le\e^q\int_{\R^d}\Big(\sup_{B(x)}|\chi|\Big)|D\zeta(\e x)|^qdx\le\e^{q-d}\int_{\R^d}\Big(\sup_{B(\frac x\e)}|\chi|\Big)|D\zeta(x)|^qdx,
\end{align}
rescale the integrals, estimate $\mu_d(|\frac\cdot\e|)\le \mu_d(\frac1\e)\mu_d(|\cdot|)$, and the conclusion~\eqref{eq:bound-T1} follows.

\medskip\noindent
We turn to the second term $T_2^\e$. Since $|\Delta_b \Aa(x)|\lesssim\mathds1_{Q(z_b)}(x)$, the Cauchy-Schwarz inequality yields
\begin{align*}
T_2^\e\lesssim
\e^{d-2}\,\expec{\sum_{b\in\B}|\phi(z_b+\ee_k)|^{2}\Big(\int_{Q(z_b)}|\nabla_kg_{\e}|^2\Big)\Big(\int_{Q(z_b)}|\nabla (u_\e^b(\e\cdot))|^2\Big)}.
\end{align*}
We bound the second local integral by an integral at the scale $r_*^b$, set $B_*^b(z):=B_{r_*^b(z)}(z)$, apply the Cauchy-Schwarz inequality, appeal to the moment bounds of Lemmas~\ref{lem:bd-sigma} and~\ref{lem:max-reg}, and obtain
\begin{align*}
T_2^\e\lesssim
\e^{d-2}\,\|\mu_d(|\cdot|)^\frac12 \nabla g_\e\|_{\Ld^4(\R^d)}^2\,\expec{\sum_{b\in\B}\Big(\fint_{B_*^b(z_b)\cup Q(z_b)}|\nabla (u_\e^b(\e\cdot))|^2\Big)^2}^\frac12,
\end{align*}
which, by exchangeability of $(\Aa,\Aa^b)$, takes the form
\begin{align*}
T_2^\e\lesssim
\e^{d-2}\,\|\mu_d(|\cdot|)^\frac12 \nabla g_\e\|_{\Ld^4(\R^d)}^2\,\expec{\int_{\R^d}\Big(\fint_{B_*(z)}|\nabla (u_\e(\e\cdot))|^2\Big)^2dz}^\frac12.
\end{align*}
By large-scale (unweighted) Calder\'on-Zygmund theory (cf.\@ Lemma~\ref{lem:max-reg})
applied to equation~\eqref{e.def-ueps} for $u_\e$ in the form $-\nabla^*\cdot\Aa\nabla(u_\e(\e\cdot))=\nabla^*\cdot(\e f_\e)$, we deduce
\begin{align*}
T_2^\e\lesssim
\e^{d}\,\|\mu_d(|\cdot|)^\frac12 \nabla g_\e\|_{\Ld^4(\R^d)}^2\,\|f_\e\|_{\Ld^4(\R^d)}^2,
\end{align*}
and the conclusion~\eqref{eq:bound-T2} follows similarly as above.

\medskip\noindent
The proof of~\eqref{eq:bound-T3} for $T_3^\e$ is more direct. Indeed, using the moment bounds of Lemma~\ref{lem:bd-sigma}, decomposing $\nabla(g_{\e,i}\nabla u_\e)=\nabla g_{\e,i}\otimes \nabla u_\e+g_{\e,i}(\cdot+\ee_k)\ee_k\otimes\nabla_k\nabla u_\e$, and applying the Cauchy-Schwarz inequality, we find
\begin{eqnarray*}
T_3^\e&\lesssim&\e^{d-2}\sum_{k=1}^d\expec{\sum_{b\in\B}|\phi(z_b+\ee_k)|^2|\nabla\phi^b(z_b)+\Id|^2\Big(\int_{Q(z_b)}|\nabla(g_{\e}\nabla\bar u_\e)|\Big)^2}\\
&\lesssim&\e^{d-2}\,\int_{\R^d}\mu_d(|\cdot|)|\nabla(g_{\e}\nabla\bar u_\e)|^2\\
&\lesssim&\e^{d-2}\Big(\|\nabla\bar u_\e\|_{\Ld^4(\R^d)}^2\|\mu_d(|\cdot|)\nabla g_\e\|_{\Ld^4(\R^d)}^2+\|g_\e\|_{\Ld^4(\R^d)}^2\|\mu_d(|\cdot|)\nabla^2\bar u_\e\|_{\Ld^4(\R^d)}^2\Big),
\end{eqnarray*}
and the conclusion~\eqref{eq:bound-T3} follows as above.

\medskip\noindent
We turn to the fourth term $T_4^\e$: We smuggle in a power $\alpha\frac{p-1}{2p}$ of the weight $w_\e$, apply H\"older's inequality with exponents $(\frac{2p}{p-1},2p,2)$, appeal to the moment bounds of Lemma~\ref{lem:max-reg}, use the exchangeability of $(\Aa,\Aa^b)$, and obtain for all $p>1$ and all $\alpha>d$,
\begin{eqnarray*}
T_4^\e&\lesssim&\e^{d-2}\,\expec{\sum_{b\in\B}r_*(z_b)^dr_*^b(z_b)^d\Big(\fint_{B_*(z_b)\cup Q(z_b)}|\nabla r_\e|^2\Big)\Big(\fint_{B_*^b(z_b)\cup Q(z_b)}|\nabla(u_\e^b(\e\cdot))|^2\Big)dz}\\
&\lesssim&\e^{\frac d2(1+\frac1p)-2}\,\expec{\int_{\R^d}w_\e(z)^{\alpha(p-1)}\Big(\fint_{B_*(z)}|\nabla r_\e|^2\Big)^{2p}dz}^\frac1{2p}\\
&&\hspace{7cm}\times\,\expec{\int_{\R^d}\Big(\fint_{B_*(z)}|\nabla(u_\e(\e\cdot))|^2\Big)^2dz}^\frac12.
\end{eqnarray*}
By the large-scale weighted Calder\'on-Zygmund theory (cf.\@ Lemma~\ref{lem:max-reg})
applied to equation~\eqref{eq:aux-reps} for $r_\e$ and to equation~\eqref{e.def-ueps} for $u_\e$, and the moment bounds of Lemma~\ref{lem:bd-sigma}, we deduce for all $0<p-1\ll1$ and all $0<\alpha-d\ll1$,
\begin{align*}
T_4^\e\lesssim_{\alpha,p}\e^{\frac d2(1+\frac1p)}\,\|w_\e^{\alpha\frac{p-1}{4p}}\mu_d(|\cdot|)^\frac12\nabla g_\e\|_{\Ld^{4p}(\R^d)}^2\|f_\e\|_{\Ld^4(\R^d)}^2,
\end{align*}
and the conclusion~\eqref{eq:bound-T4} follows as before.

\medskip\noindent
Finally, we turn to the last term $T_5^\e$: We use~\eqref{eq:apriori-est-depert} in form of $|\nabla\phi^b(z_b)+\Id|\lesssim|\nabla\phi(z_b)+\Id|$, smuggle in a power $\alpha\frac{p-1}{2p}$ of the weight $w_\e$, apply H\"older's inequality with exponents $(\frac{2p}{p-1},\frac{2p}{p+1})$, appeal to the moment bounds of Lemmas~\ref{lem:bd-sigma} and~\ref{lem:max-reg}, and therefore obtain for all $p>1$ and all $\alpha>d$,
\begin{eqnarray*}
T_5^\e&\lesssim&\e^{d-2}\,\expec{\int_{\R^d}|\nabla\phi(z)+\Id|^2\Big(\int_{Q_2(z)}|\nabla R_{\e}|^2\Big)dz}\\
&\lesssim_{\alpha,p}&\e^{\frac d2(1+\frac1p)-2}\,\expec{\int_{\R^d}w_\e(z)^{\alpha\frac{p-1}{p+1}}\Big(\fint_{B_*(z)}|\nabla R_{\e}|^2\Big)^\frac{2p}{p+1}dz}^\frac{p+1}{2p}.
\end{eqnarray*}
By the large-scale weighted Calder\'on-Zygmund theory (cf.\@ Lemma~\ref{lem:max-reg})
applied to equation~\eqref{eq:aux-seps} for $R_\e$, and the moment bounds of Lemma~\ref{lem:bd-sigma}, we deduce for all $0<p-1\ll1$ and all $0<\alpha-d\ll1$,
\begin{align*}
T_5^\e\lesssim_{\alpha,p}\e^{\frac d2(1+\frac1p)-2}\,\|w_\e^{\alpha\frac{p-1}{4p}}\mu_d(|\cdot|)^{\frac12}\nabla(g_\e\nabla\bar u_\e)\|_{\Ld^{\frac{4p}{p+1}}(\R^d)}^2.
\end{align*}
Decomposing $\nabla(g_{\e,i}\nabla u_\e)=\nabla g_{\e,i}\otimes \nabla u_\e+g_{\e,i}(\cdot+\ee_k)\ee_k\otimes\nabla_k\nabla u_\e$ and suitably applying H\"older's inequality with exponents $(\frac{p+1}{p},p+1)$, the conclusion~\eqref{eq:bound-T5} follows as before.

\medskip
\step3 Conclusion.\\
Inserting estimates~\eqref{eq:bound-T1}--\eqref{eq:bound-T5} into~\eqref{eq:bound-E0-pre-SG} yields for all $0<p-1\ll1$ and all $\alpha>d\frac{p-1}{4p}$,
\begin{multline}\label{eq:bound-E0-st}
\|E^\e_0(f,g)\|_{\Ld^2(\Omega)}\,\lesssim_{\alpha,p}\,\e\mu_d(\tfrac1\e)^\frac12\\
\times\Big(\|f\|_{\Ld^{4}(\R^d)}\|w_1^{\alpha}Dg\|_{\Ld^{4p}(\R^d)}
+\|g\|_{\Ld^{4}(\R^d)}\|w_1^{\alpha}Df\|_{\Ld^{4p}(\R^d)}\Big).
\end{multline}
It remains to deduce the corresponding result for $E^\e(f,g)$, and deal with the discretization error.
In terms of $\tilde u_\e:=\bar u_\e(\frac\cdot\e)$ and $\tilde v_\e:=\bar v_\e(\frac\cdot\e)$, equations~\eqref{e.def-utildeeps} take the form
\begin{align}\label{eq:def-ueps-tilde}
-\nabla_\e^*\cdot\bar\Aa\nabla_\e\tilde u_\e=\nabla_\e^*\cdot f,\qquad -\nabla_\e^*\cdot\bar\Aa^*\nabla_\e\tilde v_\e=\nabla_\e^*\cdot g.
\end{align}
The definitions of $E^\e$ and $E_0^\e$ then lead to the relation
\begin{align}\label{eq:ident-E0-dis}
E^\e(f,g)\,=\,E^\e_0(f,g)+J_0^\e\big((\nabla_\e\tilde u_\e-D\bar u)\otimes g\big),
\end{align}
where $J_0^\e((\nabla_\e\tilde u_\e-D\bar u)\otimes g)$ is a discretization error. By Proposition~\ref{prop:bounded-I0} (and Cauchy-Schwarz' inequality), it is enough to establish for all $1<p<\infty$ and all $0\le\alpha<d\frac{p-1}{p}$,
\begin{align}\label{eq:discr-Helmholtz}
\|w_1^\alpha(\nabla_\e\tilde u_\e-D\bar u)\|_{\Ld^{p}(\R^d)}\lesssim_{\alpha,p}\e\|w_1^\alpha Df\|_{\Ld^p(\R^d)}.
\end{align}
For that purpose, we note that $\bar u$ is an approximate solution of the discrete equation~\eqref{eq:def-ueps-tilde}. Indeed, integrating equation~\eqref{e.def-ubar} for $\bar u$ on a unit cube yields for all $x\in \R^d$
\begin{align}\label{eq:rewrite-discr-def-ubar}
0=\int_{[-1,0)^d}D\cdot(\bar\Aa D\bar u+f)(x+\e y) dy=\nabla_\e^*\cdot(\bar\Aa D\bar u+f)(x)+\nabla_\e^*\cdot T_{\e}(x),
\end{align}
where the error term $T_{\e}$ is given by $T_{\e}(x):=\ee_i\int_{S_i}((\bar\Aa D\bar u+f)_i(x+\e y)-(\bar\Aa D\bar u+f)_i(x))dy$ in terms of $S_i:=\{y\in[-1,0]^d:y_i=0\}$, and
satisfies for all $1\le p<\infty$ and all $0\le \alpha<\infty$,
\begin{align}\label{eq:est-Reps-err}
\|w_1^\alpha T_{\e}\|_{\Ld^p(\R^d)}\lesssim_\alpha\e\|w_1^\alpha D(\bar\Aa D\bar u+f)\|_{\Ld^p(\R^d)}\lesssim \e\|w_1^\alpha Df\|_{\Ld^p(\R^d)}+\e\|w_1^\alpha D^2\bar u\|_{\Ld^p(\R^d)}.
\end{align}
Comparing equations~\eqref{eq:def-ueps-tilde} and~\eqref{eq:rewrite-discr-def-ubar}, the difference $\bar u-\tilde u_\e$ satisfies
\[-\nabla_\e^*\cdot\bar\Aa\nabla_\e(\bar u-\tilde u_\e)=\nabla_\e^*\cdot T_\e-\nabla_\e^*\cdot\bar\Aa(\nabla_\e\bar u-D\bar u).\]
Hence, using the standard weighted Calder\'on-Zygmund theory (e.g.~\cite[Section~V.4.2]{Stein-93}) applied to this discrete constant-coefficient equation, we obtain for all $1<p<\infty$ and all $0\le\alpha<d\frac{p-1}{p}$,
\begin{align*}
\|w_1^\alpha\nabla_\e(\bar u-\tilde u_\e)\|_{\Ld^p(\R^d)}\,\lesssim_{\alpha,p}\,\|w_1^\alpha T_\e\|_{\Ld^p(\R^d)}+\|w_1^\alpha(\nabla_\e\bar u-D\bar u)\|_{\Ld^p(\R^d)}.
\end{align*}
Since the second right-hand side term is bounded by $\e\|w_1^\alpha D^2\bar u\|_{\Ld^p(\R^d)}$, estimate~\eqref{eq:est-Reps-err} yields
\begin{align*}
\|w_1^\alpha(\nabla_\e\tilde u_\e-D\bar u)\|_{\Ld^p(\R^d)}\,\lesssim_{\alpha,p}\,\e\|w_1^\alpha Df\|_{\Ld^p(\R^d)}+\e\|w_1^\alpha D^2\bar u\|_{\Ld^p(\R^d)}.
\end{align*}
The claim~\eqref{eq:discr-Helmholtz} then follows from the standard weighted Calder\'on-Zygmund theory applied to the constant-coefficient equation~\eqref{e.def-ubar} for $\bar u$.

\subsection{Proof of Corollary~\ref{cor:pathwise}}

We start with the proof of~\eqref{eq:pathwise-I1-st} for $I_1^\e$. By integration by parts, equations~\eqref{eq:def-ueps-tilde} and~\eqref{e.def-ueps} for $\tilde v_\e$, $\tilde u_\e$, and $u_\e$ lead to
\begin{eqnarray*}
\int g\cdot\nabla_\e(u_\e-\tilde u_\e)~\stackrel{\eqref{eq:def-ueps-tilde}}{=}~-\int \nabla_\e\tilde v_\e\cdot \bar\Aa\nabla_\e(u_\e-\tilde u_\e)&\stackrel{\eqref{eq:def-ueps-tilde}}{=}&-\int \nabla_\e\tilde v_\e\cdot f-\int \nabla_\e\tilde v_\e\cdot \bar\Aa\nabla_\e u_\e\\
&\stackrel{\eqref{e.def-ueps}}{=}&\int \nabla_\e\tilde v_\e\cdot (\Aa_\e\nabla_\e u_\e-\bar\Aa\nabla_\e u_\e).
\end{eqnarray*}
Subtracting the expectation of both sides yields a discrete version of identity~\eqref{eq:rel-I1}. In terms of $J_0^\e$, $I_1^\e$, and $E_0^\e$ (cf.\@ Section~\ref{chap:notations} and~\eqref{eq:def-E0}), this takes on the following guise,
\begin{align}\label{eq:ident-I1-dis}
I_1^\e(f,g)-J_0^\e(D\bar u\otimes D\bar v)
=J_0^\e(\nabla_\e\tilde u_\e\otimes \nabla_\e\tilde v_\e-D\bar u\otimes D\bar v)+E_0^\e(f,\nabla_\e\tilde v_\e).
\end{align}
Using~\eqref{eq:discr-Helmholtz} and the standard weighted Calder\'on-Zygmund theory (e.g.~\cite[Section~V.4.2]{Stein-93}) applied to the constant-coefficient equations~\eqref{e.def-ubar} and~\eqref{eq:def-ueps-tilde},
the conclusion~\eqref{eq:pathwise-I1-st} for $I_1^\e$ follows from~\eqref{eq:bound-E0-st} together with Proposition~\ref{prop:bounded-I0}.

\medskip\noindent
We turn to the proof of~\eqref{eq:pathwise-I1-st} for $I_2^\e$.
By definition of $J_0^\e$, $I_1^\e$, $I_2^\e$, and $E^\e$ (cf.\@ Section~\ref{chap:notations}), we find
\begin{align*}
I_2^\e(f,g)
=E^\e(f,g)+I_1^\e(f,\bar\Aa^*g)+J_0^\e(D\bar u\otimes g).
\end{align*}
Inserting identities~\eqref{eq:ident-E0-dis} and~\eqref{eq:ident-I1-dis} (with $g$ replaced by $\bar\Aa^*g$ and thus $\bar v$ replaced by the solution $\bar v^\circ$ of $-D\cdot\bar\Aa^*D\bar v^\circ=D\cdot\bar\Aa^*g$, so that $\bar\calP_L^*g=D\bar v^\circ+g$),
the conclusion~\eqref{eq:pathwise-I1-st} follows similarly as for $I_1^\e$.

\medskip\noindent
We now turn to the proof of~\eqref{eq:pathwise-I23-st}.
Let $\calS(\R^d)$ denote the Schwartz space of rapidly decaying functions, and consider the subspace $\calK_\e:=\{g\in \calS(\R^d)^d:\bar v_\e\in \calS(\R^d)\}$, cf.~\eqref{e.def-utildeeps}.
Given some fixed $\chi\in C^\infty_c(\R^d)$, set $\chi_L:=\chi(L\cdot)$ for $L\ge1$.
For $g\in\calK_\e$, we compute by integration by parts, using equation~\eqref{e.corr} for $\phi_j$ and equation~\eqref{e.def-utildeeps} for $\bar v_\e$, together with the discrete Leibniz rule~\eqref{pr:cov-str-2.4-bis},
\begin{eqnarray*}
\lefteqn{\int_{\R^d}\chi_L\nabla \bar v_\e\cdot\Xi_i
\,=\,\int_{\R^d}\chi_L\nabla\bar v_\e\cdot\big(\Aa(\nabla\phi_i+\ee_i)-\bar\Aa(\nabla\phi_i+\ee_i)\big)}\\
&\hspace{0.5cm}\stackrel{\eqref{e.corr}}{=}&-\int_{\R^d}\nabla(\bar v_\e \chi_L)\cdot\bar\Aa\nabla\phi_i-\int_{\R^d}\bar v_\e(\cdot+\ee_j)\nabla_{j}\chi_L\,\Xi_{ij}\\
&\hspace{0.5cm}\stackrel{\eqref{e.def-utildeeps}}{=}&\e\int_{\R^d} \chi_Lg_\e\cdot\nabla\phi_i+\e\int_{\R^d} \phi_i(\cdot+\ee_j) g_{\e,j}\nabla_j\chi_L+\int_{\R^d}\phi_i(\cdot+\ee_j)\nabla_j\chi_L\ee_j\cdot\bar\Aa\nabla\bar v_\e\\
&&\hspace{1.5cm}-\int_{\R^d}\bar v_\e(\cdot+\ee_j)\nabla_j\chi_L \ee_j\cdot\bar\Aa\nabla\phi_i-\int_{\R^d}\bar v_\e(\cdot+\ee_j)\nabla_{j}\chi_L\,\Xi_{ij}.
\end{eqnarray*}
For fixed $\e$ and $g\in\calK_\e$, using the moment bounds of Lemma~\ref{lem:bd-sigma} and the rapid decay at infinity of $g$ and $\bar v_\e$, we may pass to the limit $L\uparrow\infty$ in both sides in $\Ld^2(\Omega)$, and we deduce almost surely
\begin{align*}
\int_{\R^d}\nabla \bar v_\e:\Xi_j
\,=\,\e\int_{\R^d} g_\e\cdot\nabla\phi_j,
\end{align*}
that is, after rescaling,
\begin{align}\label{eq:equal-J1I0-pre}
J_1^\e(\ee_j\otimes g)=J_0^\e(\ee_j\otimes\nabla_\e\tilde v_\e).
\end{align}
We now argue that for all $\e>0$ this almost sure identity can be extended to hold in $\Ld^2(\Omega)$ for all $g\in C_c^\infty(\mathbb{R}^d)^d$ {\color{black}(even though $\nabla_\e\tilde v_\e$ is not integrable)}.
First note that Proposition~\ref{prop:bounded-I0} combined with the standard weighted Calder\'on-Zygmund theory (e.g.~\cite[Section~V.4.2]{Stein-93}) for the constant-coefficient equation~\eqref{eq:def-ueps-tilde} yields for all $0<p-1\ll1$ and all $d\frac{p-1}{4p}<\alpha<d\frac{2p-1}{4p}$,
\[\expec{|J_0^\e(\ee_j\otimes\nabla_\e\tilde v_\e)|^2}^\frac12\,\lesssim_{\alpha,p}\,\|w_1^{2\alpha}\nabla_\e\tilde v_\e\|_{\Ld^{2p}(\R^d)}\,\lesssim_{\alpha,p}\,\|w_1^{2\alpha}g\|_{\Ld^{2p}(\R^d)},\]
and in addition,
\[\expec{|J_1^\e(\ee_j\otimes g)|^2}^\frac12\,\lesssim\,\|w_1^{2\alpha}g\|_{\Ld^{2p}(\R^d)}.\]
Hence, it suffices to check the following density result: for all test functions $g\in C_c^\infty(\mathbb{R}^d)^d$ there exist a sequence $(g_n)_n$ of elements of $\calK_\e$ such that $\|w_1^{2\alpha}(g_n-g)\|_{\Ld^{2p}(\R^d)}\to 0$ holds for some $0<p-1\ll1$ and some $\alpha>d\frac{p-1}{4p}$.
Let $g\in C^\infty_c(\R^d)^{d}$ be fixed. Up to a convolution argument on large scales, we may already assume that the Fourier transform $\hat g$ has compact support, say contained in $B_R$.
Since the (continuum) Fourier symbol of the discrete Helmholtz projection 
$\nabla_\e(\nabla_\e^*\cdot\bar a\nabla_\e)^{-1}\nabla_\e^*\cdot$ is bounded and smooth outside of the dual lattice $(\frac{2\pi}{\e}\mathbb{Z})^d$, a function $g_n\in \calS(\mathbb{R}^d)^d$ actually belongs to ${\mathcal K}_\e$ whenever its Fourier transform $\hat g_n$ vanishes in a neighborhood of $(\frac{2\pi}{\e}\mathbb{Z})^d$.
Choosing $\chi\in C^\infty_c(\R^d)$ with $\chi=1$ in a neighborhood of $0$, and defining
\[\chi_n\,:=\,1-\sum_{z\in(\frac{2\pi}\e\Z)^d}\chi(n(\cdot-z)),\]
the function $g_n\in \calS(\R^d)^d$ defined by $\hat g_n:=\chi_n \hat g$ thus belongs to $\calK_\e$.
{\color{black}For $p\ge1$, setting $q:=\frac{2p}{2p-1}$, the Hausdorff-Young inequality yields
\begin{multline*}
\|w_1^{2\alpha}(g_n-g)\|_{\Ld^{2p}(\R^d)}\,\lesssim_p\,\|\mathcal F(w_1^{2\alpha}(g_n-g))\|_{\Ld^{q}(\R^d)}\\
\,\lesssim\,\|\hat g_n-\hat g\|_{W^{2\alpha,q}(\R^d)}\,=\,\|(\chi_n-1)\hat g\|_{W^{2\alpha,q}(\R^d)},
\end{multline*}
where $\mathcal F(h)=\hat h$ denotes the Fourier transform.
Hence, since $\hat g$ is supported in~$B_R$,}
\begin{multline*}
\|w_1^{2\alpha}(g_n-g)\|_{\Ld^{2p}(\R^d)}\,\lesssim_\alpha\,\|\chi_n-1\|_{W^{2\alpha,q}(B_R)}\|\hat g\|_{W^{2\alpha,\infty}(\R^d)}\\
\lesssim\,\|\chi_n-1\|_{W^{2\alpha,q}(B_R)}\|w_1^{2\alpha}g\|_{\Ld^1(\R^d)}.
\end{multline*}
For $2\alpha<\frac dq=d\frac{2p-1}{2p}$, reflecting the fact that the Sobolev space $W^{2\alpha,q}(\R^d)$ fails to embed into the space of continuous functions, there holds $\chi_n\to1$ in $W^{2\alpha,q}_{\operatorname{loc}}(\R^d)$ as $n\uparrow\infty$, and hence $\|w_1^{2\alpha}(g_n-g)\|_{\Ld^{2p}(\R^d)}\to 0$.
This establishes the claimed density result, 
and we conclude that identity~\eqref{eq:equal-J1I0-pre} can be extended in $\Ld^2(\Omega)$ to all $g\in C^\infty_c(\R^d)^d$. 
The estimate~\eqref{eq:pathwise-I23-st} for $J_1^\e$ then follows from the discretization error estimate~\eqref{eq:discr-Helmholtz} together with Proposition~\ref{prop:bounded-I0}. The estimate~\eqref{eq:pathwise-I23-st} for $J_2^\e$ is obtained in a similar way.

\section{Asymptotic normality}\label{sec:approx-normal-I0}

We turn to the normal approximation result for the homogenization commutator $\Xi$ as stated in Proposition~\ref{prop:norm-dens}.

\subsection{Structure of the proof and auxiliary results}

The main tool to prove Proposition~\ref{prop:norm-dens} is the following suitable form of a second-order Poincar\'e inequality \`a la Chatterjee~\cite{C1,LRP-15}. Based on Stein's method, it can be shown to hold 
for any product measure $\Pm$ on $\Omega$. (The proof follows from~\cite[Theorem~2.2]{C1} and~\cite[Theorem~4.2]{LRP-15} in the case of the Wasserstein and of the Kolmogorov metric, respectively, combined with the spectral gap estimate of Lemma~\ref{lem:var}.)
The first use of such functional inequalities in stochastic homogenization is due to Nolen~\cite{N,Nolen-16}.
Let us first fix some more notation. Let $X=X(\Aa)$ be a random variable on $\Omega$, that is, a measurable function of $(a(b))_{b\in\B}$. For all $E\subset\B$ we denote by $\Aa^E$ the random field that coincides with $\Aa$ on all edges $b\notin E$ and with the iid copy $\Aa'$ on all edges $b\in E$. In particular, $\Aa$ and $\Aa^E$ always have the same law.
We use the abbreviation $X^E:=X(\Aa^E)$ and define $\Delta_bX^E:=X^E-X^{E\cup\{b\}}$.
As before, we write for simplicity $X^b:=X^{\{b\}}$, and similarly $X^{b,b'}:=X^{\{b,b'\}}$. In particular, $\Delta_b\Delta_{b'}X=X-X^b-X^{b'}+X^{b,b'}$.

\begin{lem}[\cite{C1,LRP-15}]\label{lem:chat0}
Let $\Pm$ be a product measure and let $\delta_\calN$ be defined in~\eqref{eq:def-Delta-N}.
For all $X=X(\Aa)\in \Ld^2(\Omega)$, we have
\begin{align*}
&\delta_\calN(X)
\,\lesssim\,\frac1{\var{X}^\frac 32}\sum_{b\in\B_L}\expec{|\Delta_bX|^6}^\frac12\\
&\hspace{4cm}+\frac1{\var{X}}\bigg(\sum_{b\in\B}\Big(\sum_{e'\in\B}\expec{|\Delta_{b'}X|^4}^\frac14\expec{|\Delta_b\Delta_{b'}X|^4}^\frac14\Big)^2\bigg)^\frac12.\qedhere
\end{align*}
\end{lem}

In addition, we make crucial use of the following optimal annealed estimate on the mixed gradient of the Green's function, first proved by Marahrens and the third author~\cite{MaO}. {\color{black}Note that the product space assumption can be substantially relaxed (cf.~\cite{Bella-Giunti-18} and~\cite[Section~8.5]{AKM-book}).}
\begin{lem}[\cite{MaO}]\label{lem:ann-Green}
Let $d\ge2$ and let $\Pm$ be a product measure. For all $y\in\Z^d$ there exists a function $\nabla G(\cdot,y)$ that is the unique decaying solution in $\Z^d$ of
\[-\nabla^*\cdot\Aa\nabla G(\cdot,y)=\delta(\cdot-y).\]
It satisfies the following moment bound: for all $q<\infty$ and all $x,y\in \Z^d$,
\begin{gather*}
\expec{|\nabla\nabla G(x,y)|^q}^{\frac{1}{q}} \,\lesssim_q\, (1+|x-y|)^{-d},
\end{gather*}
where $\nabla\nabla$ denotes the mixed second gradient.
\end{lem}

%%%%%%%%%%%%%%%%%%%

\subsection{Proof of Proposition~\ref{prop:norm-dens}}

Let $\F \in C^\infty_c(\R^d)^{d\times d}$, and set $\F_\e:=\F(\e\cdot)$.
By Lemma~\ref{lem:chat0}, it is enough to estimate the following two contributions,
\begin{gather*}
K_1^\e:=\,\sum_{b\in\B}\expec{|\Delta_b J_0^{\e}(\F)|^6}^\frac12,\quad
K_2^\e:=\sum_{b\in\B}\left(\sum_{b'\in\B}\expec{|\Delta_{b'} J_0^{\e}(\F)|^{4}}^{\frac{1}{4}}\expec{|\Delta_b\Delta_{b'} J_0^{\e}(\F)|^{4}}^{\frac{1}{4}}\right)^2.
\end{gather*}
We split the proof into three steps: we start with an auxiliary estimate, and then estimate $K_1^\e$ and $K_2^\e$ separately.

\medskip
\step1 Auxiliary estimate: for all $\zeta\in C^\infty_c(\R^d)$, all $1\le p<\infty$, and all $r\ge0$,
\begin{align}\label{eq:auxiliary-result-CLT}
\int_{\R^d}\log^r(2+|z|)\bigg(\int_{\R^d}\frac{|\zeta(x)|}{(1+|x-z|)^d}dx\bigg)^pdz~\lesssim_{p,r}~\int_{\R^d}\log^{p+r}(2+|x|)\,|\zeta(x)|^p\,dx.
\end{align}
Let $\alpha>0$ be fixed. Smuggling in a power $\alpha\frac{p-1}{p}$ of the weight $1+|x|$, and applying H\"older's inequality with exponent $p$, we find
\begin{multline*}
\int_{\R^d}\log^r(2+|z|)\bigg(\int_{\R^d}\frac{|\zeta(x)|}{(1+|x-z|)^d}dx\bigg)^pdz\\
\le\int_{\R^d}\log^r(2+|z|)\bigg(\int_{\R^d}\frac{(1+|x|)^{\alpha(p-1)}|\zeta(x)|^p}{(1+|x-z|)^d}dx\bigg)\bigg(\int_{\R^d}\frac{dx}{(1+|x-z|)^d(1+|x|)^{\alpha}}\bigg)^{p-1}dz.
\end{multline*}
The last integral is controlled by $C(d,\alpha)\frac{\log(2+|z|)}{(1+|z|)^{\alpha}}$, hence by Fubini's theorem
\begin{eqnarray*}
\lefteqn{\int_{\R^d}\log^r(2+|z|)\bigg(\int_{\R^d}\frac{|\zeta(x)|}{(1+|x-z|)^d}dx\bigg)^pdz}\\
&\lesssim_{\alpha,p}&\int_{\R^d}\frac{\log^{p+r-1}(2+|z|)}{(1+|z|)^{\alpha(p-1)}}\bigg(\int_{\R^d}\frac{(1+|x|)^{\alpha(p-1)}|\zeta(x)|^p}{(1+|x-z|)^d}dx\bigg)dz\\
&=&\int_{\R^d}(1+|x|)^{\alpha(p-1)}|\zeta(x)|^p \bigg(\int_{\R^d}\frac{\log^{p+r-1}(2+|z|)}{(1+|x-z|)^d(1+|z|)^{\alpha(p-1)}}dz\bigg)dx.
\end{eqnarray*}
Since the last integral is controlled by $C(d,p,r,\alpha)\frac{\log^{p+r}(2+|x|)}{(1+|x|)^{\alpha(p-1)}}$, the conclusion~\eqref{eq:auxiliary-result-CLT} follows.

\medskip
\step{2} Proof of 
$$
K_1^\e\,\lesssim\, \e^{\frac{d}{2}}\big(\|\F\|_{\Ld^3(\R^d)}^3+\|\log(2+|\cdot|)\,\mu_d(|\cdot|)^\frac12D\F\|_{\Ld^3(\R^d)}^3\big).
$$
After integration by parts, the representation formula for the vertical derivative $\Delta_b\Xi$ in~\eqref{eq:der-Xi-claim} leads to
\begin{multline*}
\Delta_bJ_0^\e(\F)=\e^{\frac d2}\int_{\R^d}\F_{\e,ij}(\nabla\phi_j^*+\ee_j)\cdot\Delta_b\Aa(\nabla\phi_i^b+\ee_i)\\
+\e^{\frac d2}\int_{\R^d}\phi_j^*(\cdot+\ee_k)\nabla_k\F_{\e,ij}\ee_k\cdot\Delta_b\Aa(\nabla\phi_i^b+\ee_i)\\
+\e^{\frac d2}\int_{\R^d}\big(\phi_j^*(\cdot+\ee_k)\Aa_{kl}\nabla_k\F_{\e,ij}+\sigma_{jkl}^*(\cdot-\ee_k)\nabla_k^*\F_{\e,ij}\big)\nabla_l\Delta_b\phi_i.
\end{multline*}
For $b=(z_b,z_b+\ee_b)$, the Green representation formula applied to equation~\eqref{eq:corr-eq-der} takes the following form, for all $x\in\Z^d$,
\begin{align}\label{eq:formvertdernablaphi}
\nabla \Delta_b \phi_i(x)\,=\,-\nabla \nabla G(x,z_b) \Delta_b\Aa(z_b)(\nabla \phi_i^b(z_b)+\ee_i).
\end{align}
Inserted into the above representation formula for $\Delta_bJ_0^\e(\F)$, and combined with $|\Delta_b\Aa(x)|\lesssim\mathds1_{Q(z_b)}(x)$ and the moment bounds of Lemmas~\ref{lem:bd-sigma} and~\ref{lem:ann-Green}, it yields for all $q<\infty$,
\begin{align}\label{eq:bound-pre-K1eps}
\expec{|\Delta_{b} J_0^{\e}(\F)|^{q}}^{\frac{1}{q}}\,\lesssim_q\, \e^\frac d2\int_{Q(z_{b})}|\F_\e|+\e^\frac d2\int_{\R^d}\frac{\mu_d(|x|)^\frac12|\nabla\F_\e(x)|}{(1+|x-z_{b}|)^d}dx.
\end{align}
Summing the cube of this estimate over $b \in \B$ for $q=6$, and using~\eqref{eq:auxiliary-result-CLT} with $p=3$, $r=0$, and $\zeta=\mu_d(|\cdot|)^\frac12\nabla\F_\e$, we obtain
\begin{eqnarray*}
K_1^\e&\lesssim&\e^{\frac{3d}{2}} \int_{\R^d} |\F_\e|^3
+\e^{\frac{3d}{2}} \int_{\R^d} \bigg(\int_{\R^d}
\frac{\mu_d (|x|)^\frac{1}{2}|\nabla \F_\e(x)| }{(1+|x-z|^d)}dx\bigg)^3dz\\
&\lesssim& \e^{\frac{3d}{2}} \int_{\R^d} |\F_\e|^3
+\e^{\frac{3d}{2}} \int_{\R^d} \log^3(2+|\cdot|) \,\mu_d (|\cdot|)^\frac{3}{2}|\nabla \F_\e|^3.
\end{eqnarray*}
Rescaling the integrals and using~\eqref{eq:nabla-D}, the conclusion follows.

\medskip
\step{3} Proof of 
\begin{multline*}
K_2^\e\lesssim\e^{d}\log^2(2+\tfrac1\e)\Big(\|\F\|_{\Ld^4(\R^d)}^4+\|\log(2+|\cdot|)\F\|_{\Ld^4(\R^d)}^4\Big)\\
+\e^{d+2}\log^4(2+\tfrac1\e)\mu_d(\tfrac1\e)\Big(\|\log(2+|\cdot|)\F\|_{\Ld^4(\R^d)}^4+\|\log^2(2+|\cdot|)\mu_d(|\cdot|)^\frac12D\F\|_{\Ld^4(\R^d)}^4\Big)\\
+\e^{d+4}\log^6(2+\tfrac1\e)\mu_d(\tfrac1\e)^2\|\log^2(2+|\cdot|)\mu_d(|\cdot|)^\frac12D\F\|_{\Ld^4(\R^d)}^4.
\end{multline*}
We need to iterate the vertical derivative and estimate $\Delta_{b} \Delta_{b'} I_0^{\e}(\F)$.
By definition of the homogenization commutator, we find
\begin{eqnarray}
\lefteqn{\Delta_{b}\Delta_{b'}J_0^\e(\F)=\e^{\frac d2}\Delta_{b}\int_{\R^d}\F_{\e,ij}\ee_j\cdot\Delta_{b'}\Aa(\nabla\phi_i^{b'}+\ee_i)+\e^{\frac d2}\Delta_{b}\int_{\R^d}\F_{\e,ij}\ee_j\cdot(\Aa-\bar\Aa)\nabla\Delta_{b'}\phi_i}\nonumber\\
&=&\e^{\frac d2}\int_{\R^d}\F_{\e,ij}\ee_j\cdot\Delta_b\Delta_{b'}\Aa(\nabla\phi^{b,b'}_i+\ee_i)
+\e^{\frac d2}\int_{\R^d}\F_{\e,ij}\ee_j\cdot\big(\Delta_b\Aa\nabla\Delta_{b'}\phi_i^b+\Delta_{b'}\Aa\nabla\Delta_{b}\phi_i^{b'}\big)\nonumber\\
&&\hspace{6cm}
+\e^{\frac d2}\int_{\R^d}\F_{\e,ij}\ee_j\cdot(\Aa-\bar\Aa)\nabla\Delta_b\Delta_{b'}\phi_i.\label{eq:pre-form-der3-I0}
\end{eqnarray}
In order to avoid additional logarithmic factors, we need to suitably rewrite the last right-hand side term, and we argue similarly as in the proof of~\eqref{eq:der-Xi-claim}.
Using the definition~\eqref{f.20} of $\sigma_j^*$ in the form $(\Aa^*-\bar\Aa^*)\ee_j=-\Aa^*\nabla\phi_j^*+\nabla^*\cdot\sigma_j^*$, applying the discrete Leibniz rule~\eqref{pr:cov-str-2.4}, and using the skew-symmetry~\eqref{f.19} of $\sigma_i$, we obtain
\begin{multline*}
\ee_j\cdot(\Aa-\bar\Aa)\nabla\Delta_b\Delta_{b'}\phi_i=(\nabla^*\cdot\sigma_j^*)\cdot\nabla\Delta_b\Delta_{b'}\phi_i-\nabla\phi_j^*\cdot\Aa\nabla\Delta_b\Delta_{b'}\phi_i\\
=-\nabla_k\big(\sigma_{jkl}^*(\cdot-\ee_k)\nabla_l\Delta_b\Delta_{b'}\phi_i\big)
-\nabla_k^*\big(\phi_j^*(\cdot+\ee_k)\ee_k\cdot\Aa\nabla\Delta_b\Delta_{b'}\phi_i\big)
+\phi_j^*\nabla^*\cdot\Aa\nabla\Delta_b\Delta_{b'}\phi_i.
\end{multline*}
The vertical derivative of equation~\eqref{eq:corr-eq-der} takes the form 
\begin{align}\label{eq:corr-eq-der2}
-\nabla^*\cdot\Aa\nabla\Delta_b\Delta_{b'}\phi_i
=\nabla^*\cdot\Delta_{b'}\Aa\nabla\Delta_b\phi_i^{b'}+\nabla^*\cdot\Delta_b\Delta_{b'}\Aa(\nabla\phi_i^{b,b'}+\ee_i)+\nabla^*\cdot\Delta_b\Aa\nabla\Delta_{b'}\phi_i^b,
\end{align}
which, combined with the above, yields
\begin{multline*}
\ee_j\cdot(\Aa-\bar\Aa)\nabla\Delta_b\Delta_{b'}\phi_i
=-\nabla_k\big(\sigma_{jkl}^*(\cdot-\ee_k)\nabla_l\Delta_b\Delta_{b'}\phi_i\big)
-\nabla_k^*\big(\phi_j^*(\cdot+\ee_k)\ee_k\cdot\Aa\nabla\Delta_b\Delta_{b'}\phi_i\big)\\
-\phi_j^*\nabla^*\cdot\Delta_{b'}\Aa\nabla\Delta_b\phi_i^{b'}-\phi_j^*\nabla^*\cdot\Delta_b\Delta_{b'}\Aa(\nabla\phi_i^{b,b'}+\ee_i)-\phi_j^*\nabla^*\cdot\Delta_b\Aa\nabla\Delta_{b'}\phi_i^b.
\end{multline*}
Inserting this into~\eqref{eq:pre-form-der3-I0}, integrating by parts, and applying the discrete Leibniz rule~\eqref{pr:cov-str-2.4-bis}, we obtain the following representation formula
\begin{multline}\label{eq:vert-der-Xi-applied-second}
\Delta_b\Delta_{b'}J_0^\e(\F)=\e^{\frac d2}\int_{\R^d}\F_{\e,ij}(\nabla\phi_j^*+\ee_j)\cdot\Delta_b\Delta_{b'}\Aa(\nabla\phi^{b,b'}_i+\ee_i)\\
+\e^{\frac d2}\int_{\R^d}\F_{\e,ij}(\nabla\phi_j^*+\ee_j)\cdot\big(\Delta_b\Aa\nabla\Delta_{b'}\phi_i^b+\Delta_{b'}\Aa\nabla\Delta_{b}\phi_i^{b'}\big)\\
+\e^{\frac d2}\int_{\R^d}\big(\sigma_{jkl}^*(\cdot-\ee_k)\nabla_k^*\F_{\e,ij}+\phi_j^*(\cdot+\ee_k)\Aa_{kl}\nabla_k\F_{\e,ij}\big)\nabla_l\Delta_b\Delta_{b'}\phi_i\\
+\e^{\frac d2}\int_{\R^d}\phi_j^*(\cdot+\ee_k)\nabla_k\F_{\e,ij}\ee_k\cdot\big(\Delta_b\Aa\nabla\Delta_{b'}\phi_i^b+\Delta_{b'}\Aa\nabla\Delta_b\phi_i^{b'}\big)\\
+\e^{\frac d2}\int_{\R^d}\phi_j^*(\cdot+\ee_k)\nabla_k\F_{\e,ij}\ee_k\cdot\Delta_b\Delta_{b'}\Aa(\nabla\phi_i^{b,b'}+\ee_i).
\end{multline}
We need to estimate the moment of each right-hand side term. Fix momentarily $b=(z_b,z_b+\ee_b)$ and $b'=(z_{b'},z_{b'}+\ee_{b'})$.
Applying Lemmas~\ref{lem:bd-sigma} and~\ref{lem:ann-Green} to the Green representation formula~\eqref{eq:formvertdernablaphi} for $\nabla\Delta_b\phi$, we find for all $q<\infty$,
\[\expec{|\nabla\Delta_b\phi(x)|^q}^\frac1{q}\lesssim_q(1+|x-z_b|)^{-d}.\]
We then turn to the second vertical derivatives. We obviously have $\Delta_b\Delta_{b'}\Aa=\mathds1_{b=b'}\Delta_b\Aa$.
Next, the Green representation formula applied to equation~\eqref{eq:corr-eq-der2} yields
\begin{multline*}
\nabla\Delta_b\Delta_{b'}\phi_j(x)=-\nabla\nabla G(x,z_{b'})\cdot\Delta_{b'}\Aa(z_{b'})\nabla\Delta_b\phi_j^{b'}(z_{b'})
-\nabla\nabla G(x,z_b)\cdot\Delta_b\Aa(z_b)\nabla\Delta_{b'}\phi_j^b(z_b)\\
-\mathds1_{b=b'}\nabla\nabla G(x,z_b)\cdot\Delta_b\Aa(z_b)(\nabla\phi_j^{b}+\ee_j),
\end{multline*}
so that, for all $q<\infty$, Lemmas~\ref{lem:bd-sigma} and~\ref{lem:ann-Green} lead to
\begin{align*}
\expec{|\nabla \Delta_b\Delta_{b'} \phi(x)|^q}^\frac1q\,\lesssim_q\,(1+|z_b-z_{b'}|)^{-d}\big((1+|x-z_{b'}|)^{-d}+(1+|x-z_b|)^{-d}\big).
\end{align*}
Inserting these estimates into~\eqref{eq:vert-der-Xi-applied-second}, we obtain
\begin{multline*}
\expec{|\Delta_b\Delta_{b'}J_0^\e(\F)|^4}^\frac14
\lesssim\frac{\e^{\frac d2}}{(1+|z_b-z_{b'}|)^{d}}\bigg(\int_{Q(z_b)}|F_\e|+\int_{Q(z_{b'})}|F_\e|\\
+\int_{\R^d}\frac{\mu_d(|x|)^\frac12|\nabla F_{\e}(x)|}{(1+|x-z_{b'}|)^{d}}dx
+\int_{\R^d}\frac{\mu_d(|x|)^\frac12|\nabla F_{\e}(x)|}{(1+|x-z_{b}|)^{d}}dx\bigg).
\end{multline*}
Combining this with~\eqref{eq:bound-pre-K1eps} and with the definition of $K_2^\e$,
with the short-hand notation
\[I(\zeta)(z):=\int_{\R^d}\frac{|\zeta(x)|}{(1+|x-z|)^d}dx,\]
and $G_\e:=\mu_d(|\cdot|)^\frac12\nabla \F_\e$, we deduce
\begin{multline*}
K_2^\e\lesssim\e^{2d}\int_{\R^d}\Big(|\F_\e|^2|I(\F_\e)|^2+|I(|\F_\e|^2)|^2+|I(\F_\e)|^2|I(G_\e)|^2+|I(F_\e I(G_\e))|^2\\
+|F_\e|^2|I(I(G_\e))|^2+|I(|I(G_\e)|^2)|^2+|I(G_\e)|^2|I(I(G_\e))|^2\Big).
\end{multline*}
By the Cauchy-Schwarz inequality and  a multiple use of~\eqref{eq:auxiliary-result-CLT} in the form
\[\|\log^r(2+|\cdot|)\,I(\zeta)\|_{\Ld^p(\R^d)}\lesssim_{p,r}\|\log^{r+1}(2+|\cdot|)\,\zeta\|_{\Ld^p(\R^d)},\]
we are led to
\begin{multline*}
K_2^\e\lesssim\e^{2d}\Big(\|\F_\e\|_{\Ld^4(\R^d)}^2\|\log(2+|\cdot|)\F_\e\|_{\Ld^4(\R^d)}^2+\|\log^\frac12(2+|\cdot|)\F_\e\|_{\Ld^4(\R^d)}^4\\
+\|\log(2+|\cdot|)\F_\e\|_{\Ld^4(\R^d)}^2\|\log(2+|\cdot|)G_\e\|_{\Ld^4(\R^d)}^2+\|F_\e\|_{\Ld^4(\R^d)}^2\|\log^2(2+|\cdot|)G_\e\|_{\Ld^4(\R^d)}^2\\
+\|\log^\frac32(2+|\cdot|)G_\e\|_{\Ld^4(\R^d)}^4+\|\log(2+|\cdot|)G_\e\|_{\Ld^4(\R^d)}^2\|\log^2(2+|\cdot|)G_\e\|_{\Ld^4(\R^d)}^2\Big).
\end{multline*}
Inserting the definition of $G_\e$, rescaling the integrals, and using~\eqref{eq:nabla-D}, the conclusion follows.

\section{Covariance structure}\label{sec:cov-struc}

In this section, we turn to the limiting covariance structure of the homogenization commutator, as stated in Proposition~\ref{prop:conv-cov}.

\subsection{Structure of the proof and auxiliary results}

The main tool to prove Proposition~\ref{prop:conv-cov} is the following stronger version of the spectral gap estimate of Lemma~\ref{lem:var}, which gives an identity (rather than a bound) for the variance of a random variable in terms of its variations. This is an iid\@ version of the so-called Helffer-Sj\"ostrand representation formula~\cite{HS-94,Sjostrand-96} (see also~\cite{NS2,MO}), which holds for any product measure $\Pm$ on~$\Omega$. A proof is included for completeness in Subsection~\ref{sec:proof-HS} below.
{\color{black}Instead of using the difference operator $\triangle_b$ as in Lemma~\ref{lem:var}, this result is more conveniently formulated in terms of $\tilde\Delta_bX:=X-\E_{a(b)}[X]$, with the notation $\E_{a(b)}[\cdot]:=\expeCm{\cdot}{\!(a(b'))_{b'\ne b}}$. Note that $\tilde\Delta_bX=\expeCm{\Delta_bX}{\!\Aa}$ and $\expecm{|\tilde\Delta_bX|^2}=\frac12\expecm{|\Delta_bX|^2}$.}

\begin{lem}\label{lem:cov-inequ-improved}
Let $\Pm$ be a product measure. For all $X=X(\Aa)\in\Ld^2(\Omega)$ we have
\[\var{X}=\sum_{b\in\B}\expec{ (\tilde\Delta_bX)\,\calT\,(\tilde\Delta_bX)},\]
where $\calT:=(\sum_{b\in\B}\tilde\Delta_b\tilde\Delta_b)^{-1}$ is a self-adjoint positive operator on $\Ld^2(\Omega)/\R:=\{X\in\Ld^2(\Omega):\expec{X}=0\}$ with operator norm bounded by $1$.
In particular, it implies the following covariance inequality: for all $X,Y\in \Ld^2(\Omega)$ we have
\[\cov{X}{Y}\,\le \, \frac12\sum_{b\in\B}\expec{|\Delta_b X|^2}^{\frac{1}{2}}
\expec{|\Delta_b Y|^2}^{\frac{1}{2}}.\qedhere\]
\end{lem}

The proof of Proposition~\ref{prop:conv-cov}(i) below further implies that the effective fluctuation tensor~$\calQ$ is given by the following formula, with the notation $b_n:=(0,\ee_n)$,
\begin{eqnarray}\label{eq:representation-K-tocalQ-state}
\calQ_{ijkl}&:=&\sum_{n=1}^d\expec{\big(M_{ij}^n-\expecm{M_{ij}^n}\big)\,\calT\,\big(M_{kl}^n-\expecm{M_{kl}^n}\big)},\\
M_{ij}^n&:=&\E\Big[(a(b_n)-a^{b_n}(b_n))\big(\ee_n\cdot(\nabla\phi_j^*(0)+\ee_j)\big)\big(\ee_n\cdot(\nabla\phi_i^{b_n}(0)+\ee_i)\big)\,\Big|\,\Aa\Big],\nonumber
\end{eqnarray}
in terms of the abstract operator $\calT$ defined above.
Although not convenient for numerical approximation of $\calQ$, this formula allows to easily deduce the non-degeneracy result contained in Proposition~\ref{prop:conv-cov}(ii). In addition, this is key to the proof of Theorem~\ref{th:main-2} on the RVE method.

\subsection{Proof of Lemma~\ref{lem:cov-inequ-improved}}\label{sec:proof-HS}

We start with some observations on the difference operator $\tilde\Delta_b$ on $\Ld^2(\Omega)$. For all $X,Y\in\Ld^2(\Omega)$, by exchangeability of $(\Aa,\Aa^b)$, we find
\begin{multline*}
\expecm{X\tilde\Delta_bY}=\expecm{XY}-\expecm{X\E_{a(b)}[Y]}=\expecm{XY}-\expecm{\E_{a(b)}[X]\E_{a(b)}[Y]}
\\
=\expecm{XY}-\expecm{Y\E_{a(b)}[X]}=\expecm{Y\tilde\Delta_b X},
\end{multline*}
so that $\tilde\Delta_b$ is symmetric on $\Ld^2(\Omega)$. In addition, we easily compute, for all $b,b'\in\B$,
\begin{align}\label{eq:commut-deltae}
\tilde\Delta_b\tilde\Delta_b=\tilde\Delta_b,\qquad\tilde\Delta_b\tilde\Delta_{b'}=\tilde\Delta_{b'}\tilde\Delta_b.
\end{align}
With these observations at hand, we now turn to the study of the (densely defined) operator $\calS:=\sum_{b\in\B}\tilde\Delta_b\tilde\Delta_b$ on $\Ld^2(\Omega)$.
More precisely, we consider the space $\Ld^2(\Omega)/\R:=\{X\in\Ld^2(\Omega):\expec{X}=0\}$ of mean-zero square-integrable random variables, and we show that $\calS$ is an essentially self-adjoint, non-negative operator on $\Ld^2(\Omega)/\R$ with dense image.
First, since $\expecm{\tilde\Delta_bX}=0$ for all $b\in\B$ and $X\in\Ld^2(\Omega)$, the image $\text{Im}\,\calS$ is clearly contained in $\Ld^2(\Omega)/\R$.
Second, for all $X\in\Ld^2(\Omega)$ in the domain of $\calS$, we compute
\[\expec{X\calS X}=\sum_{b\in\B}\expecm{|\tilde\Delta_bX|^2}\ge0,\]
which shows that $\calS$ is non-negative. Third, if $X\in\Ld^2(\Omega)/\R$ in the domain of $\calS$ is orthogonal to the image $\text{Im}\,\calS$, then we deduce
\[0=\expec{X\calS X}=\sum_{b\in\B}\expecm{|\tilde\Delta_bX|^2},\]
so that $\tilde\Delta_bX=0$ almost surely for all $b\in\B$, which implies that $X$ is constant.

\medskip\noindent
These properties of $\calS$ allow us to define (densely) the inverse $\calT:=\calS^{-1}$ as an essentially self-adjoint, non-negative operator on $\Ld^2(\Omega)/\R$.
Finally, the spectral gap of Lemma~\ref{lem:var} implies, for all $X\in\Ld^2(\Omega)/\R$ in the domain of $\calS$,
\[\|X\|_{\Ld^2(\Omega)}^2=\var{X}\le\sum_{b\in\B}\expecm{|\tilde\Delta_bX|^2}=\expec{X\calS X}\le\|X\|_{\Ld^2(\Omega)}\|\calS X\|_{\Ld^2(\Omega)},\]
and hence $\|X\|_{\Ld^2(\Omega)}\le \|\calS X\|_{\Ld^2(\Omega)}$, which implies that $\calT=\calS^{-1}$ on $\Ld^2(\Omega)/\R$ has operator norm bounded by $1$.

\medskip\noindent
It remains to establish the representation formula for the variance. By density, it suffices to prove it for all $X\in\text{Im}\,\calS$. Writing $X=\calS Y$ for some $Y\in\Ld^2(\Omega)/\R$, we decompose
\begin{align*}
\var{X}=\expec{X\calS Y}=\sum_{b\in\B}\expecm{\tilde\Delta_bX\tilde\Delta_bY}=\sum_{b\in\B}\expecm{(\tilde\Delta_bX)(\tilde\Delta_b\calT X)}.
\end{align*}
Since the commutation relations~\eqref{eq:commut-deltae} ensure that $\tilde\Delta_b\calS=\calS\tilde\Delta_b$ holds on the domain of $\calS$ in $\Ld^2(\Omega)$, we deduce $\tilde\Delta_b\calT=\calT\tilde\Delta_b$ on $\Ld^2(\Omega)/\R$, and the above then leads to the desired representation
\begin{align*}
\var{X}=\sum_{b\in\B}\expecm{(\tilde\Delta_bX)\calT(\tilde\Delta_bX)}.
\end{align*}

\subsection{Proof of Proposition~\ref{prop:conv-cov}(i)}

By polarization and linearity,
it is enough to prove~\eqref{eq:conv-cov-main-I0} with $\F=\G\in C^\infty_c(\R^d)^{d\times d}$.
We thus need to establish the convergence of the variance
\[\nu_\e:=\var{\e^{-\frac d2}\int_{\R^d}\F:\Xi(\tfrac \cdot\e)}=\var{\e^{\frac d2}\int_{\R^d}\F_\e:\Xi},\]
where we have set $\F_\e:=\F(\e\cdot)$.
We split the proof into two steps.

\medskip
\step1 Proof of~\eqref{eq:conv-cov-main-I0}.\\
The Helffer-Sj\"ostrand representation of Lemma~\ref{lem:cov-inequ-improved} applied to the variance $\nu_\e$ yields
\begin{align}\label{eq:HS-first-appl}
\nu_\e&=\e^d\sum_{b\in\B}\expec{\Big(\tilde\Delta_b\int_{\R^d}\F_\e: \Xi\Big)\,\calT\,\Big(\tilde\Delta_b\int_{\R^d}\F_\e:\Xi \Big)}.
\end{align}
We now appeal to~\eqref{eq:der-Xi-claim} in the form
\begin{multline*}
\Delta_b\int_{\R^d}\F_\e: \Xi=\int_{\R^d}\F_{\e,ij}(\nabla\phi_j^*+\ee_j)\cdot\Delta_b\Aa(\nabla\phi_i^b+\ee_i)\\
+\int_{\R^d}\phi_j^*(\cdot+\ee_k)\nabla_k\F_{\e,ij}\ee_k\cdot\Delta_b\Aa(\nabla\phi_i^b+\ee_i)
+\int_{\R^d}\nabla h_{\e,i}\cdot\Delta_b\Aa(\nabla\phi_i^b+\ee_i),
\end{multline*}
where the auxiliary field $h_{\e,i}$ is the unique Lax-Milgram solution in $\R^d$ of
\begin{align}\label{eq:def-reps-reutil}
-\nabla^*\cdot\Aa^*\nabla h_{\e,i}=\nabla_l^*\big(\phi_j^*(\cdot+\ee_k)\Aa_{kl}\nabla_k\F_{\e,ij}+\sigma_{jkl}^*(\cdot-\ee_k)\nabla_k^*\F_{\e,ij}\big).
\end{align}
Recalling that $\tilde\Delta_bX=\E_{a^b(b)}[\Delta_bX]$, inserting this representation formula into~\eqref{eq:HS-first-appl}, extracting the first term $U_\e$ defined below, and using that $\calT$ on $\Ld^2(\Omega)/\R$ has operator norm bounded by $1$, we find
\begin{align}\label{eq:decomp-rep-diff-nu12+cov}
|\nu_\e-\e^dU_\e|\le {\e^d}\sum_{b\in\B}\Big(S_{\e}^b T_{\e}^b+ \frac12(T_{\e}^b)^2\Big),
\end{align}
where for convenience we define
\begin{eqnarray*}
U_\e&:=&\sum_{b\in\B}\expec{\,(V_\e^b-\expecm{V_\e^b})\,\calT\,(V_\e^b-\expecm{V_\e^b})\,},\\
V_\e^b&:=&\E_{a^b(b)}\Big[\int_{\R^d}\F_{\e,ij}(\nabla\phi_j^*+\ee_j)\cdot\Delta_b\Aa(\nabla\phi_i^b+\ee_i)\Big],
\end{eqnarray*}
while for all $b\in\B$ the error terms are given by
\begin{eqnarray*}
S_{\e}^b&:=&\expec{\Big(\int_{\R^d}|\Delta_b\Aa||\nabla\phi^*+\Id||\nabla\phi^b+\Id||\F_\e|\Big)^2}^\frac12,
\end{eqnarray*}
and by $T_{\e}^b:=T_{\e,1}^b+T_{\e,2}^b$ with
\begin{eqnarray*}
T_{\e,1}^b&:=&\sum_{k=1}^d\expec{\Big(\int_{\R^d}|\Delta_b\Aa||\phi^*(\cdot+\ee_k)||\nabla\phi^b+\Id||\nabla\F_{\e}|\Big)^2}^\frac12,\\
T_{\e,2}^b&:=&\expec{\Big(\int_{\R^d}|\Delta_b\Aa||\nabla\phi^b+\Id||\nabla h_{\e}|\Big)^2}^\frac12.
\end{eqnarray*}
We start with the analysis of $U_\e$. Writing $\Delta_b \Aa(x)=(a(b)-a^b(b))\mathds1_{Q(z_b)}(x) \ee_b\otimes \ee_b$ for $b=(z_b,z_b+\ee_b)$, we may compute
\begin{align*}
V_\e^b=\Big(\int_{Q(z_b)}\F_{\e,ij}\Big)\,\E_{a^b(b)}\Big[(a(b)-a^b(b))\big(\ee_b\cdot(\nabla\phi_j^*(z_b)+\ee_j)\big)\big(\ee_b\cdot(\nabla\phi_i^b(z_b)+\ee_i)\big)\Big],
\end{align*}
so that,
by stationarity,
\begin{align*}
\e^d U_\e=\calQ_{ijkl}\,\e^d\sum_{z\in\Z^d}\Big(\int_{Q(z)}\F_{\e,ij}\Big)\Big(\int_{Q(z)}\F_{\e,kl}\Big),
\end{align*}
where the coefficient $\calQ_{ijkl}$ is defined in~\eqref{eq:representation-K-tocalQ-state} above.
Since $\calT$ on $\Ld^2(\Omega)/\R$ has operator norm bounded by $1$, the moment bounds of Lemma~\ref{lem:bd-sigma} yield
\begin{align*}
|\calQ_{ijkl}| \,\lesssim\, \sum_{n=1}^d\expec{|\nabla\phi^*+\Id|^2|\nabla\phi^{b_n}+\Id|^2}\,\lesssim\,1.
\end{align*}
We may then estimate the discretization error
\begin{eqnarray}
\Big|\e^dU_\e-\calQ_{ijkl}\,\int_{\R^d}\F_{ij}\F_{kl}\Big|&=&\Big|\e^dU_\e-\calQ_{ijkl}\,\e^d\int_{\R^d}\F_{\e,ij}\F_{\e,kl}\Big|\nonumber\\
&\lesssim&\e^d\sum_{z\in\Z^d}\int_{Q(z)}\Big|\F_\e(x)-\int_{Q(z)}\F_\e\Big|^2dx\nonumber\\
&\lesssim&\e^d\int_{\R^d}|D\F_\e|^2=\e^2\int_{\R^d}|D\F|^2.\label{eq:conv-main-term-Ueps}
\end{eqnarray}
We now turn to the estimate of the right-hand side of~\eqref{eq:decomp-rep-diff-nu12+cov}.
Using $|\Delta_b\Aa(x)|\lesssim\mathds1_{Q(z_b)}(x)$ and the moment bounds of Lemma~\ref{lem:bd-sigma}, we obtain
\begin{align*}
S_{\e}^b\,\lesssim\, \expec{|\nabla\phi^*+\Id|^2|\nabla\phi^b+\Id|^2}^\frac12\int_{Q(z_b)}|\F_\e|\,\lesssim\,\int_{Q(z_b)}|\F_\e|.
\end{align*}
Hence, by the Cauchy-Schwarz inequality,
\begin{align}
\sum_{b\in\B}S_\e^bT_{\e}^b\,\lesssim\,\sum_{b\in\B} T_\e^b \int_{Q(z_b)}|\F_\e|\,&\lesssim\,\|\F_\e\|_{\Ld^2(\R^d)}\Big(\sum_{b\in\B}(T_\e^b)^2\Big)^\frac12\nonumber\\
&\lesssim\,\e^{-\frac d2}\|\F\|_{\Ld^2(\R^d)}\Big(\sum_{b\in\B}(T_\e^b)^2\Big)^\frac12,\label{eq:bound-Seps-e}
\end{align}
and it remains to estimate
$$\sum_{b\in\B}(T_{\e}^b)^2\le 2\sum_{b\in\B}(T_{\e,1}^b)^2+2\sum_{b\in\B}(T_{\e,2}^b)^2.$$
First, using $|\Delta_b\Aa(x)|\lesssim\mathds1_{Q(z_b)}(x)$ and the moment bounds of Lemma~\ref{lem:bd-sigma}, we find
\begin{align}
\e^d\sum_{b\in\B}(T_{\e,1}^b)^2\,\lesssim_{\alpha,p}\,\e^d\|\mu_d(|\cdot|)^\frac12\nabla\F_\e\|_{\Ld^2(\R^d)}^2.\label{eq:estim-sumT1}
\end{align}
Second, arguing as in the proof of Proposition~\ref{prop:bounded-I0} (cf.~\eqref{eq:main-passage-first-lem}), using the large-scale weighted Calder\'on-Zygmund theory (cf.\@ Lemma~\ref{lem:max-reg}) applied to equation~\eqref{eq:def-reps-reutil} for $h_\e$, we obtain for all $0<p-1\ll1$ and all $\alpha>d$,
\begin{align}
\e^d\sum_{b\in\B}(T_{\e,2}^b)^2\,\lesssim_{\alpha,p}\,\e^{\frac{d}p}\|w_\e^{\alpha\frac{p-1}{2p}}\mu_d(|\cdot|)^\frac12\nabla\F_\e\|_{\Ld^{2p}(\R^d)}^2.\label{eq:estim-sumT2}
\end{align}
Rescaling the integrals and using~\eqref{eq:nabla-D} and H\"older's inequality, we find
\begin{align*}
\e^d\sum_{b\in\B}(T_{\e}^b)^2\,\lesssim_{\alpha,p}\,\e^{2}\mu_d(\tfrac1\e)\|w_1^{\alpha\frac{p-1}{2p}}\mu_d(|\cdot|)^\frac12D\F\|_{\Ld^{2p}(\R^d)}^2,
\end{align*}
and the conclusion~\eqref{eq:conv-cov-main-I0} follows.

\medskip
\step2 Proof of the Green-Kubo formula~\eqref{e.cov-struc.2}.\\
In order to establish~\eqref{e.cov-struc.2}, it suffices to repeat the argument of Step~1 with the test function $\F=\mathds1_{Q}\,\ee_i\otimes\ee_j$ (hence~$\F_\e=\mathds1_{\frac1\e Q}\,\ee_i\otimes\ee_j$), for some fixed $1\le i,j\le d$. Lemma~\ref{lem:cov-inequ-improved} again leads to~\eqref{eq:decomp-rep-diff-nu12+cov}, and we briefly indicate how to analyze the different terms in the present setting. First, the estimate~\eqref{eq:conv-main-term-Ueps} is replaced by the following (no summation over repeated indices),
\begin{eqnarray*}
|\e^dU_\e-\calQ_{ijij}|&\lesssim&\e^d\sum_{z\in\Z^d}\int_Q\Big(\mathds1_{\frac1\e Q}(z+x)-\int_Q\mathds1_{\frac1\e Q}(z+y)dy\Big)^2dx\\
&\le& \e^d\sum_{z\in\Z^d}\mathds1_{(z+Q)\cap\partial(\frac1\e Q)\ne\varnothing}
~\lesssim~\e.
\end{eqnarray*}
Second, the estimate~\eqref{eq:bound-Seps-e} remains unchanged. Third, using estimates~\eqref{eq:estim-sumT1} and~\eqref{eq:estim-sumT2},
and noting that $|\nabla\F_\e|\lesssim\mathds1_{A_\e}$ with $A_\e:=B+\partial Q_{\frac1\e}$ and that $w_\e\lesssim1$ and $\mu_d(|\cdot|)\lesssim\mu_d(\frac1\e)$ on $A_\e$, we deduce
\begin{eqnarray*}
\sum_{b\in\B}(T_{\e,1}^b)^2&\lesssim&\mu_d(\tfrac1\e)|A_\e|\,\lesssim\,\e^{1-d}\mu_d(\tfrac1\e),\\
\sum_{b\in\B}(T_{\e,2}^b)^2&\lesssim&\e^{-d\frac{p-1}p}\mu_d(\tfrac1\e)|A_\e|^\frac1p\,\lesssim\,\e^{\frac{1}p-d}\mu_d(\tfrac1\e),
\end{eqnarray*}
and the conclusion~\eqref{e.cov-struc.2} follows.

%%%%%%%%%%%%%%%%%%%%%%
%%%%%%%%%%%%%%%%%%%%%%

\subsection{Proof of Proposition~\ref{prop:conv-cov}(ii)}

The following proof of the non-degeneracy of $\calQ$ is based on the Helffer-Sj\"ostrand representation formula (see also~\cite[Remark~2.3]{MO}),
and constitutes a shorter alternative to the corresponding proof in \cite{GN}. Given a fixed direction $\ee\in\R^d\setminus\{0\}$, and letting $\phi_\ee$ denote the corrector in this direction, we may write, in view of formula~\eqref{eq:representation-K-tocalQ-state} (with $\phi_e^*=\phi_e$ by symmetry of the coefficients),
\begin{gather}
(\ee\otimes\ee):\calQ\,(\ee\otimes \ee) =\sum_{n=1}^d\expecm{(\ee\cdot M^n\ee)\calT (\ee\cdot M^n\ee)},\label{eq:representation-K-tocalQ-state-utilis}\\
\ee\cdot M^n\ee:=\E\Big[(a(b_n)-a^{b_n}(b_n))\big(\ee_n\cdot(\nabla\phi_\ee(0)+\ee)\big)\big(\ee_n\cdot(\nabla\phi_\ee^{b_n}(0)+\ee)\big)\,\Big|\,\Aa\Big],\nonumber
\end{gather}
since the exchangeability of $(\Aa,\Aa^{b_n})$ indeed yields $\expecm{\ee\cdot M^n\ee}=0$ for all $n$.
We start with a suitable reformulation of $\ee\cdot M^n\ee$. Considering the difference of the corrector equation~\eqref{e.corr} for $\phi_\ee$ and $\phi_\ee^{b_n}$ in the form $-\nabla^*\cdot\Aa^{b_n}\nabla(\phi_\ee^{b_n}-\phi_\ee)=\nabla^*\cdot(\Aa^{b_n}-\Aa)(\nabla\phi_\ee+\ee)$,
an integration by parts yields
\begin{multline*}
\int_{\R^d}\nabla(\phi_\ee^{b_n}-\phi_\ee)\cdot\Aa^{b_n}\nabla(\phi_\ee^{b_n}-\phi_\ee)=-\int_{\R^d}\nabla(\phi_\ee^{b_n}-\phi_\ee)\cdot(\Aa^{b_n}-\Aa)(\nabla\phi_\ee+\ee)\\
=(a(b_n)-a^{b_n}(b_n))(\ee_n\cdot\nabla(\phi_\ee^{b_n}-\phi_\ee)(0))(\ee_n\cdot(\nabla\phi_\ee(0)+\ee)).
\end{multline*}
Hence, by definition of $\ee\cdot M^n\ee$,
\begin{multline}\label{eq:bound-gen-Lek}
\ee\cdot M^n\ee=\E\bigg[\int_{\R^d}\nabla(\phi^{b_n}_\ee-\phi_\ee)\cdot\Aa^{b_n}\nabla(\phi_\ee^{b_n}-\phi_\ee)\,\bigg|\,\Aa\bigg]\\
+(a(b_n)-\expec{a(b_n)})(\ee_n\cdot(\nabla\phi_\ee(0)+\ee))^2.
\end{multline}
We now argue by contradiction. If $(\ee\otimes\ee):\calQ\,(\ee\otimes \ee)=0$, then by formula~\eqref{eq:representation-K-tocalQ-state-utilis} and by the non-negativity of~$\calT$ we would have $\expecm{(\ee\cdot M^n\ee)\calT (\ee\cdot M^n\ee)}=0$ for all $n$. Let $1\le n\le d$ be momentarily fixed.
Recalling that $\calT=\calS^{-1}$ with $\calS=\sum_{b\in\B}\tilde\Delta_b\tilde\Delta_b$, this would imply
\[0=\expec{(\calT (\ee\cdot M^n\ee))\calS(\calT (\ee\cdot M^n\ee))}=\sum_{b\in\B}\expec{\big|\tilde\Delta_b\calT (\ee\cdot M^n\ee)\big|^2},\]
hence $\calT (\ee\cdot M^n\ee)=0$, and thus $\ee\cdot M^n\ee=0$ almost surely. Formula~\eqref{eq:bound-gen-Lek} would then imply
\begin{multline}\label{eq:pre-contrad-Qlim}
(a(b_n)-\expec{a(b_n)})(\ee_n\cdot(\nabla\phi_\ee(0)+\ee))^2\\
=-\E\bigg[\int_{\R^d}\nabla(\phi^{b_n}_\ee-\phi_\ee)\cdot\Aa^{b_n}\nabla(\phi_\ee^{b_n}-\phi_\ee)\,\bigg|\,\Aa\bigg],
\end{multline}
almost surely.
Since the law of $a(b_n)$ is non-degenerate, the event $a(b_n)>\expec{a(b_n)}$ occurs with a positive probability. Conditioning on this event, the left-hand side in~\eqref{eq:pre-contrad-Qlim} is non-negative, and the non-positivity of the right-hand side would then imply that both sides vanish, that is,
\[\ee_n\cdot(\nabla\phi_\ee(0)+\ee)=0\qquad\text{and}\qquad\E\bigg[\int_{\R^d}\nabla(\phi^{b_n}_\ee-\phi_\ee)\cdot\Aa^{b_n}\nabla(\phi_\ee^{b_n}-\phi_\ee)\,\bigg|\,\Aa\bigg]=0,\]
almost surely.
Since the integrand in this last expectation is non-negative, we would deduce that the event $a(b_n)>\expec{a(b_n)}$ entails $\ee_n\cdot(\nabla\phi_\ee(0)+\ee)=0$ and $\nabla\phi_\ee(0)=\nabla\phi_\ee^{b_n}(0)$, and thus also $\ee_n\cdot(\nabla\phi_\ee^{b_n}(0)+\ee)=0$ almost surely. Since this last event is independent of $a(b_n)$, hence of the conditioning event, we would deduce unconditionally $\ee_n\cdot(\nabla\phi_\ee^{b_n}(0)+\ee)=0$ almost surely. By exchangeability of $(\Aa,\Aa^{b_n})$, this means $\ee_n\cdot(\nabla\phi_\ee(0)+\ee)=0$ almost surely.  As this holds for any $n$, we would conclude $\nabla\phi_\ee(0)+\ee=0$ almost surely, and taking the expectation would lead to a contradiction.

%%%%%%%%%%%%%%%%%%%%%%%%%
%%%%%%%%%%%%%%%%%%%%%%%%%
%%%%%%%%%%%%%%%%%%%%%%%%%

\section{Approximation of the fluctuation tensor}\label{app:RVE-Q}

In this section, we analyze the RVE method for the approximation of the fluctuation tensor $\calQ$ as stated in Theorem~\ref{th:main-2}.

\subsection{Structure of the proof and auxiliary results}

The estimate on the standard deviation is obtained similarly as the CLT scaling in Proposition~\ref{prop:bounded-I0}, noting that the large-scale Calder\'on-Zygmund result of Lemma~\ref{lem:max-reg} also holds for the periodized operator $-\nabla^*\cdot\Aa_L\nabla$ on $Q_L$.\footnote{The only issue concerns the moment bounds for the corresponding minimal radius $r_{*,L}$ associated with the periodized operator, $\expec{r_{*,L}^q}\lesssim_q1$ for all $q<\infty$. By definition of $r_{*,L}$ in~\cite{GNO-reg}, this is a consequence of a sup-bound based on the version of Lemma~\ref{lem:bd-sigma} for the periodized correctors $(\phi_L,\sigma_L)$ (cf.~\cite{GNO1}).}
The characterization~\eqref{eq:calQ-charact} of $\calQ$ and the estimate of the systematic error are deduced as corollaries of formula~\eqref{eq:representation-K-tocalQ-state} for the fluctuation tensor $\calQ$, together with the following crucial estimates on the periodized corrector $\phi_L$.
The first estimate on $\nabla\phi_L$ is stated as such in~\cite[Proposition~1]{GNO1},
and the second one follows from a decomposition of the difference $\nabla\phi_L-\nabla\phi$ via massive approximation of the corrector
and Richardson extrapolation, applying~\cite[Lemma~2.8 and estimate (2.68)]{GN}, and optimizing the mass.

\begin{lem}[\cite{GNO1,GN}]\label{lem:corr-est-L0}
Let $d\ge2$ and let $\Pm$ be a product measure. For all $L\ge2$ and all $q<\infty$ we have
\[\expec{|\nabla\phi_L|^q}^\frac1q\lesssim_q1,\qquad\text{and}\qquad
%and in addition,
\expec{|\nabla(\phi_L-\phi)(0)|^q}^\frac1q\lesssim_qL^{-\frac d2}\log^\frac d2L.\qedhere\]
\end{lem}

\subsection{Proof of Theorem~\ref{th:main-2}}
We split the proof into two steps: we first estimate the variance of the RVE approximation, and then we turn to the characterization~\eqref{eq:calQ-charact} of $\calQ$ and to the systematic error of the RVE approximation.

\medskip
\step1 Proof of the random error estimate $|\var{\calQ_{L,N}}|^\frac12\lesssim N^{-\frac12}$.\\
Since the realizations $\bar\Aa_{L}^{(n)}$ are iid copies of $\bar\Aa_{L}$, the definition~\eqref{eq:def-QL-per} of $\calQ_{L,N}$ leads after straightforward computations to
\begin{align*}
\var{\calQ_{L,N}}=N^{-1}\var{\big(L^\frac d2\bar\Aa_{L}^*-\expecm{L^\frac d2\bar\Aa_{L}^*}\big)^{\otimes2}},
\end{align*}
and hence,
\begin{align*}
|\var{\calQ_{L,N}}\!|\lesssim N^{-1}\,\expec{\big|L^{\frac d2}(\bar\Aa_{L}-\expec{\bar\Aa_{L}})\big|^4}.
\end{align*}
Arguing as in~\cite[Lemma~2]{GNO1}, the spectral gap estimate of Lemma~\ref{lem:var} is seen to imply the following inequality: for all $X=X(\Aa)\in\Ld^4(\Omega)$,
\begin{align*}
\expec{(X-\expec{X})^{4}}
\,\le\, 4\, \expec{\bigg(\sum_{b\in\B}|\Delta_bX|^2\bigg)^2}.
\end{align*}
Applying this inequality to (each component of) $X=\bar\Aa_{L}$, we deduce
\begin{align*}
|\var{\calQ_{L,N}}\!|\lesssim
N^{-1}\,\expec{\bigg(\sum_{b\in\B_L}\Big(L^{-\frac d2}\int_{Q_L}\Delta_b\big(\Aa_L(\nabla\phi_L+\Id)\big)\Big)^2\bigg)^2}.
\end{align*}
Arguing as in the proof of Proposition~\ref{prop:bounded-I0} (with $\e$ replaced by $\frac1L$ and $F_\e$ replaced by $\Id$), using the periodized version of Lemma~\ref{lem:max-reg} and the moment bounds of Lemma~\ref{lem:corr-est-L0}, the conclusion follows.

\medskip
\step2 Proof of~\eqref{eq:calQ-charact} and of the systematic error estimate $|\expec{\calQ_{L,N}}-\calQ|\lesssim L^{-\frac d2}\log^{\frac d2} L$.\\
Since the realizations $\bar\Aa_{L}^{(n)}$ are iid copies of $\bar\Aa_{L}$, the definition~\eqref{eq:def-QL-per} of $\calQ_{L,N}$ yields after straightforward computations $\expec{\calQ_{L,N}}=\varm{L^{\frac d2}\bar\Aa_{L}^*}$, that is,
\begin{eqnarray}
\expec{(\calQ_{L,N})_{ijkl}}&=&\cov{L^{\frac d2}\bar\Aa_{L,ji}}{L^{\frac d2}\bar\Aa_{L,lk}}\nonumber \\
&=&L^{-d}\,\cov{\int_{Q_L}\ee_j\cdot\Aa_L(\nabla\phi_{L,i}+\ee_i)}{\int_{Q_L}\ee_l\cdot\Aa_L(\nabla\phi_{L,k}+\ee_k)}.\label{eq:rewrite-EQL}
\end{eqnarray}
For $b\in\B$, we write $b=(z_b,z_b+\ee_b)$.
Using the periodized corrector equation~\eqref{eq:corr-per-def} and its vertical derivative, and recalling that $\Delta_b\Aa_L(x)=(a(b)-a^b(b))\mathds1_{Q(z_b)}(x)\ee_b\otimes\ee_b$ for $b\in\B_L$ and $x\in Q_L$, we find
\begin{eqnarray}
\lefteqn{\Delta_b\int_{Q_L}\ee_j\cdot\Aa_L(\nabla\phi_{L,i}+\ee_i)=\int_{Q_L}\ee_j\cdot\Delta_b\Aa_L(\nabla\phi_{L,i}^b+\ee_i)+\int_{Q_L}\ee_j\cdot\Aa_L\nabla\Delta_b\phi_{L,i}}\nonumber\\
&\hspace{2cm}=&\int_{Q_L}\ee_j\cdot\Delta_b\Aa_L(\nabla\phi_{L,i}^b+\ee_i)-\int_{Q_L}\nabla\phi_{L,j}^*\cdot\Aa_L\nabla\Delta_b\phi_{L,i}\nonumber\\
&\hspace{2cm}=&\int_{Q_L}(\nabla\phi_{L,j}^*+\ee_j)\cdot\Delta_b\Aa_L(\nabla\phi_{L,i}^b+\ee_i)\nonumber\\
&\hspace{2cm}=&(a(b)-a^{b}(b))(\ee_b\cdot(\nabla\phi_{L,j}^*(z_b)+\ee_j))(\ee_b\cdot(\nabla\phi_{L,i}^b(z_b)+\ee_i)).\label{eq:form-De-preQLN}
\end{eqnarray}
Applying the Helffer-Sj\"ostrand representation formula of Lemma~\ref{lem:cov-inequ-improved} to the covariance in~\eqref{eq:rewrite-EQL}, we obtain by stationarity, as in the proof of Proposition~\ref{prop:conv-cov}(i),
\begin{align*}
\expec{(\calQ_{L,N})_{ijkl}}\,=\,\sum_{n=1}^d\expec{M_{ij,L}^n\,\calT\, M_{kl,L}^n},
\end{align*}
where we have set
\[M_{ij,L}^n\,:=\,\E\Big[(a(b_n)-a^{b_n}(b_n))(\ee_n\cdot(\nabla\phi_{L,j}^*(0)+\ee_j))(\ee_n\cdot(\nabla\phi_{L,i}^{b_n}(0)+\ee_i))\,\Big|\,\Aa\Big].\]
Noting that~\eqref{eq:form-De-preQLN} implies $\expecm{M_{ij,L}^n}=0$, comparing the above identity for $\expec{(\calQ_{L,N})_{ijkl}}$ with formula~\eqref{eq:representation-K-tocalQ-state} for $\calQ$, and using that the operator $\calT$ on $\Ld^2(\Omega)/\R$ has operator norm bounded by $1$, we deduce
\begin{align*}
\big|\expec{(\calQ_{L,N})_{ijkl}}-\calQ_{ijkl}\big|\,\lesssim\,\expec{|\nabla(\phi_L-\phi)(0)|^4}^\frac14\big(\,\expec{|\nabla\phi_L|^4}+\expec{|\nabla\phi|^4}\big)^\frac34,
\end{align*}
and the conclusion follows from Lemmas~\ref{lem:bd-sigma} and~\ref{lem:corr-est-L0}.

%%%%%%%%%%%%%

\section*{Acknowledgements}
The work of MD is supported by F.R.S.-FNRS through a Research Fellowship and by the CNRS-Momentum program.
AG acknowledge financial support from the European Research Council under
the European Community's Seventh Framework Programme (FP7/2014-2019 Grant Agreement
QUANTHOM 335410). The authors acknowledge the hospitality of IH\'ES, where this work was initiated in February 2015, and the support
of the Chaire Schlumberger.

%%%%%%%%%%%%%

%\input{bib}
\bibliographystyle{plain}
%\bibliography{biblio}

\def\cprime{$'$} \def\cprime{$'$} \def\cprime{$'$}

%%%%%%%%%%%%

\end{document}